\documentclass[12pt]{amsart}

\usepackage{bbm} % for identity functors
\usepackage{enumitem}
\usepackage[mode=buildnew]{standalone}% requires -shell-escape
\usepackage{tikz}
\usepackage{mathrsfs}
\usepackage{color}
\usepackage[all]{xy}
\usepackage{amsmath, amssymb}
\usepackage{graphics}
\usepackage{graphicx}
\usepackage{mathtools}
\usepackage[pdftex,hypertexnames=false]{hyperref}
\hypersetup{
	colorlinks=true,
	citecolor=blue,
	linkcolor=blue,
	urlcolor=blue
	}
	\usepackage{cleveref}
	\usepackage[square,numbers]{natbib}
	\usepackage{autonum}
\usepackage{geometry}
\numberwithin{equation}{subsection}
\theoremstyle{plain}

\newtheorem{theorem}{Theorem}[section]
\crefname{theorem}{Theorem}{Theorems}

\crefname{mainresult}{Main Result}{Main Results}
\newtheorem{lemma}[theorem]{Lemma}
\crefname{lemma}{Lemma}{Lemmas}
\newtheorem{proposition}[theorem]{Proposition}
\crefname{proposition}{Proposition}{Propositions}
\newtheorem{proposition-definition}[theorem]{Proposition-Definition}
\crefname{proposition-definition}{Proposition-Definition}{}
\newtheorem{corollary}[theorem]{Corollary}
\crefname{corollary}{Corollary}{Corollaries}
\newtheorem{definition}[theorem]{Definition}
\crefname{definition}{Definition}{Definitions}
\newtheorem{example}[theorem]{Example}
\crefname{example}{Example}{Examples}
\newtheorem{remark}[theorem]{Remark}
\crefname{remark}{Remark}{Remarks}

\crefname{criterion}{Criterion}{Criterions}

\crefname{observation}{Observation}{Observations}
\newtheorem{question-aim}[theorem]{\textcolor{magenta}{Question-Aim}}

\crefname{question}{Question}{Questions}

\crefname{problem}{Problem}{Problems}

\crefname{assumption}{Assumption}{Assumptions}
\newtheorem{notation}{Notation}
\crefname{notation}{Notation}{Notations}

\newcommand{\End}{\operatorname{End}}

\newcommand{\sEndZ}{\underline{\mathrm{End}}^{\mathbb{Z}}}
\newcommand{\Hom}{\operatorname{Hom}}
\newcommand{\HomZ}{\operatorname{Hom}^{\mathbb{Z}}}
\newcommand{\sHom}{\underline{\mathrm{Hom}}}
\newcommand{\sHomZ}{\underline{\mathrm{Hom}}^{\mathbb{Z}}}
\renewcommand{\Im}{\operatorname{Im}}
\newcommand{\Ker}{\operatorname{Ker}}

\newcommand{\id}{\operatorname{id}}

\newcommand{\Aut}{\operatorname{Aut}\nolimits}

\renewcommand{\mod}{\operatorname{mod}}
\newcommand{\proj}{\operatorname{proj}}
\newcommand{\Gproj}{\operatorname{Gproj}}
\newcommand{\sGproj}{\underline{\operatorname{Gproj}}\,}
\newcommand{\ind}{\operatorname{ind}}
\newcommand{\thick}{\operatorname{thick}}
\newcommand{\projZ}{\operatorname{proj}^{\mathbb{Z}}}
\newcommand{\GprojZ}{\operatorname{Gproj}^{\mathbb{Z}}\nolimits}
\newcommand{\sGprojZ}{\underline{\operatorname{Gproj}}^{\mathbb{Z}}}
\newcommand{\modZ}{\operatorname{mod}^{\mathbb{Z}}}
\newcommand{\smod}{\underline{\operatorname{mod}}\,}
\newcommand{\smodZ}{\underline{\operatorname{mod}}^{\mathbb{Z}}}

\newcommand{\injdim}{\operatorname{inj.dim}}

\renewcommand{\L}{\Lambda}

\newcommand{\G}{\Gamma}

\newcommand{\F}{\operatorname{\mathbb{F}}\nolimits}

\renewcommand{\P}{\mathcal{P}}
\newcommand{\B}{\mathcal{B}}
\newcommand{\C}{\mathcal{C}}
\newcommand{\T}{\mathcal{T}}

\newcommand{\D}{\operatorname{\mathsf{D}}}
\newcommand{\K}{\operatorname{\mathsf{K}}}

\newcommand{\Z}{\mathbb{Z}}

\newcommand{\Ac}{\mathbb{A}_c}
\newcommand{\A}{\mathbb{A}}

\newcommand{\degQ}{\deg_{KQ}}

\begin{document}

%%%%%%%%%%%%%%%%%%%%%%%%%%%%%%%%%%%%%%%%%%%%%
%    Title        
%%%%%%%%%%%%%%%%%%%%%%%%%%%%%%%%%%%%%%%%%%%%%
\title[The stable category of Gorenstein-projective modules]{The stable category of Gorenstein-projective modules over a monomial algebra}
%%%%%%%%%%%%%%%%%%%%%%%%%%%%%%%%%%%%%%%%%%%%%
%    Title        
%%%%%%%%%%%%%%%%%%%%%%%%%%%%%%%%%%%%%%%%%%%%%

%%%%%%%%%%%%%%%%%%%%%%%%%%%%%%%%%%%%
%       author information
%%%%%%%%%%%%%%%%%%%%%%%%%%%%%%%%%%%%
\author[T. Honma]{Takahiro Honma}
\address{%
{\rm T. Honma:} Department of General Education \\
National Institute of Technology (KOSEN), Yuge College \\ 
Ochi \\
Ehime 794-2593\\ 
Japan}
\email{99cfqqc9@gmail.com}
%%%%%%%%%%%%%%%%%%%%%%%%%%%%%%%%%%%%
\author[S. Usui]{Satoshi Usui}
\address{%
{\rm S. Usui:} Department of General Education \\
Tokyo Metropolitan College of Industrial Technology \\ 
8-17-1, Minamisenju, Arakawa-ku \\
Tokyo 116-8523\\ 
Japan}
\email{usui@metro-cit.ac.jp}
%%%%%%%%%%%%%%%%%%%%%%%%%%%%%%%%%%%%
%       author information
%%%%%%%%%%%%%%%%%%%%%%%%%%%%%%%%%%%%

\thanks{}



\date{}

\dedicatory{}

%%%%%%%%%%%%%%%%%%%%%%%%%%%%%%%%%%%%
%           Abstract
%%%%%%%%%%%%%%%%%%%%%%%%%%%%%%%%%%%%
\begin{abstract}
Let $\L$ be an arbitrary monomial algebra.
We investigate the stable category $\sGprojZ\L$ of graded Gorenstein-projective $\L$-modules and the orbit category $\sGprojZ\L/(1)$ induced by $\sGprojZ\L$ and the degree shift functor $(1)$. 
We prove that $\sGprojZ\L$ is triangle equivalent to the bounded derived category of a path algebra of Dynkin type $\A$ and that $\sGprojZ\L/(1)$ is triangle equivalent to the stable module category of a self-injective Nakayama algebra.
Both the path algebra and the self-injective Nakayama algebra will be given explicitly.
The latter result provides an explicit description of the stable category of (ungraded) Gorenstein-projective $\L$-modules.
\end{abstract}
%%%%%%%%%%%%%%%%%%%%%%%%%%%%%%%%%%%%
%           Abstract
%%%%%%%%%%%%%%%%%%%%%%%%%%%%%%%%%%%%

\maketitle
\tableofcontents

%%%%%%%%%%%%%%%%%%%%%%%%%%%%%%%%%%%%
%           Introduction ↓
%%%%%%%%%%%%%%%%%%%%%%%%%%%%%%%%%%%%
\section*{Introduction}
%%%%%%%%%%%%%%%%%%%%%%%%%%%%%%%%%%%%

In the 1960s, Auslander and Bridger \cite{Auslander-Bridger_1969} introduced the notion of the Gorenstein dimension of a finitely generated module over a two-sided Noetherian ring $R$ and studied those of Gorenstein dimension zero, which are sometimes called totally reflexive modules  \cite{Avramov-Martsinkovsky_2002}.
When $R$ is an Iwanaga-Gorenstein ring, these are also called maximal Cohen-Macaulay modules \cite{Buchweitz_1986}. 
In the 1990s, as a generalization of modules of Gorenstein dimension zero, Enochs and Jenda \cite{Enochs-Jenda_1995_MathZ} introduced the notion of Gorenstein-projective modules over an arbitrary ring. 
This terminology is appropriate because in Gorenstein homological algebra, Gorenstein-projective modules play the same role as projective modules in classical homological algebra.
Finally, as mentioned by Chen in \cite[Introduction]{X-WChen_2024}, Gorenstein-projective modules have applications in many areas.

It is well known that the stable category $\sGproj R$ of finitely generated Gorenstein-projective $R$-modules carries the structure of a triangulated category. 
A theorem of Buchweitz \cite[Theorem 4.4.1]{Buchweitz_1986} and Happel \cite[4.6]{Happel_1991} states that $\sGproj R$ is triangle equivalent to the singularity category $\D_{\rm sg}(\mod R)$ of $R$  in the case where $R$ is Iwanaga-Gorenstein. 
Consequently, many authors have investigated the stable categories of Gorenstein-projective modules over Iwanaga-Gorenstein rings and their singularity categories; see \cite{Chen-Geng-Lu_2015,Chen-Lu_2015,Iyama-Kimura-Ueyama_2024,Kalck_2014,Kalck-Iyama-Wemyss-Yang_2015,Kimura-Minamoto-Yamaura_2025,Liu-Lu_2015,MingLu_2016,Lu-Zhu_2021,Matsui-Takahashi_2017,Zhang_2011} for example.

On the other hand, among the results on $\sGproj R$ for more general rings $R$, such as \cite{MingLu_2019,Qin-Shen_2023,Zhang_2013}, Ringel has shown in \cite[Proposition 1]{Ringel_2013} that $\sGproj R$ is triangle equivalent to the stable module category of a connected self-injective Nakayama algebra when $R$ is a connected Nakayama algebra without simple projective modules. 
Moreover, Chen, Shen and Zhou have proven in \cite[Proposition 5.9]{Chen-Shen-Zhou_2018} that $\sGproj R$ is triangle equivalent to the stable module category of a radical square zero self-injective Nakayama algebra when $R$ is a monomial algebra without overlaps between any perfect paths. 
The notion of perfect paths was introduced in their paper, and they proved that the cyclic modules generated by perfect paths form a complete set of pairwise non-isomorphic finitely generated indecomposable non-projective Gorenstein-projective modules.

In this paper, we will extend the above result on $\sGproj R$ by Chen, Shen and Zhou to arbitrary monomial algebras.
In light of their result, it remains to consider the case where there exists an overlap between some perfect paths.
Therefore, we first examine perfect paths in this situation. 
Further, we will also introduce two partial orders $\preceq$ and $\leq$ on the set of perfect paths. 
It turns out that the sinks and sources in the Hasse quiver with respect to $\preceq$ play important roles in our development. 
Indeed, it will be proved that any perfect path can be uniquely decomposed into perfect paths that are sinks in the Hasse quiver with respect to $\preceq$; see Theorem \ref{u_claim_25}. This structure theorem is the first main result of this paper.

For a monomial algebra $\L = KQ/I$, we can consider it as a positively graded algebra by setting the degree of each arrow to one.
Hence we can define the category $\GprojZ\L$ of finitely generated graded Gorenstein-projective $\L$-modules and its stable category $\sGprojZ\L$. 
According to Lu and Zhu \cite[Lemma 4.2.1]{Lu-Zhu_2021}, there exists a $G$-covering $\sGprojZ\L \to \sGproj\L$, where $G$ stands for the cyclic group generated by the degree shift functor $(1)$. 
Note that $G$ acts freely on $\sGprojZ\L$.

First, by applying tilting theory, we show that $\sGprojZ\L$ is triangle equivalent to the bounded derived category of the path algebra of a disjoint union of Dynkin quivers of type $\A$; see Theorem \ref{u_claim_34}. 
This is the second main result of this paper. 
We remark that for another positive grading on $\L$ with a certain condition, one can obtain a similar result to the second main result; see Theorem \ref{u_claim_59}.
Second, through an investigation of the orbit category $\sGprojZ\L/(1)$, we show that $\sGproj\L$ is triangle equivalent to the stable module category of a self-injective Nakayama algebra; see Theorem \ref{u_claim_43}.
This is the third main result of this paper.
The path algebra in the second main result and the self-injective Nakayama algebra in the third main result are both determined explicitly.
Finally, we apply the third main result to certain classes of monomial algebras to recover and improve on previous results.

This paper is organized as follows.
%%%
In Section \ref{preliminaries}, we review basic facts related to Gorenstein-projective modules, positively graded algebras and monomial algebras.
%%%
In Section \ref{section_1}, we delve into the study of perfect paths with elementary proofs, examining perfect paths that have overlaps and presenting a structure theorem for perfect paths.
%%%
In Section \ref{section_3}, employing tilting theory, we investigate the stable category of graded Gorenstein-projective modules over a monomial algebra $\L$.
%%%
Finally, in Section \ref{section_4}, we study the stable category of Gorenstein-projective $\L$-modules via orbit categories.

Throughout this paper, let $K$ be a field.
By an algebra, we mean a finite dimensional associative $K$-algebra with a unit (except when considering the path algebra $KQ$ of a finite quiver $Q$ that contains a cycle). 
Further, a module means a finitely generated right module. 
For an algebra $\L$, we denote by $\mod \L$ the category of $\L$-modules and by $\proj\L$ the full subcategory of $\mod \L$ consisting of projective $\L$-modules.
For a Krull-Schmidt category $\mathcal{C}$, we denote by $\ind \mathcal{C}$ the set of indecomposable objects of $\mathcal{C}$ up to isomorphism.
Finally, we assume that any full subcategory of an additive category is closed under isomorphisms.

%%%%%%%%%%%%%%%%%%%%%%%%%%%%%%%%%%%%
%           Introduction ↑
%%%%%%%%%%%%%%%%%%%%%%%%%%%%%%%%%%%%

%%%%%%%%%%%%%%%%%%%%%%%%%%%%%%%%%%%%
%           Section ↓
%%%%%%%%%%%%%%%%%%%%%%%%%%%%%%%%%%%%
\section{Preliminaries} \label{preliminaries}
%%%%%%%%%%%%%%%%%%%%%%%%%%%%%%%%%%%%

%%%%%%%%%%%%%%%%%%%%%%%%%%%%%%%%%%%%
\subsection{Gorenstein-projective modules}  \label{preliminaries_1}
%%%%%%%%%%%%%%%%%%%%%%%%%%%%%%%%%%%%

In this subsection, we will recall some basic notions and facts related to Gorenstein-projective modules.
Throughout this subsection, we let $\L$ be an algebra.

A cochain complex  $P^{\bullet} : \cdots \to P^{i-1} \xrightarrow{d^{i-1}} P^{i} \xrightarrow{d^{i}} P^{i+1} \to \cdots$ of projective $\L$-modules is called {\it totally acyclic} \cite{Avramov-Martsinkovsky_2002} if both $P^{\bullet}$ and the Hom complex $\Hom_{\L}(P^\bullet, \L)$ are acyclic. 
We say that a $\L$-module $M$ is  {\it Gorenstein-projective} \cite{Enochs-Jenda_1995_MathZ} if there exists a totally acyclic complex $P^{\bullet}$ such that $\Ker d^{0}$ is isomorphic to $M$ in $\mod\L$.
We refer to \cite{X-WChen_2017} for their basic properties.

Let $\Gproj \L$ be the full subcategory of $\mod \L$ consisting of Gorenstein-projective $\L$-modules.
Since projective modules are Gorenstein-projective, we have that $\proj\L \subseteq \Gproj \L \subseteq \mod \L$.
We say that $\L$ is {\it CM-free} if $\proj\L = \Gproj \L$.
For example, algebras of finite global dimension are CM-free.
We say that $\L$ is {\it CM-finite} if there are only finitely many pairwise non-isomorphic indecomposable Gorenstein-projective $\L$-modules.
CM-free algebras and representation-finite algebras are both examples of CM-finite algebras. 
As will be mentioned in Remark \ref{remark_6}, monomial algebras are also CM-finite.
On the other hand, it is easily seen that $\Gproj \L = \mod \L$ if and only if $\L$ is self-injective (i.e.~$\L$ is injective as a one-sided $\L$-module).

Recall that the {\it stable category} $\smod\L$ of $\mod\L$ is defined as the category whose objects are the same as $\mod\L$ and whose morphisms are given by 
\begin{align}
    \sHom_{\L}(M, N) := \Hom_{\smod\L}(M, N) = \Hom_{\L}(M, N)/\P(M, N)
\end{align}
for any $M$ and $N \in \smod \L$, where $\P(M, N)$ denotes the space of morphisms from $M$ to $N$ that factor through a projective $\L$-module.
Let $\sGproj \L$ be the full subcategory of $\smod\L$ consisting of Gorenstein-projective $\L$-modules.
It is clear that $\sGproj \L = 0$ if and only if  $\L$ is CM-free.
On the other hand, the category $\Gproj \L$ is known to be a Frobenius exact category whose projective objects are precisely the projective $\L$-modules.
Therefore, the stable category $\sGproj \L$ carries a structure of a triangulated category; see \cite{Happel_Book}.

Recall that $\L$ is called {\it $d$-Iwanaga-Gorenstein} \cite{Iwanaga_1979} (or simply {\it Iwanaga-Gorenstein}) if both $\injdim_{\L}\L$ and $\injdim \L_{\L}$ are finite and at most $d$, where $\injdim_{\L}M$ denotes the injective dimension of $M$ in $\mod\L$.
Then there exists a triangle equivalence from the stable category $\sGproj \L$ to the {\it singularity category} $\D_{\rm sg}(\mod \L)$ of $\L$, where $\D_{\rm sg}(\mod \L)$ is defined to be the Verdier quotient of the bounded derived category $\D^{\rm b}(\mod \L)$ of $\L$ by (the thick subcategory equivalent to) the perfect derived category $\K^{\rm b}(\proj \L)$; see \cite{Buchweitz_1986}.

%%%%%%%%%%%%%%%%%%%%%%%%%%%%%%%%%%%%
\subsection{Positively graded algebras and forgetful functors}  \label{preliminaries_2}
%%%%%%%%%%%%%%%%%%%%%%%%%%%%%%%%%%%%

This subsection is devoted to recalling some basic facts in the representation theory of finite dimensional positively graded algebras from \cite{Gordon-Green_1982_Graded,Lu-Zhu_2021,Yamaura_2013} and presenting an observation that will be used in this paper.
Throughout this subsection, let $\L = \bigoplus_{i \geq 0} \L_i$ be a positively graded algebra (i.e.~a $\mathbb{Z}$-graded algebra satisfying $\L_i = 0$ for $i < 0$).

Let $\modZ\L$ be the category of graded $\L$-modules and $\projZ\L$ its full subcategory consisting of graded projective $\L$-modules.
Recall that the space of morphisms from $M$ to $N$ in $\modZ\L$ is defined by 
\begin{align}
    \HomZ_{\L}(M, N) 
    &:= \Hom_{\modZ\L}(M, N) \\
    &= \{\, f \in \Hom_\L(M, N) \mid f(M_i) \subseteq N_i \mbox{ for all } i\in \mathbb{Z} \,\}.
\end{align} 
For a graded $\L$-module $M = \bigoplus_{i \in \mathbb{Z}}M_{i}$ and an integer $j$, we define the degree shift $M(j) \in \modZ\L$ by $M(j)_i = M_{i+j}$ for $i\in \mathbb{Z}$.
This operation induces an automorphism $(j): \modZ\L \to \modZ\L$.
In particular, we have 
\begin{align} \label{eq_12}
    \Hom_{\L}^{\mathbb{Z}}(M(j), N(k)) \ = \      
    \Hom_{\L}^{\mathbb{Z}}(M, N(k-j))  \  = \ 
    \Hom_{\L}^{\mathbb{Z}}(M(j-k), N).
\end{align}
for any $j, k \in \mathbb{Z}$.

Let $F: \modZ\L \to \mod\L$ be the forgetful functor, and let $\smodZ \L$ be the stable category of $\modZ \L$, defined in a similar way to the case of $\smod \L$.
For any $M, N \in \smodZ\L$, we denote 
\begin{align}
    \sHomZ_{\L}(M, N) := \Hom_{\smodZ\L}(M, N).
\end{align}
For any $i \in \mathbb{Z}$, the automorphism $(i): \modZ \L \to \modZ \L$ induces an automorphism $(i) : \smodZ \L \to \smodZ \L$.
The following description of morphism spaces in $\modZ\L$ and $\smodZ\L$ is fundamental.

%%%%%%%%%%%%%%%%%%%%%%%%%%%%%%%%%%%%
%           Statement ↓
%%%%%%%%%%%%%%%%%%%%%%%%%%%%%%%%%%%%
\begin{lemma}[(cf.~\cite{Yamaura_2013}] \label{u_claim_37}
\begin{enumerate}

    \item For any $M, N \in \modZ\L$, we have  
    \begin{align}
    \Hom_{\L}(FM, FN) = \bigoplus_{i \in  \mathbb{Z}} \HomZ_{\L}(M, N(i)).
    \end{align} 

    \item For any $M, N \in \smodZ\L$, we have  
    \begin{align}
    \sHom_{\L}(FM, FN) = \bigoplus_{i \in  \mathbb{Z}} \sHomZ_{\L}(M, N(i)).
    \end{align} 
    
\end{enumerate}
\end{lemma}
\begin{proof}
See \cite[Section 2]{Gordon-Green_1982_Graded} for the proof of (1).
We now prove (2). 
Denote by $\P^{\mathbb{Z}}(M, N)$ be the subspace of $\HomZ_{\L}(M, N)$ consisting of morphisms factoring through a graded projective $\L$-module.
It is enough to show that 
\begin{align}
    \P(FM, FN) = \bigoplus_{i \in  \mathbb{Z}} \P^{\mathbb{Z}}(M, N(i)).
\end{align}
It follows from \cite[Proposition 1.3]{Gordon-Green_1982_Rep} that a graded $\L$-module $M$ is projective if and only if $FM$ is projective.
Then the inclusion $(\supseteq)$ follows from the\,\lq\lq only if\rq\rq part, and the reverse inclusion $(\subseteq)$ follows from the\,\lq\lq if\rq\rq part and \cite[Corollary 3.4]{Gordon-Green_1982_Graded}.
\end{proof}
%%%%%%%%%%%%%%%%%%%%%%%%%%%%%%%%%%%%
%           Statement ↑
%%%%%%%%%%%%%%%%%%%%%%%%%%%%%%%%%%%%

Replacing $\mod\L$ by $\modZ\L$, one can define the notion of {\it graded Gorenstein-projective $\L$-modules}.
Denote by $\GprojZ \L$ (resp.~$\sGprojZ \L$) the full subcategory of $\modZ\L$ (resp.~$\smodZ\L$) consisting of graded Gorenstein-projective $\L$-modules.
As in Section \ref{preliminaries_1}, we have that $\projZ\L \subseteq \GprojZ \L \subseteq \modZ \L$.
We say that $\L$ is {\it graded CM-free} if $\projZ\L = \GprojZ \L$, or equivalently, $\sGprojZ\L$ vanishes. 
The algebra $\L$ is said to be {\it graded CM-finite} if the number of the isomorphism classes of indecomposable graded Gorenstein-projective $\L$-modules is finite up to degree shift.
By \cite[Theorem 3.2]{Gordon-Green_1982_Graded},  a graded $\L$-module $M$ is indecomposable in $\modZ\L$ if and only if $FM$ is indecomposable in $\mod\L$.
Moreover, as mentioned in the proof of Lemma \ref{u_claim_37}, the graded $\L$-module $M$ is projective in $\modZ\L$ if and only if $FM$ is projective in $\mod\L$.
Thus if $\Lambda$ is graded CM-free, then it is graded CM-finite. 
It was proved in \cite[Lemma 2.2.6]{Lu-Zhu_2021} that the forgetful functor $F: \modZ\L \to \mod\L$ restricts to an additive functor $F_{G} : \GprojZ\L \to \Gproj\L$. 
Also, for any $i \in  \mathbb{Z}$, the automorphism $(i): \modZ\L \to \modZ\L$ restricts to an automorphism $(i): \GprojZ\L \to \GprojZ\L$.
Combining the above observations and \cite[Theorem 4.1]{Gordon-Green_1982_Graded}, one can obtain the following proposition, which explains a connection between CM-freeness and graded CM-freeness as well as CM-finiteness and graded CM-finiteness; compare  \cite[Theorem 6.1]{Minamoto-Yamaura_2020}.
Recall that a functor $F: \C \to \C^\prime$ is called {\it dense} if for any object $x^\prime \in \C^\prime$, there exists an object $x \in \C$ such that $x^\prime$ is isomorphic to $ F(x)$ in $\C^\prime$.

%%%%%%%%%%%%%%%%%%%%%%%%%%%%%%%%%%%%
%           Statement ↓
%%%%%%%%%%%%%%%%%%%%%%%%%%%%%%%%%%%%
\begin{proposition} \label{u_claim_46}
We have the following statements.
\begin{enumerate}
    \item If $\L$ is CM-free, then $\L$ is graded CM-free.
    \item If $\L$ is CM-finite, then $\L$ is graded CM-finite.
    \item If the forgetful functor $F_{G} : \GprojZ\L \to \Gproj\L$ is dense, then we have 
\begin{align}
    \ind \GprojZ\L = \{\,M (i) \mid M \in \ind \Gproj\L, i \in \mathbb{Z} \,\},
\end{align}
and the converses of (1) and (2) hold.
\end{enumerate}
\end{proposition}
\begin{proof}
We prove only (3); the proofs of (1) and (2) are left to the reader. 
To prove (3), it suffices to show the equality
\begin{align}
    \ind \GprojZ\L = \{\,M (i) \mid M \in \ind \Gproj\L, i \in \mathbb{Z} \,\}.
\end{align}
Since $F_{G} : \GprojZ\L \to \Gproj\L$ is dense by assumption, the representatives of the isomorphism classes of pairwise non-isomorphic indecomposable Gorenstein-projective $\L$-modules can be chosen to be gradable in the sense of \cite{Gordon-Green_1982_Graded}.
Now, let $X \in \ind \GprojZ\L$ be arbitrary.
Since $F_G(X) \in \ind \Gproj\L$, there exists a gradable module $M \in \ind \Gproj\L$ such that $F_G(X) \cong M = F(M)$ in $\Gproj\L$.
Note that $M$ is indecomposable as an object of $\modZ\L$.
Thus \cite[Theorem 4.1]{Gordon-Green_1982_Graded} implies that $X \cong M(i)$ in $\modZ\L$ for some $i \in \mathbb{Z}$.
This implies that 
\begin{align}
    \ind \GprojZ\L \subseteq \{\,M (i) \mid M \in \ind \Gproj\L, i \in \mathbb{Z} \,\}.
\end{align}
On the other hand, since $\GprojZ\L$ is closed under taking isomorphisms, the isomorphism $M \cong X(-i)$ in $\modZ\L$ implies  that the indecomposable graded $\L$-module $M$ is graded Gorenstein-projective.
This yields the desired reverse inclusion $(\supseteq)$.
\end{proof}
%%%%%%%%%%%%%%%%%%%%%%%%%%%%%%%%%%%%
%           Statement ↑
%%%%%%%%%%%%%%%%%%%%%%%%%%%%%%%%%%%%

We say that $\L$ is {\it graded Iwanaga-Gorenstein} (resp.~{\it graded self-injective}) if $\operatorname{gr.inj.dim}_{\L}\L < \infty$ and $\operatorname{gr.inj.dim} \L_{\L} < \infty$ (resp.~$\L$ is graded injective on one side), where $\mathrm{gr.inj.dim}_{\L}M$ denotes the injective dimension of $M$ in $\modZ\L$.
As mentioned in \cite[Section 2.1]{Lu-Zhu_2021}, the positively graded algebra $\L$ is graded Iwanaga-Gorenstein (resp.~graded self-injective) if and only if $\L$ is Iwanaga-Gorenstein (resp.~self-injective) as an ungraded algebra.
Thus we do not distinguish between being graded Iwanaga-Gorenstein (resp.~graded self-injective) and being Iwanaga-Gorenstein (resp.~self-injective).
As in the ungraded case, $\GprojZ\L = \modZ \L$ if and only if $\L$ is  self-injective.
The category $\GprojZ\L$ is a Frobenius exact category such that $\projZ\L$ is the class of projective objects, so that the stable category $\sGprojZ\L$ has a structure of a triangulated category.
If $\L$ is Iwanaga-Gorenstein, then $\sGprojZ\L$ is triangle equivalent to the {\it graded singularity category} $\D_{\rm sg}(\modZ\L)$ of $\L$, where $\D_{\rm sg}(\modZ\L)$ is defined by the Verdier quotient of $\D^{\rm b}(\modZ\L)$ by $\K^{\rm b}(\projZ \L)$.

We end this subsection with a preparation of a covering technique for  $\sGproj \L$.
We follow the terminology of \cite{Asashiba_2011}.

A pair $(\C, A)$ of a $K$-linear category $\C$ (i.e.~a category such that the morphism sets are $K$-vector spaces, and the composition are $K$-bilinear) and a group homomorphism $A : G \to \Aut(\C)$ is called a category with a {\it $G$-action} or a {\it $G$-category}, where $\Aut(\C)$ is the group of automorphisms of $\C$. 
If there is no confusion, we simply write $\C = (\C, A)$. 
For a $G$-category $\C$, we will write $\alpha x = A(\alpha)(x)$ for all $\alpha \in G$ and $x \in \C$. 
We denote by $1$ the identity of $G$.
Let $F:\C^\prime \to \C $ be a $K$-linear functor between $K$-linear categories with $\C^\prime$ a $G$-category. 
An {\it invariance adjuster} of $F$ is a family $\phi = (\phi_\alpha)_{\alpha\in G}$ of functorial isomorphisms $\phi_\alpha : F \to  F A(\alpha)$ with  $\alpha  \in G$ satisfying the following two conditions:
\begin{enumerate}
    \item $\phi_1: F \to F$ is the identify functor $\mathbbm{1}_F$ of $F$.
    \item The following diagram is commutative for all $\alpha,\beta \in G$:
    \begin{align}
        \xymatrix{
        F \ar[r]^-{\phi_\alpha} \ar[d]_-{\phi_{\beta\alpha}} & FA(\alpha) \ar[d]^-{\phi_\beta A(\alpha)}\\
        FA(\beta\alpha) \ar@{=}[r]&   FA(\beta) A(\alpha).
        }
    \end{align}
\end{enumerate}
We call the pair $(F, \phi)$ a {\it $G$-invariant functor}. 
As mentioned in \cite[Remark 1.2]{Asashiba_2011}, the condition (1) automatically follows from (2).
Let $F= (F, \phi): \C \to \C^\prime$ be a $G$-invariant functor.
For any $x, y  \in \C$, we define the following two $K$-linear maps $F^{(1)}_{x,y}$ and $F^{(2)}_{x,y}$:\vspace{1mm}
\begin{align}
    &F^{(1)}_{x,y}: \bigoplus_{\alpha \in G}\Hom_{\C}(\alpha x, y) \to  \Hom_{\C^\prime}(F(x), F(y)), \quad \  ( f_\alpha)_{\alpha \in G} \mapsto \sum_{\alpha \in G} F(f_\alpha)\cdot \phi_{\alpha}(x);\\[2mm]
    &F^{(2)}_{x,y} : \bigoplus_{\beta \in G}\Hom_{\C}(x,\beta y) \to \Hom_{\C^\prime}(F(x), F(y)), \quad \ ( f_\beta)_{\beta \in G} \mapsto \sum_{\beta\in G}\phi_{\beta^{-1}}(\beta y) \cdot F(f_\beta),
\end{align}
where 
\begin{align}
    F(f_\alpha)\cdot \phi_{\alpha}(x): F(x) \to F(y)
\end{align}
is the composition of $\phi_{\alpha}(x) : F(x) \to FA(\alpha)(x)$ and $F(f_\alpha): F(\alpha x) \to F(y)$, and
\begin{align}
    \phi_{\beta^{-1}}(\beta y) \cdot F(f_\beta): F(x) \to F(y)
\end{align}
is the composition of $F(f_\beta) : F(x) \to F(\beta y)$ and $\phi_{\beta^{-1}}(\beta y): F(\beta y) \to FA(\beta^{-1})(\beta y) = F(y)$.
We know from \cite[Proposition 1.6]{Asashiba_2011} that $F^{(1)}_{x,y}$ is bijective if and only if $F^{(2)}_{x,y}$ is bijective for any $x,y \in \C$.
We say that a $G$-invariant functor $F= (F, \phi): \C \to \C^\prime$ is a {\it $G$-precovering} if for any $x, y  \in \C$, the $K$-linear map 
\begin{align}
    &F^{(1)}_{x,y}: \bigoplus_{\alpha \in G}\Hom_{\C}(\alpha x, y) \to  \Hom_{\C^\prime}(F(x), F(y))
\end{align}
is bijective.
Further, the $G$-invariant functor $F= (F, \phi): \C \to \C^\prime$ is called a {\it $G$-covering} if $F$ is a $G$-precovering and  is dense.

The forgetful functor $F: \modZ\L \to \mod\L$ induces an additive functor $\Tilde{F}: \smodZ\L \to \smod\L$, which restricts to an additive functor $\Tilde{F}_G : \sGprojZ\L\to\sGproj\L$.
Note that $\Tilde{F}_G$ is also induced by the functor $F_G : \GprojZ\L\to\Gproj\L$.
Since $F: \modZ\L \to \mod\L$ is an exact functor and preserves projective objects,  $\Tilde{F}_G : \sGprojZ\L\to\sGproj\L$ becomes a triangulated functor.
Also,  the automorphism $(i): \modZ\L \to \modZ\L$ with $i \in \mathbb{Z}$ induces an automorphism $(i): \smodZ\L \to \smodZ\L$, which restricts to an automorphism $(i): \sGprojZ\L \to \sGprojZ\L$.
Note that this automorphism is also induced by the automorphism $(i): \GprojZ\L \to \GprojZ\L$.
Then, since $F = F \circ (i)$ for any $i \in \mathbb{Z}$, together with Lemma \ref{u_claim_37}, one concludes that $F, F_G, \Tilde{F}$ and $\Tilde{F}_G$ are all $G$-precoverings, where $G$ is the cyclic group generated by the automorphism $(1)$.
Note that the invariance adjusters of $F_G, \Tilde{F}$ and $\Tilde{F}_G$ are obtained by replacing $F$ with the corresponding functor among $F_G, \Tilde{F}$ and $\Tilde{F}_G$ in the formula $F = F \circ (i)$.
We also note that $G$ freely acts on all the categories $\modZ\L, \GprojZ\L, \smodZ\L$ and $\sGprojZ\L$.

The following lemma is an easy consequence of \cite[Corollary 3.4]{Gordon-Green_1982_Graded}, which asserts that any projective $\L$-modules are gradable.

%%%%%%%%%%%%%%%%%%%%%%%%%%%%%%%%%%%%
%           Statement ↓
%%%%%%%%%%%%%%%%%%%%%%%%%%%%%%%%%%%%
\begin{lemma} \label{claim_38} 
The following conditions are equivalent.
\begin{enumerate}
    \item $F_G: \GprojZ\L \to \Gproj\L$ is dense.
    \item $F_G: \GprojZ\L \to \Gproj\L$ is a $G$-covering.
    \item $\Tilde{F}_G: \sGprojZ\L \to \sGproj\L$ is a $G$-covering.
\end{enumerate}
\end{lemma}
%%%%%%%%%%%%%%%%%%%%%%%%%%%%%%%%%%%%
%           Statement ↑
%%%%%%%%%%%%%%%%%%%%%%%%%%%%%%%%%%%%

For a $G$-category $\mathcal{C}$,  when $G$ is the cyclic group generated by an automorphism $F:\mathcal{C} \to \mathcal{C}$, we denote by $\mathcal{C}/F$ the orbit category $\mathcal{C}/G$.

Assume that $\Tilde{F}_G: \sGprojZ\L \to \sGproj\L$ is a $G$-covering.
Then, by \cite[Theorem 2.9]{Asashiba_2011}, there exists an equivalence $H: \sGprojZ\L/(1) \xrightarrow{\,\sim\,} \sGproj\L$ of additive categories
such that $\Tilde{F}_G = HP$ (as $G$-invariant functors), where
$P: \sGprojZ\L \to \sGprojZ\L/(1)$ 
is the canonical functor:
\begin{align}
\xymatrix@M=2mm@R=5mm@C=0mm{
    \sGprojZ\L \ar[rr]^-{\Tilde{F}_G} \ar[rd]_-{P} & & \sGproj\L \\
     & \sGprojZ\L/(1) \ar[ru]_-{H}^{\sim} & 
}\end{align}
Thus the orbit category $\sGprojZ\L/(1)$ becomes a triangulated category whose triangulated structure is derived from that of the triangulated category $\sGproj\L$ via the equivalence $H$.
Finally, since $\Tilde{F}_G$ and $H$ are both triangulated functors, the canonical functor $P$ is also triangulated.

%%%%%%%%%%%%%%%%%%%%%%%%%%%%%%%%%%%%
\subsection{Monomial algebras and their Gorenstein-projective modules} \label{monomial algebras} \label{preliminaries_3}
%%%%%%%%%%%%%%%%%%%%%%%%%%%%%%%%%%%%

In this subsection, we review some notations and facts related to monomial algebras, including an explicit description of Gorenstein-projective modules by Chen, Shen and Zhou \cite{Chen-Shen-Zhou_2018}.
For the basics on quivers and bound quiver algebras, we refer to \cite{ASS06}.

Let $Q = (Q_0, Q_1, s, t)$ be a finite quiver. 
A {\it path of length $n$} in $Q$ is a sequence $p = a_1 a_2 \cdots a_n$ of arrows $a_i \in Q_1$ with $t(a_i) = s(a_{i+1})$ for all $1 \leq i \leq n-1$.  
We define its source $s(p) := s(a_1)$ and its target $t(p) := t(a_n)$.
The length of a path $p$ is denoted by $l(p)$.
For a vertex $v$ in $Q$, we associate a {\it trivial path} $e_{v}$ of length zero. 
It is understood that $s(e_{v}) = v = t(e_{v})$.

We denote by $\B$ the set of paths, by $\mathcal{B}_i$ the set of paths of length $i$ and by $\B_{\geq i}$  the set of paths of length $\geq i$. 
Thus $\B_{>0} := \B_{\geq 1}$ is the set of non-trivial paths, and we have $\B_{\geq i} = \bigcup_{j\geq i} \B_j$ and $\B = \B_{\geq 0}$.

The concatenation of two paths $p$ and $q$ with $t(p) = s(q)$ is denoted by $pq$.
For any $p$ and $q \in \B$, we say that $q$ is a {\it subpath} of $p$ if $p = p^{\prime} q p^{\prime\prime}$ for some $p^{\prime}, p^{\prime\prime} \in \B$.
When $p^{\prime}$ is trivial, we refer to the subpath $q$ as a {\it left divisor} of $p$. Dually, when $p^{\prime\prime}$ is trivial, we refer to the subpath $q$ as a {\it right divisor} of $p$.
A subpath $q$ of $p$ is called {\it proper} if $q \not=p$.

A {\it cycle} in $Q$ is a path $c$ such that $s(c) = t(c)$.
We say that two cycles $c$ and $c^\prime$ are {\it equivalent} if  $c$ coincides with $c^\prime$ up to cyclic permutation.

Let $S$ be a subset of $\B$. We say that $p \in S$ is {\it left minimal in $S$} if there exists no proper left divisor of $p$ that belongs to $S$ (i.e.~there exists no path $s \in S$ such that $p = s p^\prime$ for some $p^\prime \in \B_{>0}$). Dually, we say that $p \in S$ is {\it right minimal in $S$} if there exists no proper right divisor of $p$ that belongs to $S$. 
We say that $p \in S$ is {\it minimal in $S$} if there exists no proper subpath of $p$ that belongs to $S$.

A bound quiver algebra $KQ/I$ is called {\it monomial} if the admissible ideal $I$ of the path algebra $KQ$ is generated by paths. 
Note that such paths are of length at least two. 
When $I$ is generated by paths of length two, the algebra $KQ/I$ is called {\it quadratic monomial}.

Let $\L = KQ/I$ be a monomial algebra.
We denote by $\F$ the set of minimal paths in the set $\{$paths belonging to $I\}$.
Note that $\F$ generates $I$ as an ideal of $KQ$.
In particular, $\F$ is a finite set.
We say that $p \in \B$ is {\it non-zero in $\L$} if the canonical image $p+I$ in $\L$ is non-zero, or equivalently, $p$ does not belong to $I$, which is equivalent to saying that $p$ does not contain any path of $\F$ as a subpath. 
It is easily seen that the non-zero paths form a $K$-basis of $\L$.
Finally, for each $p \in \B$, we simply denote the canonical image $p+I$ in $\L$ by $p$.
Moreover, we will write $p=0$ in $\L$ when $p$ lies in $I$.

Let $p$ be a non-trivial path that is non-zero in $\L$.
Then we define $L(p)$ as the set of right minimal paths in the set $\{$non-zero paths $q \mid t(q) = s(p)$ and $qp = 0$ in $\L \}$.  
Dually, we define $R(p)$ as the set of left minimal paths in the set $\{$non-zero paths $q \mid t(p) = s(q)$ and $pq = 0$ in $\L \}$.

%%%%%%%%%%%%%%%%%%%%%%%%%%%%%%%%%%%%
%           Statement ↓
%%%%%%%%%%%%%%%%%%%%%%%%%%%%%%%%%%%%
\begin{definition}[{cf.~\cite[Definitions 3.3 and 3.7]{Chen-Shen-Zhou_2018}}]{\rm  \label{definition_1}
Let $\L$ be a monomial algebra.
\begin{enumerate}
\item A pair $(p, q)$ of non-zero paths in $\L$ is said to be {\it perfect} if the following conditions are satisfied:
    \begin{enumerate}[label=(P\arabic*)]
    
        \item $p$ and $q$ are both non-trivial and satisfy $t(p) = s(q)$ and $pq=0$ in $\L$;
        
        \item  If $pq^\prime = 0$ for a non-zero path $q^\prime$ with $t(p) = s(q^\prime)$, then $q$ is a left divisor of $q^\prime$  (in other words, $R(p) =\{q\}$);
        
        \item  If $p^\prime q = 0$ for a non-zero path $p^\prime$ with $t(p^\prime) = s(q)$, then $p$ is a right divisor of $p^\prime$ (in other words,  $L(q)=\{p\}$).   
        
    \end{enumerate}

\item A sequence $(p_1, \ldots, p_n, p_{n+1}=p_1)$ of non-zero paths in $\L$ is called a {\it perfect path sequence} if $(p_i, p_{i+1})$ is a perfect pair for all $ 1 \leq i \leq n$. 
A path appearing in a perfect path sequence is called a {\it perfect path}.

\item A perfect path sequence $(p_1, \ldots, p_n, p_{n+1}=p_1)$ is called {\it minimal} if $p_i \not= p_j$ for any $1\leq i\not=j \leq n$.

\end{enumerate} 
}\end{definition}
%%%%%%%%%%%%%%%%%%%%%%%%%%%%%%%%%%%%
%           Statement ↑
%%%%%%%%%%%%%%%%%%%%%%%%%%%%%%%%%%%%

%%%%%%%%%%%%%%%%%%%%%%%%%%%%%%%%%%%%
%           Statement ↓
%%%%%%%%%%%%%%%%%%%%%%%%%%%%%%%%%%%%
\begin{remark}{\rm \label{remark_4}
    For a perfect pair $(p, q)$, we have the following short exact sequence of $\L$-modules:
    \begin{align} \label{eq_1}
        0 \rightarrow q\L \xrightarrow{\,\iota\,} e_{t(p)}\L \xrightarrow{\,\pi\,} p\L \rightarrow 0,
    \end{align}

    where $\iota$ is the canonical inclusion, and $\pi$ is the left multiplication by $p$, that is, $\pi(x) = px$ for all $x \in e_{t(p)}\L$. 

}\end{remark}
%%%%%%%%%%%%%%%%%%%%%%%%%%%%%%%%%%%%
%           Statement ↑
%%%%%%%%%%%%%%%%%%%%%%%%%%%%%%%%%%%%

For a monomial algebra $\L$, we denote by $\mathbb{P}_{\L}$ the set of perfect paths.
The following result describes indecomposable non-projective Gorenstein-projective $\L$-modules.
Note that $\ind \sGproj\L$ can be identified with the set of the isomorphism classes of indecomposable non-projective Gorenstein-projective $\L$-modules.

%%%%%%%%%%%%%%%%%%%%%%%%%%%%%%%%%%%%
%           Statement ↓
%%%%%%%%%%%%%%%%%%%%%%%%%%%%%%%%%%%%
\begin{theorem}[{\cite[Theorem 4.1]{Chen-Shen-Zhou_2018}}] \label{Theorem4.1_Chen-Shen-Zhou_2018}
Let $\L$ be a monomial algebra.
Then there is a bijection 
\begin{align}
    \mathbb{P}_{\L}
     \ \xleftrightarrow{\ 1:1\ } \  
    \ind \sGproj\L,
\end{align}
where a perfect path $p$ is sent to the cyclic $\L$-module $p\L$.
\end{theorem}
%%%%%%%%%%%%%%%%%%%%%%%%%%%%%%%%%%%%
%           Statement ↑
%%%%%%%%%%%%%%%%%%%%%%%%%%%%%%%%%%%%

%%%%%%%%%%%%%%%%%%%%%%%%%%%%%%%%%%%%
%           Statement ↓
%%%%%%%%%%%%%%%%%%%%%%%%%%%%%%%%%%%%
\begin{remark}{\rm \label{remark_6}
%%%% 
%%%%
%It follows from the theorem that monomial algebras $\L$ are always CM-finite.  
%%%%
%%%%
The set $\mathbb{P}_\L$ is finite for any monomial algebra $\L$ since the number of paths that are non-zero in $\L$ is finite.
Thus the theorem implies that $\L$ is CM-finite. 
%%%% 
%%%%
It also follows that $\L$ is CM-free precisely when $\mathbb{P}_\L$ is empty.  
}\end{remark}
%%%%%%%%%%%%%%%%%%%%%%%%%%%%%%%%%%%%
%           Statement ↑
%%%%%%%%%%%%%%%%%%%%%%%%%%%%%%%%%%%%

Let $a$ be an element of an algebra $\L$.
Then there exists a monomorphism of left $\L$-modules
\begin{align}
        \theta_{a}^\prime:\L a \rightarrow \Hom_{\L}(a\L, \L)
\end{align}
given by $\theta_{a}^\prime(xa)(ay)= xay$ for all $x, y \in \L$. 
By \cite[Lemma 2.3]{Chen-Shen-Zhou_2018}, if $\theta_{a}^\prime$ is an isomorphism, then we obtain the following isomorphisms of $K$-vector spaces for any $b\in \L$:
\begin{align}  \label{eq_16} 
\theta_{a, b}^\prime : b\L \cap \L a \xrightarrow{\ \sim\ }  \Hom_{\L}(a\L, b\L) \quad \mbox{and} \quad  \dfrac{b\L \cap \L a}{b\L a} \xrightarrow{\ \sim\ }  \sHom_{\L}(a\L, b\L).
\end{align}
Note that $\theta_{a, b}^\prime$ is the restriction of $\theta_{a}^\prime$ to the subspace $ b\L \cap \L a$.
It also follows from \cite[Lemma 3.6]{Chen-Shen-Zhou_2018} that if $\L$ is monomial, then the monomorphism $\theta^\prime_p$ is an isomorphism for any perfect path $p$.

Let $p$ and $q$ be perfect paths in a monomial algebra $\L$.
We say that {\it there exists an overlap between $p$ and $q$} (and that {\it $p$ overlaps $q$}) if one of the following conditions is satisfied:
\begin{enumerate}
    \item [(O1)] $p=q$, and $p=p^\prime x$ and $q=x q^\prime$ for some $x, p^\prime$ and $q^\prime \in \B_{>0}$ satisfying $p^\prime x q^\prime \not = 0$.
    \item [(O2)] $p\not=q$, and $p=p^\prime x$ and $q=x q^\prime$ for some $x \in \B_{>0}$ and $p^\prime$ and $q^\prime \in \B$ satisfying  $p^\prime x q^\prime \not = 0$.
\end{enumerate}

%%%%%%%%%%%%%%%%%%%%%%%%%%%%%%%%%%%%
%           Statement ↓
%%%%%%%%%%%%%%%%%%%%%%%%%%%%%%%%%%%%
\begin{proposition}[{\cite[page 1130]{Chen-Shen-Zhou_2018}}] \label{Chen-Shen-Zhou_2018_page1130} 
Let $\L$ be a monomial algebra. 
Then the following statements hold for any perfect paths $p$ and $q$.
\begin{enumerate}

\item Condition $({\rm O}1)$ holds for $p$ and $q$ if and only if the one-dimensional vector space $K \id_{p\L}$ is properly contained in  $\sHom_{\L}(p\L, p\L)$.

\item Condition $({\rm O}2)$ holds for $p$ and $q$ if and only if $\sHom_{\L}(p\L, q\L) \not = 0$.

\end{enumerate}
\end{proposition}
%%%%%%%%%%%%%%%%%%%%%%%%%%%%%%%%%%%%
%           Statement ↑
%%%%%%%%%%%%%%%%%%%%%%%%%%%%%%%%%%%%

We say that {\it there exists no overlap in a monomial algebra $\L$} if there exists no overlap between any two perfect paths. 
In this case, it follows from the proposition that there exists no non-trivial homomorphism in $\sGproj\L$. 
This class of monomial algebras includes quadratic monomial algebras, such as gentle algebras.
In \cite[Proposition 5.9]{Chen-Shen-Zhou_2018}, Chen, Shen and Zhou provided necessary and sufficient conditions under which there exists no overlap in $\L$ in terms of the triangulated structure of $\sGproj \L$, and moreover they described $\sGproj \L$ explicitly.  
Motivated by this result, we focus on perfect paths that have overlaps in the next section.

%%%%%%%%%%%%%%%%%%%%%%%%%%%%%%%%%%%%
%           Section ↑
%%%%%%%%%%%%%%%%%%%%%%%%%%%%%%%%%%%%

%%%%%%%%%%%%%%%%%%%%%%%%%%%%%%%%%%%%
%           Section ↓
%%%%%%%%%%%%%%%%%%%%%%%%%%%%%%%%%%%%
\section{Perfect paths} \label{section_1}
%%%%%%%%%%%%%%%%%%%%%%%%%%%%%%%%%%%%

This section is divided into three subsections. 
In the first, as mentioned at the end of the previous section, we investigate perfect paths having overlaps. 
In the second, we introduce two partial orders on the set of perfect paths and examine perfect paths with respect to these partial orders.
In the third, we establish a structure theorem for perfect paths, which is the main result of this section.

In the rest of this paper, let $\L = KQ/I$ be a monomial algebra and $\F$ the set of minimal paths in $I$.

%%%%%%%%%%%%%%%%%%%%%%%%%%%%%%%%%%%%
%           Subsection ↓
%%%%%%%%%%%%%%%%%%%%%%%%%%%%%%%%%%%%
\subsection{Perfect paths with overlaps} \label{section_1_1}
%%%%%%%%%%%%%%%%%%%%%%%%%%%%%%%%%%%%

A perfect path sequence $(p_{1}, \ldots, p_{n}, p_{n+1} = p_{1})$ gives rise to a cycle $p_1 \cdots p_n$. 
We refer to the shortest cycle $c$ such that $p_1 \cdots p_n = c^l$ for some $l>0$ as the {\it underlying cycle associated with the perfect path $p_1$} and denote it by $c_{p_1}$. 
To check the well-definedness of the notation $c_{p_1}$, one can use the observation that any perfect path sequence $(p_1, \ldots, p_n, p_{n+1} = p_1)$ can be recovered from any given path $p_i$ in the sequence up to cyclic permutation by inductively setting $p_{i+1}$ to be the unique element of $R(p_i)$.

Cyclic permutations define an equivalence relation on the set of underlying cycles. 
We denote by $\C(\L)$ the set of the equivalence classes of underlying cycles modulo this equivalence relation.
It follows from the above observation that $c_{p_i} = c_{p_j}$ in $\C(\L)$ for any $1 \leq i, j \leq n$.

Let us start with the following lemma.

%%%%%%%%%%%%%%%%%%%%%%%%%%%%%%%%%%%%
%           Statement ↓
%%%%%%%%%%%%%%%%%%%%%%%%%%%%%%%%%%%%
\begin{lemma} \label{u_claim_49}
Let $(p_{1}, \ldots, p_{n}, p_{n+1} = p_{1})$ and $(q_{1}, \ldots, q_{m}, q_{m+1} = q_{1})$ be perfect path sequences.
Assume that $p_1 = q_0^\prime x_1, q_1 = x_1 q_1^\prime$ and $q_0^\prime x_1 q_1^\prime \not=0$ in $\L$ for some $x_1 \in \B_{>0}$ and $q_0^\prime, q_1^\prime \in \B$.
Then there exist families $\{x_l \in \B_{>0}\}_{l\geq 2}$  and $\{q_l^\prime \in \B\}_{l\geq 2}$ such that 
\begin{align}
    p_{l} =  q_{l-1}^\prime x_l \quad \mbox{and} \quad q_{l} = x_l q_{l}^\prime \quad \mbox{for $l \geq 2$,}
\end{align}
where we consider the index $l$ of $p_l$ (resp.~$q_l$) modulo $n$ (resp.~$m$).
Moreover, the following statements hold.
\begin{enumerate}

    \item The pairs $(p_l q_{l}^\prime, x_{l+1})$ and $(x_{l+1}, p_{l+2} q_{l+2}^\prime)$ are perfect for $l \geq 1$.

    \item If $q_0^\prime$ is trivial, then the families  satisfy the following conditions:
    \begin{enumerate}
        \item[$($a$)$] $x_{2i} = q_{2i}$ and $x_{2i+1} = p_{2i+1}$ for $i \geq 1$.
        \item[$($b$)$] $q_{2i}^\prime$ is trivial for $i \geq 1$.
    \end{enumerate}
    If, furthermore, $q_1^\prime$ is non-trivial, then $(q_{2i+1}^\prime, q_{2i+2}p_{2i+3})$ and $(q_{2i+2}p_{2i+3}, q_{2i+3}^\prime)$ are perfect for $i \geq 0$.

    \item If $q_1^\prime$ is trivial, then the families  satisfy the following conditions:
    \begin{enumerate}
        \item[$($a$)$] $x_{2i} = p_{2i}$ and $x_{2i+1} = q_{2i+1}$ for $i \geq 1$.
        \item[$($b$)$] $q_{2i+1}^\prime$ is trivial for $i \geq 1$.
    \end{enumerate}
    If, furthermore, $q_0^\prime$ is non-trivial, then $(q_{2i}^\prime, q_{2i+1}p_{2i+2})$ and $(q_{2i+1}p_{2i+2}, q_{2i+2}^\prime)$ are perfect for $i \geq 0$.
\end{enumerate} 
\end{lemma}

\begin{proof}
Since $p_1(q_1^\prime q_2) = q_0^\prime x_1(q_1^\prime q_2) = q_0^\prime q_1q_2 = 0$ in  $\L$,  the fact that $(p_1, p_2)$ is perfect implies that $q_1^\prime q_2$ contains $p_2$ as a left divisor.
Since $p_1 q_1^\prime = q_0^\prime x_1 q_1^\prime  \not= 0$ in $\L$, $p_2$ must contain $q_1^\prime$ as a proper left divisor, so that $p_2 = q_1^\prime x_2$ for some $x_2 \in \B_{>0}$.
Now, since $p_2 = q_1^\prime x_2$ is a left divisor of $q_1^\prime q_2$, it follows that  $q_2 = x_2 q_2^\prime$ for some $q_2^\prime \in \B$.
On the other hand, we have $q_1^\prime x_2 q_2^\prime = q_1^\prime q_2 \not=0$ (otherwise, $q_1^\prime$ contains $q_1=x_1q_1^\prime$ as a right divisor).
So far our situation is as follows: 
%%%%%%%%%%%%%%%%%%%%%%%%%%%%%
%   diagram ↓
%%%%%%%%%%%%%%%%%%%%%%%%%%%%%
\[\begin{tikzpicture}[scale=1]
%%%%%%%%%%%%%
% paths
%%%%%%%%%%%%%
\node at (2,1.25) {$p_{1}$};
\node at (6,1.25) {$p_{2}$};
\node at (1.5,0.35) {$q_{0}^\prime$};
\node at (3.25,0.25) {$x_{1}$};
\node at (5.3,0.35) {$q_{1}^\prime$};
\node at (7.3,0.25) {$x_2$};
\node at (9,0.35) {$q_{2}^\prime$};
\node at (4.5,-0.75) {$q_{1}$};
\node at (8,-0.75) {$q_{2}$};
%%%%%%%%%%%%%
% vertices
%%%%%%%%%%%%%
%\node (a); 
\node (a1) at (0,1) {$\cdot$};
\node (a2) at (4,1) {$\cdot$};
\node (a3) at (8.2,1) {$\cdot$};
%\node (b);
\node (b1) at (0,0) {$\cdot$};
\node (b2) at (2.5,0) {$\cdot$};
\node (b3) at (4,0) {$\cdot$};
\node (b4) at (6.3,0) {$\cdot$};
\node (b5) at (8.2,0) {$\cdot$};
\node (b6) at (10,0) {$\cdot$};
%\node (c);
\node (c1) at (2.5,-1) {$\cdot$};
\node (c2) at (6.3,-1) {$\cdot$};
\node (c3) at (10,-1) {$\cdot$};
%%%%%%%%%%%%%
% arrows
%%%%%%%%%%%%%
% a
\draw[->] (a1) -- (a2) ;
\draw[->] (a2) -- (a3) ;
% b 
\draw[->] (b1) -- (b2) ;
\draw[->] (b2) -- (b3) ;
\draw[->] (b3) -- (b4) ;
\draw[->] (b4) -- (b5) ;
\draw[->] (b5) -- (b6) ;
% c
\draw[->] (c1) -- (c2) ;
\draw[->] (c2) -- (c3) ;
\draw[dashed] (a1) -- (b1) ;
\draw[dashed] (a2) -- (b3) ;
\draw[dashed] (a3) -- (b5) ;
\draw[dashed] (b2) -- (c1) ;
\draw[dashed] (b4) -- (c2) ;
\draw[dashed] (b6) -- (c3) ;
\end{tikzpicture}
\]
%%%%%%%%%%%%%%%%%%%%%%%%%%%%%
% diagram ↑
%%%%%%%%%%%%%%%%%%%%%%%%%%%%%

Applying the above argument repeatedly, we obtain infinitely many non-trivial paths $x_2, x_3, \ldots$ and paths $q_2^\prime, q_3^\prime, \ldots$ with the desired property.

For (1), we only address the pair $(p_l q_{l}^\prime, x_{l+1})$; the proof for the other is analogous.  
First, we have $(p_l q_{l}^\prime)x_{l+1} = p_{l}p_{l+1} = 0$ in $\L$.
Second, let $u$ be a non-zero path with $s(u) = t(p_l q_{l}^\prime)$ such that $(p_l q_{l}^\prime) u = 0$ in $\L$.
Since $(p_{l}, p_{l+1})$ is a perfect pair, $q_{l}^\prime u$ contains $p_{l+1}=q_{l}^\prime x_{l+1}$ as a left divisor, which yields that $u$ contains $x_{l+1}$ as a left divisor.
Third, take a non-zero path $v$  with $t(v) = s(x_{l+1})$ satisfying $v x_{l+1}= 0$ in $\L$. 
Since $v q_{l+1} = v(x_{l+1} q_{l+1}^\prime) = 0$ in $\L$, together with the fact that $(q_{l}, q_{l+1})$ is perfect, it follows that $v = v^\prime q_{l} = v^\prime x_{l} q_{l}^\prime$ for some $v^\prime \in \B$.
Furthermore, since $(v^\prime x_{l}) p_{l+1} = (v^\prime x_{l})q_{l}^\prime x_{l+1} = v x_{l+1} = 0$, we can deduce that $v^\prime x_{l}$ contains $p_{l}$ as a right divisor.
As a result, $v = v^\prime x_{l} q_{l}^\prime$ contains $p_{l}q_{l}^\prime$ as a right divisor.
Consequently, we see from Definition \ref{definition_1} that $(p_l q_{l}^\prime, x_{l+1})$ is a perfect pair.

For (2), suppose that $q_0^\prime$ is trivial, that is, $q_1 = p_1 q_1^\prime$. 
Then $p_1(q_1^\prime q_2) = x_1(q_1^\prime q_2) = q_1q_2 \in \F$.
Here, we use the fact that $pq \in \F$ for any perfect pair $(p, q)$; see \cite[Lemma 3.4]{Chen-Shen-Zhou_2018}.
Thus there exists an $r \in \B$ such that $q_1^\prime q_2= p_2 r$. 
But then, since both $p_1(q_1^\prime q_2)$ and $p_1 p_2$ belong to $\F$, we must have that $r$ is trivial, that is, $q_1^\prime q_2 = p_2$, so that $x_2 = q_2$, and $q_2^\prime$ is trivial. 
Moreover, $p_3 = q_2^\prime x_3 = x_3$.
Applying this argument repeatedly, one can check that $\{x_l \in \B_{>0}\}_{l\geq 2}$ and $\{q_l^\prime \in \B\}_{l\geq 2}$ have the desired property.
Finally, an argument as in the proof of (1) shows the last statement.

For (3), suppose that $q_1^\prime$ is trivial.
Then we see that $p_2 = x_2$ and $q_2 = x_2 q_2^\prime$ with $x_2 \in \B_{>0}$.
Thus (3) is an immediate consequence of (2).
\end{proof}
%%%%%%%%%%%%%%%%%%%%%%%%%%%%%%%%%%%%
%           Statement ↑
%%%%%%%%%%%%%%%%%%%%%%%%%%%%%%%%%%%%

%%%%%%%%%%%%%%%%%%%%%%%%%%%%%%%%%%%%
%           Statement ↓
%%%%%%%%%%%%%%%%%%%%%%%%%%%%%%%%%%%%
\begin{proposition} \label{u_claim_9}
If there exists an overlap between perfect paths $p$ and $q$, then the underlying cycles $c_{p}$ and $c_{q}$ are equivalent.
\end{proposition}

\begin{proof}
Let $( p=p_{1}, \ldots, p_{n}, p_{n+1} = p_{1})$ and $(q = q_{1}, \ldots, q_{m}, q_{m+1} = q_{1})$ be the associated perfect path sequences.
Without loss of generality, we may assume that $n \geq 2$.
By assumption, we have $p_1 = q_0^\prime x_1$ and $q_1 = x_1 q_1^\prime$ for some $x_1 \in \B_{>0}$ and $q_0^\prime, q_1^\prime \in \B$ with $q_0^\prime x_1 q_1^\prime \not = 0$.
By Lemma \ref{u_claim_49}, there exist non-trivial paths $x_1, x_2,\ldots, x_{mn+1}$ and paths $ q_1^\prime, q_2^\prime, \ldots, q_{mn+1}^\prime$ such that $p_l =  q_{l-1}^\prime x_l$ and $q_l = x_l q_{l}^\prime$ for $1 \leq l \leq mn$.
In particular, we have 
\begin{align}
    p_1 p_2 \cdots p_{mn+1}q_{mn+1}^\prime 
    = (q_{0}^\prime x_1)(q_{1}^\prime x_2) \cdots (q_{mn+1}^\prime x_{mn+1}) q_{mn+1}^\prime 
    = q_0^\prime q_1 q_2  \cdots q_{mn+1}
\end{align}
where $p_{mn+1} = p_1$ and $q_{mn+1} = q_1$.
On the other hand, the associated underlying cycles $c_{p}$ and $c_{q}$ satisfy $p_1 \cdots p_n = c_{p}^h$ and $q_1 \cdots q_m = c_{q}^k$ for some $h$ and $k>0$, respectively.
Thus we have 
\begin{align}
    c_p^{mh} p_{mn+1}q_{mn+1}^\prime 
    = p_1 p_2 \cdots p_{mn+1}q_{mn+1}^\prime
    = q_0^\prime q_1 q_2  \cdots q_{mn+1}
    = q_0^\prime c_q^{nk} q_{mn+1}.
\end{align}
Since 
\begin{align}
    c_p^{mh} p_{mn+1}q_{mn+1}^\prime 
    &= q_0^\prime c_q^{nk} q_{mn+1} \\
    & = q_0^\prime c_q^{(n-1)k} c_q^k q_{mn+1} \\
    &= q_0^\prime c_q^{(n-1)k}  q_1 \cdots q_m x_{mn+1}q_{mn+1}^\prime \\
    &= q_0^\prime c_q^{(n-1)k} q_1 \cdots q_{m-1} x_m p_{mn+1}q_{mn+1}^\prime,
\end{align}
$c_p^{mh}$ contains $q_0^\prime c_q^{(n-1)k}$ as a left divisor, which implies that $c_q$ is a subpath of $c_p^{mh}$.
Now, write $c_p = a_1 a_2 \cdots a_{n_p}$, where $a_i \in Q_1$.
Since $c_q$ is a subpath of $c_p^{mh}$, we have  
\begin{align}
    c_q = a_i \cdots a_{n_p} c_p^{s} a_1 \cdots a_{i-1}
\end{align}
for some $1 \leq i \leq n_p$ and some $ 0 \leq s \leq mh$, where  $l(a_i \cdots a_{n_p}) < l(c_p) $ and $ l(a_1 \cdots a_{i-1}) < l(c_p)$.
Then we have $s=0$ since, by definition, $c_q$ contains no proper subcycle $c^\prime$ such that $(c^\prime)^t = c_q$ for some $t \geq 1$.  
Therefore, we can conclude that the cycles $c_p = a_1  \cdots a_{n_p}$ and $c_q = a_i \cdots a_{n_p} a_1 \cdots a_{i-1}$ are equivalent.
\end{proof}
%%%%%%%%%%%%%%%%%%%%%%%%%%%%%%%%%%%%
%           Statement ↑
%%%%%%%%%%%%%%%%%%%%%%%%%%%%%%%%%%%%

We then obtain the following corollary.

%%%%%%%%%%%%%%%%%%%%%%%%%%%%%%%%%%%%
%           Statement ↓
%%%%%%%%%%%%%%%%%%%%%%%%%%%%%%%%%%%%
\begin{corollary} \label{u_claim_19}
Let $p$ and $q$ be perfect paths with underlying cycles $c_p$ and $c_q$, respectively.
If $c_p \not= c_q$ in $\C(\L)$, then $\sHom_\L(p\L, q\L) = 0$.
\end{corollary}
\begin{proof}
It is a consequence of Propositions \ref{Chen-Shen-Zhou_2018_page1130} and \ref{u_claim_9}.
\end{proof}
%%%%%%%%%%%%%%%%%%%%%%%%%%%%%%%%%%%%
%           Statement ↑
%%%%%%%%%%%%%%%%%%%%%%%%%%%%%%%%%%%%

For any perfect paths $p$ and $q$, we know that $p\L \cap \L q/p \L q \cong \sHom_\L(p\L, q\L)$ as $K$-vector spaces.
If the quotient $p\L \cap \L q/p \L q$ is non-zero,  it has a $K$-basis consisting of non-zero paths. 
The following proposition shows that every basis element of such a $K$-basis is perfect.

%%%%%%%%%%%%%%%%%%%%%%%%%%%%%%%%%%%%
%           Statement ↓
%%%%%%%%%%%%%%%%%%%%%%%%%%%%%%%%%%%%
\begin{proposition} \label{u_claim_14}
Let $p$ and $q \in \mathbb{P}_\L$ satisfy that $p$ overlaps $q$.
Then any non-zero path in $p\L \cap \L q$ that does not belong to $p \L q$ is perfect.
\end{proposition}
\begin{proof} 
Let $( p = p_{1}, \ldots, p_{n}, p_{n+1} = p_{1})$ and $(q = q_{1}, \ldots, q_{m}, q_{m+1} = q_{1})$ be the associated perfect path sequences. 
By assumption, there exist paths $x_1 \in \B_{>0}$ and $q_0^\prime, q_1^\prime \in \B$ such that $p_1 = q_0^\prime x_1$, $q_1 = x_1 q_1^\prime$ and $q_0^\prime x_1 q_1^\prime \not =0$.
Then $q_0^\prime x_1 q_1^\prime$ lies in $ p\L \cap \L q \backslash p \L q$, and any non-zero path in $p\L \cap \L q \backslash p \L q$ can be expressed in this manner.
We now show that $q_0^\prime x_1 q_1^\prime = p_1 q^\prime_1$ is a perfect path.
Since there is nothing to prove for the case that $q_0^\prime$  or $q_1^\prime$ is trivial, suppose that $q_0^\prime$ and  $q_1^\prime \in \B_{>0}$.
By Lemma \ref{u_claim_49}, there exist families $\{x_l \in \B_{>0}\}_{l\geq 2}$  and $\{q_l^\prime \in \B\}_{l\geq 2}$ such that $p_{l} = q_{l-1}^\prime x_{l}$ and $q_{l} = x_{l} q_{l}^\prime$ for any $l \geq 2$ and such that $(p_l q_{l}^\prime, x_{l+1})$ and $(x_{l+1}, p_{l+2} q_{l+2}^\prime)$ are perfect for $l \geq 1$.
Now, since $q_m = q_{2im} = x_{2im} q_{2im}^\prime$ for all $i > 0$,  the fact that $l(q_m) < \infty$ implies that $x_{2jm} = x_{2km}$ for some $j < k$.
Hence the sequence
\begin{align}
    (x_{2jm}, \ p_{2jm+1} q_{2jm+1}^\prime, \  x_{2jm+2}, \  \ldots, \  p_{2km-1} q_{2km-1}^\prime,  \ x_{2km} = x_{2jm})
\end{align}
is a perfect path sequence.
Then $x_{2jm}$ and hence $p_{2jm-1} q_{2jm-1}^\prime, x_{2jm-2}, \ldots, x_2, p_1 q_{1}^\prime$ are all perfect paths.
\end{proof}
%%%%%%%%%%%%%%%%%%%%%%%%%%%%%%%%%%%%
%           Statement ↑
%%%%%%%%%%%%%%%%%%%%%%%%%%%%%%%%%%%%

%%%%%%%%%%%%%%%%%%%%%%%%%%%%%%%%%%%%
%           Subsection ↓
%%%%%%%%%%%%%%%%%%%%%%%%%%%%%%%%%%%%
\subsection{Two partial orders on the set of perfect paths} \label{section_1_2}
%%%%%%%%%%%%%%%%%%%%%%%%%%%%%%%%%%%%

In this subsection, we will introduce and study two partial orders $\preceq$ and $\leq$ on perfect paths.
Let us start by defining them.

%%%%%%%%%%%%%%%%%%%%%%%%%%%%%%%%%%%%
%           Statement ↓
%%%%%%%%%%%%%%%%%%%%%%%%%%%%%%%%%%%%
\begin{definition}{\rm 
For two perfect paths $p$ and $q$, we write $p \preceq q$ (resp.~$p \prec q$) if $p$ is a (resp.~proper) left divisor of $q$.
Dually, we write $p \leq q$ (resp.~$p < q$) if $q$ is a (resp.~proper) right divisor of $p$.
}\end{definition}
%%%%%%%%%%%%%%%%%%%%%%%%%%%%%%%%%%%%
%           Statement ↑
%%%%%%%%%%%%%%%%%%%%%%%%%%%%%%%%%%%%

Recall that $\mathbb{P}_{\L}$ stands for the set of perfect paths in $\L$.
It is straightforward that $(\mathbb{P}_{\L}, \preceq)$ and $(\mathbb{P}_{\L}, \leq)$ are partially ordered sets.
We remark that the partial order which can be found in \cite[the proof of Theorem 5.2.2]{Lu-Zhu_2021} is obtained by reversing the partial order $\leq$.

Let  $p$ and $q$ be perfect paths. 
Using the isomorphism $\theta^\prime_{q,p}: p\L \cap \L q \xrightarrow{\ \sim\ }  \Hom_{\L}(q\L, p\L)$ defined as in (\ref{eq_16}), one observes that $p \preceq q$ if and only if there exists an injective homomorphism $q\L \to p\L$ and that $p \leq q$ if and only if there exists a surjective homomorphism $q\L \to p\L$.

Let us begin with the following proposition.

%%%%%%%%%%%%%%%%%%%%%%%%%%%%%%%%%%%%
%           Statement ↓
%%%%%%%%%%%%%%%%%%%%%%%%%%%%%%%%%%%%
\begin{proposition} \label{u_claim_13}
The following statements hold for perfect paths $p$ and $q$.
\begin{enumerate}
    \item If $p \prec q$, then the non-trivial path $r$ in the factorization $q=pr$ is perfect. 
    \item If $p < q$, then the non-trivial path $r$ in the factorization $p=rq$ is perfect. 
\end{enumerate}
\end{proposition}
\begin{proof}  
Let $( p = p_{1}, \ldots, p_{n}, p_{n+1} = p_{1})$ and $(q = q_{1}, \ldots, q_{m}, q_{m+1} = q_{1})$ be the associated perfect path sequences. 
We first prove (1). 
Assume that there exists a non-trivial path $q_1^\prime=r$ such that $q_1 = p_1 q_1^\prime$.
It follows from Lemma \ref{u_claim_49} that there exist non-zero non-trivial paths $q_3^\prime, \ldots, q_{2mn+1}^\prime$ such that  $(q_{2i+1}^\prime, q_{2i+2}p_{2i+3})$ and $(q_{2i+2}p_{2i+3}, q_{2i+3}^\prime)$ are perfect pairs for $i \geq 0$.
Since $p_{1} q_1^\prime = q_1 = q_{2mn+1} = p_{2mn+1} q^\prime_{2mn+1} = p_{1} q^\prime_{2mn+1}$, we have $q_1^\prime = q^\prime_{2mn+1}$.
Thus we obtain the following perfect path sequence:
\begin{align}
    (q_1^\prime,  q_2 p_3,   q_3^\prime, \ldots,   q_{2mn}p_{2mn+1},   q^\prime_{2mn+1} = q_1^\prime),
\end{align}
which shows that $r = q_1^\prime$ is a perfect path.

It remains to prove (2). 
Assume that $p_1 = q_0^\prime q_1$ for some $ q_0^\prime =r \in \B_{>0}$.
Since $(p_n q_0^\prime) q_1 = p_n p_1 \in \F$, together with the fact that $(q_m, q_1)$ is perfect, we see that $q_m = p_n q_0^\prime$, so that $q_0^\prime$ is perfect by (1).
\end{proof}
%%%%%%%%%%%%%%%%%%%%%%%%%%%%%%%%%%%%
%           Statement ↑
%%%%%%%%%%%%%%%%%%%%%%%%%%%%%%%%%%%%

We then have the following corollary.

%%%%%%%%%%%%%%%%%%%%%%%%%%%%%%%%%%%%
%           Statement ↓
%%%%%%%%%%%%%%%%%%%%%%%%%%%%%%%%%%%%
\begin{corollary} \label{u_claim_15}
Let $(p_1, \ldots, p_n, p_{n+1}=p_1)$ and $(q_1, \ldots, q_n, q_{m+1}=q_1)$ be perfect path sequences. 
\begin{enumerate}
        
    \item If $p_1 \prec q_1$, then we have the following short exact sequence of Gorenstein-projective $\L$-modules:
    \[
    0 \rightarrow q_1 \L \xrightarrow{\,\iota\,} p_1 \L \xrightarrow{\,\pi\,} q_m p_1\L \rightarrow 0,
    \] 
    where $\iota$ is the canonical inclusion, and $\pi$ is the left multiplication by $q_m$.
        
    \item If $p_1 < q_1$, then we have the following short exact sequence of Gorenstein-projective $\L$-modules:
    \[
    0 \rightarrow q_1 p_2 \L \xrightarrow{\,\iota\,} q_1 \L \xrightarrow{\,\pi\,}  p_1\L \rightarrow 0,
    \]
    where $\iota$ is the canonical inclusion, and $\pi$ is the left multiplication by the non-trivial path $q_0^\prime$ in the factorization  $p_1 = q_0^\prime q_1$.
        
\end{enumerate}
\end{corollary}

\begin{proof}
We only prove (1); the proof of (2) is dual.
Suppose that $q_1 = p_1 q_1^\prime$ for some non-trivial path $q_1^\prime$.
Then $q_1^\prime$ is a perfect path by Proposition \ref{u_claim_13}, so that $q_m p_1$ is also perfect since we know that  $(q_m p_1, q_1^\prime)$ is a perfect pair.
Now, we claim that $\Ker \pi = q_1 \L$.
The inclusion $(\supseteq)$ is clear.
So take a path  $p_1 x $ such that $q_m(p_1 x) = 0$ in $\L$. 
Using the fact that $(q_m p_1, q_1^\prime)$ is perfect, one concludes that the path $x$ contains $ q_1^\prime$ as a left divisor, which shows that $p_1 x $ belongs to $p_1 q_1^\prime \L = q_1 \L$.
\end{proof}
%%%%%%%%%%%%%%%%%%%%%%%%%%%%%%%%%%%%
%           Statement ↑
%%%%%%%%%%%%%%%%%%%%%%%%%%%%%%%%%%%%

The vertex sets of the Hasse quivers  $H(\mathbb{P}_\L, \preceq)$ and $H(\mathbb{P}_\L, \leq)$  are both the finite set $\mathbb{P}_\L$. 
See Remark \ref{remark_6} for the finiteness of $\mathbb{P}_\L$.
Moreover, there exists an arrow from $q$ to $p$ in the Hasse quiver $H(\mathbb{P}_\L, \preceq)$ if and only if $p \prec q$ and there exists no perfect path $r$ such that  $p \prec r \prec q$; we draw arrows in the Hasse quiver $H(\mathbb{P}_\L, \leq)$ in a similar way.

%%%%%%%%%%%%%%%%%%%%%%%%%%%%%%%%%%%%
%           Statement ↓
%%%%%%%%%%%%%%%%%%%%%%%%%%%%%%%%%%%%
\begin{proposition} \label{u_claim_4}
The Hasse quivers $H(\mathbb{P}_\L, \preceq)$ and $H(\mathbb{P}_\L, \leq)$ have the form of a disjoint union of linearly oriented quivers of type $\A$. 
\end{proposition}    
\begin{proof} 
We show the statement for $H(\mathbb{P}_\L, \preceq)$; the proof for $H(\mathbb{P}_\L, \leq)$ is dual.
We first claim that there exists no subquiver of the form $q \to p \leftarrow r$ or $q \leftarrow p \to r$ in $H(\mathbb{P}_\L, \preceq)$.
Assume that there exist two arrows $q \to p$ and $r \to p$ in $H(\mathbb{P}_\L, \preceq)$, where $p, q$ and $r$ are perfect paths.
Let $(q, q_2)$ and $(r, r_2)$ be the associated perfect pairs.
Then there exist non-trivial paths $q^\prime$ and $r^\prime$ such that $q = p q^\prime$ and $r = p r^\prime$.
Since $p(q^\prime q_2) = q q_2 \in \F$ and $p(r^\prime r_2) = r r_2 \in \F$, it follows that $q^\prime q_2 = r^\prime r_2$.
Without loss of generality, we may assume that $r^\prime$ is a left divisor of $q^\prime$, that is, $q^\prime = r^\prime x$ for some $x \in \B$.
Then we have $q = p q^\prime = p r^\prime x = rx$, which implies that $r \preceq q$.
Thus $q$ must coincide with $r$.
Then, since $H(\mathbb{P}_\L, \preceq)$ contains by definition no multiple arrows, the two arrows $q \to p$ and $r \to p$ must be equal. 
Similarly, one can show that $H(\mathbb{P}_\L, \preceq)$ contains no subquiver of the form $q \leftarrow p \to r$.
Now, it follows from the claim that there exist no parallel paths between any two vertices in $H(\mathbb{P}_\L, \preceq)$.
Further, since, by definition, there exists no cycle in $H(\mathbb{P}_\L, \preceq)$, we deduce that $H(\mathbb{P}_\L, \preceq)$ is a disjoint union of connected trees.
Then these connected trees are all linearly oriented quivers of type $\A$ precisely when they contain no subquiver of the form $q \to p \leftarrow r$ or $q \leftarrow p \to r$, which is what we have proved.
\end{proof}
%%%%%%%%%%%%%%%%%%%%%%%%%%%%%%%%%%%%
%           Statement ↑
%%%%%%%%%%%%%%%%%%%%%%%%%%%%%%%%%%%%

%%%%%%%%%%%%%%%%%%%%%%%%%%%%%%%%%%%%
%           Statement ↓
%%%%%%%%%%%%%%%%%%%%%%%%%%%%%%%%%%%%
\begin{definition}{\rm \label{def_2}
A perfect path $p$ is called {\it elementary} if $p$ is a source in  $H(\mathbb{P}_\L, \preceq)$. Dually, the perfect path $p$ is called {\it co-elementary} if $p$ is a sink in  $H(\mathbb{P}_\L, \preceq)$.
We refer to a module that is isomorphic to the cyclic module generated by an elementary perfect path as an {\it elementary Gorenstein-projective module}.
}\end{definition}
%%%%%%%%%%%%%%%%%%%%%%%%%%%%%%%%%%%%
%           Statement ↑
%%%%%%%%%%%%%%%%%%%%%%%%%%%%%%%%%%%%

%%%%%%%%%%%%%%%%%%%%%%%%%%%%%%%%%%%%
%           Statement ↓
%%%%%%%%%%%%%%%%%%%%%%%%%%%%%%%%%%%%
\begin{remark}{\rm 
It follows from Proposition \ref{u_claim_16} below that elementary Gorenstein-projective modules are the same as those defined in \cite{Ringel_2013}.
}\end{remark}
%%%%%%%%%%%%%%%%%%%%%%%%%%%%%%%%%%%%
%           Statement ↑
%%%%%%%%%%%%%%%%%%%%%%%%%%%%%%%%%%%%

We denote by $\mathbb{E}_\L$ (resp.~$\mathbb{E}^{\rm co}_\L$) the set of elementary (resp.~co-elementary) perfect paths. 
By Proposition \ref{u_claim_4}, 
$\mathbb{E}_{\L}$ is in a bijection with $\mathbb{E}_{\L}^{\rm co}$.

%%%%%%%%%%%%%%%%%%%%%%%%%%%%%%%%%%%%
%           Statement ↓
%%%%%%%%%%%%%%%%%%%%%%%%%%%%%%%%%%%%
\begin{example}{\rm \label{ex_3}
Let $\L$ be the monomial algebra given by the following quiver with relations:
\[
\xymatrix{
1 \ar[r]^-{a_1}  & 2 \ar[r]^-{b_2} \ar[d]^-{a_2} & 4  \ar@<0.6ex>[r]^-{a_4} & 5 \ar@<0.6ex>[l]^-{a_5}  \\
  & 3 \ar[lu]^-{a_3} & 
}
\qquad 
a_{12312312} = a_{3123123} = a_{45454545} = 0,
\]
where $a_{12312312}$ denotes the path $a_1 a_2 a_3 a_1 a_2 a_3 a_1 a_2$ for example. 
Then a direct calculation shows that up to cyclic permutation, 
$(a_{12}, a_{312312}, a_3, a_{123123}, a_{12})$, $(a_{12312},$ $a_{312}, a_{3123},$ $a_{123}, a_{12312})$, $(a_{45}, a_{454545}, a_{45})$ and $(a_{4545}, a_{4545})$ 
are the minimal perfect path sequences.
Hence $\C(\L) =\{a_{123},\, a_{45}\}$, and the Hasse quivers $H(\mathbb{P}_{\L}, \preceq)$ and $H(\mathbb{P}_{\L}, \leq)$ are given as follows:
%%%%%%%%%%%%%%%%%%%%%%%%%%%%%
% diagram ↓
%%%%%%%%%%%%%%%%%%%%%%%%%%%%%
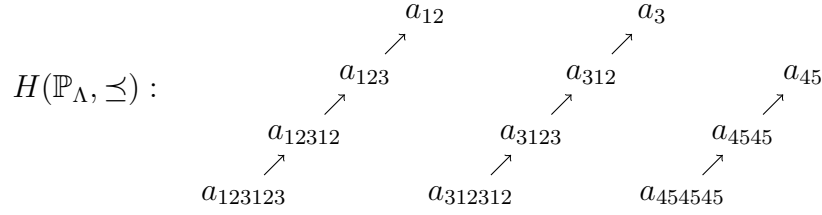
\[
\begin{tikzpicture}[scale=1]
\node (0) at (-3.7,-.2) {$H(\mathbb{P}_{\L}, \preceq):$};
\node (1) at (-1.6,-1.6) {$a_{123123}$};
\node (2) at (-.8,-.8) {$a_{12312}$};
\node (3) at (0,0) {$a_{123}$};
\node (4) at (.8,.8) {$a_{12}$};
\node (5) at (1.4,-1.6) {$a_{312312}$};
\node (6) at (2.2,-.8) {$a_{3123}$};
\node (7) at (3,0) {$a_{312}$};
\node (8) at (3.8,.8) {$a_{3}$};
\node (9) at (4.2,-1.6) {$a_{454545}$};
\node (10) at (5.0,-.8) {$a_{4545}$};
\node (11) at (5.8,0) {$a_{45}$};
\draw[->] (1) -- (2);
\draw[->] (2) -- (3);
\draw[->] (3) -- (4);
\draw[->] (5) -- (6);
\draw[->] (6) -- (7);
\draw[->] (7) -- (8);
\draw[->] (9) -- (10);
\draw[->] (10) -- (11);
\end{tikzpicture}
\]
%%%%%
and
%%%%%
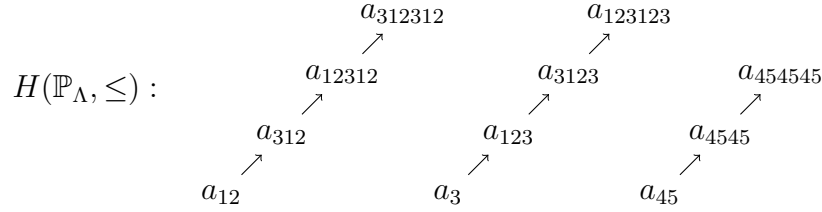
\[
\begin{tikzpicture}[scale=1]
\node (0) at (-3.4,-.2) {$H(\mathbb{P}_{\L}, \leq):$};
\node (1) at (-1.6,-1.6) {$a_{12}$};
\node (2) at (-.8,-.8) {$a_{312}$};
\node (3) at (0,0) {$a_{12312}$};
\node (4) at (.8,.8) {$a_{312312}$};
\node (5) at (1.4,-1.6) {$a_{3}$};
\node (6) at (2.2,-.8) {$a_{123}$};
\node (7) at (3.0,0) {$a_{3123}$};
\node (8) at (3.8,.8) {$a_{123123}$};
\node (9) at (4.2,-1.6) {$a_{45}$};
\node (10) at (5.0,-.8) {$a_{4545}$};
\node (11) at (5.8,0) {$a_{454545}$};
\draw[->] (1) -- (2);
\draw[->] (2) -- (3);
\draw[->] (3) -- (4);
\draw[->] (5) -- (6);
\draw[->] (6) -- (7);
\draw[->] (7) -- (8);
\draw[->] (9) -- (10);
\draw[->] (10) -- (11);
\end{tikzpicture}
\]
%%%%%%%%%%%%%%%%%%%%%%%%%%%%%
% diagram ↑
%%%%%%%%%%%%%%%%%%%%%%%%%%%%%
In particular, $\mathbb{E}_\L = \{a_{123123}, a_{312312}, a_{454545}\}$ and $\mathbb{E}^{\rm co}_\L = \{a_{12}, a_{3}, a_{45}\}$.
}\end{example}
%%%%%%%%%%%%%%%%%%%%%%%%%%%%%%%%%%%%
%           Statement ↑
%%%%%%%%%%%%%%%%%%%%%%%%%%%%%%%%%%%%

In the following, we present properties of elementary and co-elementary paths.

%%%%%%%%%%%%%%%%%%%%%%%%%%%%%%%%%%%%
%           Statement ↓
%%%%%%%%%%%%%%%%%%%%%%%%%%%%%%%%%%%%
\begin{proposition} \label{u_claim_16}
The following statements hold for a perfect path $p$. 
\begin{enumerate}
    \item $p$ is elementary if and only if $p$ is a sink in $H(\mathbb{P}_\L, \leq)$.
    \item $p$ is co-elementary if and only if $p$ is a source in $H(\mathbb{P}_\L, \leq)$.
    \item If $(p_n, p)$ and $(p, p_2)$ are perfect pairs, then the following conditions are equivalent:
        \begin{enumerate}
            \item $p$ is elementary.
            \item $p_2$ is co-elementary.
            \item $p_n$ is co-elementary.
        \end{enumerate}
\end{enumerate}
\end{proposition}
\begin{proof} 
We prove (1); and (2) is proved dually.
For the\,\lq\lq only if\,\rq\rq\,part, assume that $p$ is not a sink in  $H(\mathbb{P}_\L, \leq)$, that is, $q < p$ for some $q \in \mathbb{P}_{\L}$.
Let $(q, q_2)$ be the associated perfect pair.
Then, by Corollary \ref{u_claim_15} (2), we have the canonical inclusion $pq_2 \L \hookrightarrow p\L$ with $pq_2$ perfect, which means that $p \prec pq_2$. This implies that $p$ is not elementary.  
For the\,\lq\lq if\,\rq\rq\,part, suppose that $p$ is not elementary, or equivalently, $p \prec q$ for some $q \in \mathbb{P}_\L$.
Then Corollary \ref{u_claim_15} (1) yields an epimorphism  $p \L \twoheadrightarrow q_m p\L$ with $q_m p$ perfect, where $(q_m, q)$ is the associated perfect pair. 
Therefore, we have $q_m p < p$, which means that $p$ is not a sink in  $H(\mathbb{P}_\L, \leq)$.

For (3), we will show the equivalence of (a) and (b); the equivalence of (a) and (c) is obtained similarly. 
To show that (a) implies (b), assume that $p_2$ is not co-elementary.
Then there exists a perfect path $q$ such that $q \prec p_2$.
In particular, $q$ is a proper left divisor of $p_2$.
Applying Proposition \ref{u_claim_13} to $q_1 := p_2$ and $p_1 := q$, we see that $q_{2mn}$ and $p_{2mn+1}$ in the proof of the proposition are nothing but $p$ and  $q$, respectively. 
Hence $pq$ is a perfect path, which implies that $p \prec pq$. Therefore, $p$ is not elementary.
Conversely, suppose that $p$ is not elementary, that is, $p \prec q$ for some $q \in \mathbb{P}_{\L}$.
It then follows from Proposition \ref{u_claim_13} again that the path $r$ in the factorization $q=pr$ is perfect.
Further, since the underlying cycles $c_q$ and $c_p$ are equivalent, $r$ has to be a left divisor of the perfect path $p_2$ next to $p$.
Consequently, we have $r \prec p_2$, which shows that $p_2$ is not co-elementary.
\end{proof}
%%%%%%%%%%%%%%%%%%%%%%%%%%%%%%%%%%%%
%           Statement ↑
%%%%%%%%%%%%%%%%%%%%%%%%%%%%%%%%%%%%

For each $c \in \C(\L)$, we define 
\begin{align}
    X_c := \{ p \in \mathbb{E}_{\L} \mid c_p = c \mbox{ in } \C(\L) \}
    \quad \mbox{and} \quad 
    Y_c := \{ r \in \mathbb{E}_{\L}^{\rm co} \mid c_r = c \mbox{ in } \C(\L) \}.
\end{align}
Observe that 
\begin{align} \label{eq_13}
    \mathbb{E}_{\L}=\bigcup_{c \in \C(\L)} X_c \quad \mbox{and} \quad  \mathbb{E}_{\L}^{\rm co} = \bigcup_{c \in \C(\L)} Y_c.
\end{align}
For each $p \in X_c$, the perfect path $r$ in the perfect pair $(p, r)$ is co-elementary by Proposition \ref{u_claim_16} and belongs to $Y_c$ since $c_r = c_p = c$ in $\C(\L)$.
Thus we get a map $\varphi_c: X_c \to Y_c$ which sends $p \in  X_c$ to $r \in Y_c$.
It is easily verified that $\varphi_c: X_c \to Y_c$ is bijective.

%%%%%%%%%%%%%%%%%%%%%%%%%%%%%%%%%%%%
%           Subsection ↓
%%%%%%%%%%%%%%%%%%%%%%%%%%%%%%%%%%%%
\subsection{Main result} \label{section_1_3}
%%%%%%%%%%%%%%%%%%%%%%%%%%%%%%%%%%%%

The aim of this subsection is to establish a structure theorem for perfect paths.
For this, we need the following two lemmas.

%%%%%%%%%%%%%%%%%%%%%%%%%%%%%%%%%%%%
%           Statement ↓
%%%%%%%%%%%%%%%%%%%%%%%%%%%%%%%%%%%%
\begin{lemma} \label{u_claim_28}
There exists no overlap between any two co-elementary paths.
\end{lemma}

\begin{proof} 
Assume to the contrary that there exists an overlap between two co-elementary paths $p$ and $q$. 
Then there exist paths $p^\prime, q^\prime  \in \B$ and $x \in \B_{>0}$ such that 
$p = p^\prime x, q = x q^\prime$ and $p^\prime x q^\prime \not= 0$.
By Propositions \ref{u_claim_14} and \ref{u_claim_13}, all the non-trivial paths among the three paths $p^\prime, q^\prime$ and $x$ are perfect. 
In particular, $x$ is always perfect. 
This contradicts the fact that $p$ is co-elementary.
\end{proof}
%%%%%%%%%%%%%%%%%%%%%%%%%%%%%%%%%%%%
%           Statement ↑
%%%%%%%%%%%%%%%%%%%%%%%%%%%%%%%%%%%%

%%%%%%%%%%%%%%%%%%%%%%%%%%%%%%%%%%%%
%           Statement ↓
%%%%%%%%%%%%%%%%%%%%%%%%%%%%%%%%%%%%
\begin{lemma} \label{u_claim_17}
Let $p$ and $q$ be perfect paths.
\begin{enumerate}
    \item Assume that $p\prec q$. Then there exists an arrow $q \rightarrow p$ in $H(\mathbb{P}_\L, \preceq)$ if and only if the perfect path $r$ in $q = pr$ is co-elementary. In this case, the cokernel of the induced canonical inclusion $q\L \hookrightarrow p\L$ is elementary Gorenstein-projective.

    \item Assume that $p< q$. Then there exists an arrow $q \rightarrow p$ in $H(\mathbb{P}_\L, \leq)$ if and only if the perfect path $r$ in $p = rq$ is co-elementary. In this case, the kernel of the induced canonical surjection $q\L \twoheadrightarrow p\L$ is elementary Gorenstein-projective.
\end{enumerate}
\end{lemma}
\begin{proof} 
We only show (1); (2) is proved dually.
We first show the\,\lq\lq only if\,\rq\rq\,part. 
Assume for a contradiction that $q= pr$ for some perfect path $r$ which is not co-elementary. 
By definition, there exists a perfect path $s$ such that $s \prec r$, that is, $r = s r^\prime$ for some $r^\prime \in \B_{>0}$.
This implies that $p\prec ps \prec q$, a contradiction to the existence of the arrow $p \to q$.
Conversely, assume that $r$ is co-elementary.
If there exists a perfect path $q^\prime$ such that $p \prec q^\prime \prec q$, then $q^\prime = p r^{\prime}$ for some $r^{\prime} \in \mathbb{P}_\L$. This implies that $r^{\prime}$ is a left divisor of $r$, that is, $r^{\prime} \prec r$.
This is a contradiction.

For the last statement, let $(q_m, q)$ be the associated perfect pair.
Then it follows from Corollary \ref{u_claim_15} and its proof that the cokernel of the inclusion $q\L \hookrightarrow p\L$ is given by $q_m p\L$, and the pair $(q_m p, r)$ is perfect.
As $r$ is co-elementary, Proposition \ref{u_claim_16} shows that $q_m p$ is elementary.
\end{proof}
%%%%%%%%%%%%%%%%%%%%%%%%%%%%%%%%%%%%
%           Statement ↑
%%%%%%%%%%%%%%%%%%%%%%%%%%%%%%%%%%%%

We are now ready to show the following main result of this section.

%%%%%%%%%%%%%%%%%%%%%%%%%%%%%%%%%%%%
%           Statement ↓
%%%%%%%%%%%%%%%%%%%%%%%%%%%%%%%%%%%%
\begin{theorem} \label{u_claim_25}
Let $\L$ be a monomial algebra. 
Then we have the following assertions.
\begin{enumerate}
    
    \item For a perfect path $p$, there exist finitely many co-elementary paths $r_1, \ldots, r_n$ such that $p = r_1 \cdots r_n$. 
    Moreover, if $p = s_1 \cdots s_m$ with $s_i$ co-elementary, then $m = n$ and $r_i = s_i$ for all $i$.

    \item 
    Let $r_1, \ldots, r_n$ be co-elementary paths such that $c_{r_i} = c_{r_j}$ in $\C(\L)$ for any $i, j$.
    If the product $r_1 \cdots r_n$ is non-zero in $\L$, then it is perfect. 
    
\end{enumerate}
\end{theorem}

\begin{proof} 
Let $p\in \mathbb{P}_{\L}$.
It follows from Lemma \ref{u_claim_17} that $p$ can be written as the product $r_1 r_2 \cdots r_n$ of co-elementary paths $r_1, r_2\ldots, r_n$.
On the other hand, by definition, any co-elementary path contains no perfect left divisor. 
Consequently, if $p = s_1 s_2 \cdots s_m$ with $s_i \in \mathbb{E}^{\rm co}_{\L}$, then we must have $r_1 = s_1$,  so that $r_2 \cdots r_n = s_2 \cdots s_m$.
Repeating this argument, we obtain the desired uniqueness.

To show (2), we use an induction on $n$.
The case $n=1$ is trivial.
For $n > 1$, take a non-zero path $r_1 r_2 \cdots r_n$, where each $r_i$ is co-elementary, and assume that the right divisor $r = r_2 \cdots r_n$ is perfect.
Let $(r_1, s)$ be the associated perfect pair.
Together with the fact that $s$ is elementary, the fact that $c_r = c_{r_1}$ in $\C(\L)$ implies that $r$ is a left divisor of $s$.
Then, applying Lemma \ref{u_claim_49} (2) to $p_1 := r, q_1:= s$ and $q_m:=r_1$, we deduce that $r_1 r $ is perfect (cf.~the proof of Proposition \ref{u_claim_13}). 
\end{proof}
%%%%%%%%%%%%%%%%%%%%%%%%%%%%%%%%%%%%
%           Statement ↑
%%%%%%%%%%%%%%%%%%%%%%%%%%%%%%%%%%%%

The theorem enables us to obtain the following observation, which asserts that any underlying cycle is the concatenation of co-elementary paths.

%%%%%%%%%%%%%%%%%%%%%%%%%%%%%%%%%%%%
%           Statement ↓
%%%%%%%%%%%%%%%%%%%%%%%%%%%%%%%%%%%%
\begin{proposition-definition} \label{u_claim_26}
Let $c \in \C(\L)$. Then there uniquely exist finitely many co-elementary paths $r_1,  \ldots, r_n$ such that $c= r_1  \cdots r_n$.
We denote $|c|:=n$.
\end{proposition-definition}
\begin{proof} 
First, we claim that for any co-elementary paths $r$ and $s$ with $c_r = c_s$ in $\C(\L)$, if $r$ contains $s$ as a subpath, then $r=s$.
Assume that $s$ is a proper subpath of $r$, say $r=xsy$, where $x, y \in \B$.
We may assume that both $x$ and $y$ are non-trivial and not perfect.
Since $c_r$ and $c_s$ are equivalent by assumption, it follows that $y$ is a proper left divisor of $s_2$, where $s_2$ is the perfect path appearing in the perfect pair $(s, s_2)$: 
%%%%%%%%%%%%%%%%%%%%%%%%%%%%%
% diagram ↓
%%%%%%%%%%%%%%%%%%%%%%%%%%%%%
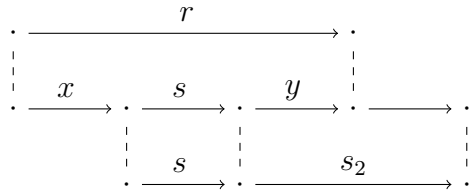
\[\begin{tikzpicture}[scale=1.0]
%%%%%%%%%%%%%
% paths
%%%%%%%%%%%%%
%
\node at (2.3,0.25) {$r$};
\node at (.7,-.75) {$x$};
\node at (2.2,-0.75) {$s$};
\node at (3.7,-0.75) {$y$};
\node at (2.2,-1.75) {$s$};
\node at (4.5,-1.75) {$s_2$};
%%%%%%%%%%%%%
% vertices
%%%%%%%%%%%%%
% a
\node (a1) at (0,0) {$\cdot$};
\node (a2) at (4.5,0) {$\cdot$};
% b
\node (b1) at (0,-1) {$\cdot$};
\node (b2) at (1.5,-1) {$\cdot$};
\node (b3) at (3,-1) {$\cdot$};
\node (b4) at (4.5,-1) {$\cdot$};
\node (b5) at (6,-1) {$\cdot$};
% c
\node (c1) at (1.5,-2) {$\cdot$};
\node (c2) at (3,-2) {$\cdot$};
\node (c3) at (6,-2) {$\cdot$};
%%%%%%%%%%%%%
% arrows
%%%%%%%%%%%%%
% a
\draw[->] (a1) -- (a2) ;
% b
\draw[->] (b1) -- (b2) ;
\draw[->] (b2) -- (b3) ;
\draw[->] (b3) -- (b4) ;
\draw[->] (b4) -- (b5) ;
% c
\draw[->] (c1) -- (c2) ;
\draw[->] (c2) -- (c3) ;
% vertical dashed lines
\draw[dashed] (a1) -- (b1) ;
\draw[dashed] (b2) -- (c1)  ;
\draw[dashed] (b3) -- (c2) ;
\draw[dashed] (a2) -- (b4) ;
\draw[dashed] (b5) -- (c3) ;
\end{tikzpicture}
\]
%%%%%%%%%%%%%%%%%%%%%%%%%%%%%
% diagram ↑
%%%%%%%%%%%%%%%%%%%%%%%%%%%%%
Since $s_2$ is the concatenation of co-elementary paths by Theorem \ref{u_claim_25}, it follows that the co-elementary path $r$ overlaps some co-elementary subpath of $s_2$, a contradiction.

Next, we claim that $l(r) \leq l(c_r)$ for any co-elementary path $r$.
Let $(p, r)$ be the associated perfect pair.
Then we have $p = c_p^l y$ for some $l \geq 0$ and some left divisor $y$ of $c_p$. 
Since $c_p y = y c_r$, we have  
\begin{align}
    c_p p =  c_p^{l+1} y = c_p^{l} y c_r = \cdots =  y c_r^{l+1} = p c_r. 
\end{align}
There are two cases to consider.
First, we deal with the case where $p$ is a left divisor of $c_p$.
Letting $c_p = p p^\prime$ with $p^\prime \in \B$, we have $p p^\prime p = c_p p = p c_r$, so that $p^\prime p = c_r$.
Therefore, $r$ is a left divisor of $c_r$; otherwise, the co-elementary path $r$ contains the co-elementary right divisor $s$ of $p$, which satisfies $c_s = c_p =c_r$ in $\C(\L)$. 
This contradicts the first claim.
Consider now the case where $c_p$ is a proper left divisor of $p$.
Then $p c_r = 0$ in $\L$.
Indeed, if this is not the case, then the non-zero path $p c_r$ belongs to $p\L\cap\L p \backslash p\L p$ and hence is perfect by Proposition \ref{u_claim_14}. 
In particular, we have $p c_r \prec p$, a contradiction. 
Therefore, the formula $R(p) =\{r\}$ concludes that $r$ is a left divisor of $c_r$.

Now, let $c \in \C(\L)$.
Then we have $c^l = p_1 p_2 \cdots p_m$  for some perfect path sequence $( p_{1}, p_{2}, \ldots, p_{m}, p_{m+1} = p_{1})$ and some positive integer $l$.
By Theorem \ref{u_claim_25}, there exist finitely many co-elementary paths $r_1, r_2, \ldots, r_t$ such that $c^l = p_1 p_2 \cdots p_m = r_1 r_2 \cdots r_t$.
First, the second claim shows that $c = r_1 r_2 \cdots r_{s-1} y$ for some $1 \leq s \leq t$ and some non-trivial left divisor $y$ of $r_s$.
Then, due to the first claim, the subpath $y^\prime$ in $r_s = y y^\prime$ becomes a left divisor of $r_1$.
But then Lemma \ref{u_claim_28} concludes that $y^\prime$ is a trivial path.

Finally, as in the proof of Theorem \ref{u_claim_25}, using the definition of co-elementary paths, one can show the uniqueness of co-elementary subpaths of $c$.
\end{proof}
%%%%%%%%%%%%%%%%%%%%%%%%%%%%%%%%%%%%
%           Statement ↑
%%%%%%%%%%%%%%%%%%%%%%%%%%%%%%%%%%%%

Recall that elementary Gorenstein-projective modules are defined to be modules isomorphic to the cyclic module generated by an elementary perfect path.
Let $\mathcal{E}_{\L}$ denote a complete set of pairwise non-isomorphic elementary Gorenstein-projective $\L$-modules.
The following is an easy consequence of Lemma \ref{u_claim_17}.

%%%%%%%%%%%%%%%%%%%%%%%%%%%%%%%%%%%%
%           Statement ↓
%%%%%%%%%%%%%%%%%%%%%%%%%%%%%%%%%%%%
\begin{proposition} \label{u_claim_20}
Let $M$ be an indecomposable non-projective Gorenstein-projective $\L$-module.
Then there exists an ascending sequence of indecomposable non-projective Gorenstein-projective $\L$-submodules of $M$
\begin{align}
    0 = M_0 \subset M_1 \subset M_2 \subset \cdots \subset M_{n-1} \subset M_{n} = M
\end{align}
such that $M_{i}/M_{i-1}$ is in $\mathcal{E}_{\L}$ for all $1 \leq i \leq n$.
\end{proposition}
\begin{proof}
It suffices to show the statement for the modules of the form $p\L$ with $p \in \mathbb{P}_{\L}$.
Take the full subquiver of the Hasse quiver $H(\mathbb{P}_\L, \preceq)$ formed by perfect paths $q$ with $p \preceq q$
\begin{align}
     q_1 \to q_2 \to \cdots \to q_{n-1} \to q_n = p.
\end{align}
Then it follows from Lemma \ref{u_claim_17} that the induced ascending chain of Gorenstein-projective $\L$-submodules of $p\L$
\begin{align}
    0 = q_0 \L \subset q_1 \L \subset q_2 \L \subset \cdots \subset q_{n-1} \L \subset q_{n} \L = p\L
\end{align}  
satisfies $q_i\L/q_{i-1}\L$ is elementary for all $ 1 \leq i \leq n$.
This finishes the proof.
\end{proof}
%%%%%%%%%%%%%%%%%%%%%%%%%%%%%%%%%%%%
%           Statement ↑
%%%%%%%%%%%%%%%%%%%%%%%%%%%%%%%%%%%%

%%%%%%%%%%%%%%%%%%%%%%%%%%%%%%%%%%%%
%           Statement ↓
%%%%%%%%%%%%%%%%%%%%%%%%%%%%%%%%%%%%
\begin{remark}{\rm 
As the following example shows, $\mathcal{E}_\L$ is not a set of orthogonal bricks in general.
Namely, the following conditions do not hold in general: $\End_\L(E)$ is a division ring; and $\Hom_\L(E, E^\prime) = 0$ for $E \not = E^\prime \in \mathcal{E}_\L$. 
}\end{remark}
%%%%%%%%%%%%%%%%%%%%%%%%%%%%%%%%%%%%
%           Statement ↑
%%%%%%%%%%%%%%%%%%%%%%%%%%%%%%%%%%%%

%%%%%%%%%%%%%%%%%%%%%%%%%%%%%%%%%%%%
%           Statement ↓
%%%%%%%%%%%%%%%%%%%%%%%%%%%%%%%%%%%%
\begin{example}{\rm \label{ex_1}
Let $\L$ be as in Example \ref{ex_3}.
As we saw there, $\mathbb{E}_{\L}=\{a_{123123}, $ $a_{312312},$ $a_{454545}\}$. 
Since $a_{123123} a_1 b_2 a_{454545}$ is non-zero in $\L$, we obtain that 
\begin{align}
    \Hom_\L(a_{454545}\L, a_{123123}\L) \cong a_{123123}\L \cap \L a_{454545} \not = 0.
\end{align}
}\end{example}
%%%%%%%%%%%%%%%%%%%%%%%%%%%%%%%%%%%%
%           Statement ↑
%%%%%%%%%%%%%%%%%%%%%%%%%%%%%%%%%%%%

%%%%%%%%%%%%%%%%%%%%%%%%%%%%%%%%%%%%
%           Section ↑
%%%%%%%%%%%%%%%%%%%%%%%%%%%%%%%%%%%%

%%%%%%%%%%%%%%%%%%%%%%%%%%%%%%%%%%%%
%           Section ↓
%%%%%%%%%%%%%%%%%%%%%%%%%%%%%%%%%%%%
\section{Realizing stable categories of graded Gorenstein-projective modules as derived categories} \label{section_3}
%%%%%%%%%%%%%%%%%%%%%%%%%%%%%%%%%%%%

In the rest of this paper, unless otherwise specified, we think of $\L$ as a positively graded algebra by setting the degree of each arrow to one.
In this section, we apply tilting theory to obtain a triangle equivalence between the stable category $\sGprojZ\L$ and the bounded derived category $\D^{\rm b}(\mod KQ)$ of a path algebra $KQ$, where $Q$ is a disjoint union of Dynkin quivers of type $\A$.
Moreover, we provide an Auslander-Reiten triangle $\tau C \to B \to C \to \Sigma \tau C$ in $\sGprojZ\L$ for any indecomposable object $C$.

%%%%%%%%%%%%%%%%%%%%%%%%%%%%%%%%%%%%
%           Subsection ↓
%%%%%%%%%%%%%%%%%%%%%%%%%%%%%%%%%%%%
\subsection{Graded Gorenstein-projective modules} \label{section_3_1}
%%%%%%%%%%%%%%%%%%%%%%%%%%%%%%%%%%%%

In this subsection, we aim to provide the main ingredients necessary for studying $\sGprojZ\L$ and to clarify the $K$-vector space structure of the Hom space between indecomposable objects in $\sGprojZ\L$.

Any non-zero path $p$ and primitive idempotent $e$ of $\L$ are homogeneous elements of degree $l(p)$ and zero, respectively.
In particular, the cyclic $\L$-modules $p\L$ and $e\L$ are both gradable. 
We will treat $p\L$ as a graded $\L$-module whose top is concentrated in degree $l(p)$. 
Hence $p\L = \bigoplus_{i \in \mathbb{Z}}p\L_{i}$ satisfies that $p\L_{i}$ is spanned by the non-zero paths of the form $px$ with $x \in \B_{i-l(p)}$ if such non-zero paths exist and otherwise $p\L_{i} = 0$. 
Moreover, we regard $e\L$ as a graded $\L$-module whose top is concentrated in degree zero.

The algebra $\L$ is always CM-finite and hence graded CM-finite (see Remark \ref{remark_6}).
Recall that $\ind\Gproj\L =  \left\{\,p\L \mid p\in \mathbb{P}_{\L}\right\} \cup \ind \proj\L.$

%%%%%%%%%%%%%%%%%%%%%%%%%%%%%%%%%%%%
%           Statement ↓
%%%%%%%%%%%%%%%%%%%%%%%%%%%%%%%%%%%%
\begin{lemma} \label{u_claim_48} 
We have the following statements.
\begin{enumerate}

    \item (\cite[Lemma 4.2.1]{Lu-Zhu_2021}). The forgetful functor $F_G : \GprojZ\L \rightarrow \Gproj\L$ is dense. 
    
    \item $\ind\GprojZ\L = \left\{\,p\L (i) \mid p\in \mathbb{P}_{\L}, i \in \mathbb{Z}\right\} \cup \ind \projZ\L$. 
    
    \item The following conditions are equivalent.
    \begin{enumerate}
        \item $\mathbb{P}_{\L}$ is empty.
        \item $\L$ is CM-free.
        \item $\L$ is graded CM-free.
    \end{enumerate} 
    
\end{enumerate}
\end{lemma}
\begin{proof} 
For (1), it suffices to show that the graded $\L$-module $p\L$ belongs to $\GprojZ\L$ for any $p \in \mathbb{P}_\L$.
Let $p$ be a perfect path.
Then we can construct a totally acyclic complex $P^\bullet$ associated to $p\L \in \Gproj\L$ by splicing the short exact sequences arising from an associated perfect path sequence $(p=p_1, \ldots, p_n, p_{n+1} = p_1)$.
But then it is easy to see that each of the short exact sequences is gradable (cf.~Lemma \ref{u_claim_24} below), so that $P^\bullet$ is also a totally acyclic complex of graded projective $\L$-modules.
This shows that $p\L \in \GprojZ\L$. 
The remaining statements follow from (1), Proposition \ref{u_claim_46} and Theorem \ref{Theorem4.1_Chen-Shen-Zhou_2018}.
\end{proof}
%%%%%%%%%%%%%%%%%%%%%%%%%%%%%%%%%%%%
%           Statement ↑
%%%%%%%%%%%%%%%%%%%%%%%%%%%%%%%%%%%%

The following lemma is an easy consequence of results in Section \ref{section_1}.
Note that the third statement is simply a graded version of Proposition \ref{u_claim_20}.

%%%%%%%%%%%%%%%%%%%%%%%%%%%%%%%%%%%%
%           Statement ↓
%%%%%%%%%%%%%%%%%%%%%%%%%%%%%%%%%%%%
\begin{lemma} \label{u_claim_24}
The following statements hold.
    \begin{enumerate}
    
        \item For a perfect pair $(p, q)$, we have the following short exact sequence in $\GprojZ\L:$
        \[
        0 \rightarrow q\L \xrightarrow{\,\iota\,} e_{t(p)}\L \xrightarrow{\,\pi\,} p\L(l(p)) \rightarrow 0,
        \]
        where $\iota$ is the canonical inclusion, and $\pi$ is the left multiplication by $p$.
        In particular, we have 
        \begin{align}
            \Sigma q \L = p\L(l(p))\quad \ \mbox{and} \quad \ \Sigma^{-1} p \L = q\L(-l(p))
        \end{align}
        in the stable category $\sGprojZ\L$ with suspension functor $\Sigma$.
    
        \item Let $(p_1, \ldots, p_n, p_{n+1}=p_1)$ and $(q_1, \ldots, q_m, q_{m+1}=q_1)$ be perfect path sequences. 
        \begin{enumerate}
        
            \item If $p_1 \prec q_1$, then we have the following short exact sequence in $\GprojZ\L:$
        \[
        0 \rightarrow q_1 \L \xrightarrow{\,\iota\,} p_1 \L \xrightarrow{\,\pi\,} q_m p_1\L(l(q_m)) \rightarrow 0,
        \]
        where $\iota$ is the canonical inclusion, and $\pi$ is the left multiplication by $q_m$.
        
        \item If $p_1 < q_1 $, then we have the following short exact sequence in $\GprojZ\L:$
        \[
        0 \rightarrow q_1 p_2 \L \xrightarrow{\,\iota\,} q_1 \L \xrightarrow{\,\pi\,}  p_1\L(l(q^\prime_0)) \rightarrow 0,
        \]
        where $\iota$ is the canonical inclusion, and $\pi$ is the left multiplication by the non-trivial path  $q^\prime_0$ in the factorization $p_1 = q^\prime_0 q_1$.
        
        \end{enumerate}

    \item Let $\mathcal{E}_\L^{\mathbb{Z}} = \{\,E(i)\mid \mbox{$E \in \mathcal{E}_\L, i \in \mathbb{Z}$}\,\}$. 
    Then any indecomposable non-projective object of $\GprojZ\L$ admits an ascending sequence of indecomposable non-projective subobjects with factors in $\mathcal{E}_\L^{\mathbb{Z}}$.
    
    \end{enumerate}
\end{lemma}
\begin{proof} 
(1) and (2) are consequences of Remark \ref{remark_4} and Corollary \ref{u_claim_15}, respectively. 
Then (3) follows from (2) and Lemma \ref{u_claim_17}.
\end{proof}
%%%%%%%%%%%%%%%%%%%%%%%%%%%%%%%%%%%%
%           Statement ↑
%%%%%%%%%%%%%%%%%%%%%%%%%%%%%%%%%%%%

Let $p$ and $q$ be perfect paths. 
As mentioned in Section \ref{preliminaries_3}, the $K$-linear map 
\begin{align}
    \theta_{p, q}^\prime: q \L \cap \L p  \rightarrow \Hom_{\L}(p\L, q\L)
\end{align}
given by $\theta_{p, q}^\prime(xp)(py)= xpy$ for $x, y \in \L$ is an isomorphism.
Recall that $\L = \bigoplus_{i \in \mathbb{Z}} \L_i$ satisfies that $\L_i = 0$ for $i < 0$.
For each $i \in \mathbb{Z}$, we define 
    \begin{align}
        (q\L \cap \L p)_{i} := q\L \cap \L_{i} p \quad \mbox{and} \quad (q\L p)_{i} := q \L_{i-l(q)} p. 
    \end{align}
From the definition, we have that $(q \L p)_{i} \subseteq (q\L \cap \L p)_{i}$ for any $i$;  $(q\L \cap \L p)_{i} = 0$ for $i <0$; and $(q \L p)_{i} = 0$ for $i < l(q)$.
We denote by $\P^{\mathbb{Z}}(X, Y)$ the subspace of $\HomZ_{\L}(X, Y)$ consisting of homomorphisms factoring through a graded projective $\L$-module.

%%%%%%%%%%%%%%%%%%%%%%%%%%%%%%%%%%%%
%           Statement ↓
%%%%%%%%%%%%%%%%%%%%%%%%%%%%%%%%%%%%
\begin{lemma} \label{u_claim_21}
Keep the above notation. We have the following statements.
\begin{enumerate}

    \item For any $i \in \mathbb{Z}$, the isomorphism $\theta_{p, q}^\prime: q\L \cap \L p \xrightarrow{\,\sim\,} \Hom_{\L}(p\L, q\L)$ induces the following two isomorphisms of $K$-vector spaces:
    \begin{align}
       (q\L \cap \L p)_{i} \xrightarrow{\,\sim\,} \HomZ_{\L}(p\L, q\L(i)) \quad \mbox{and} \quad \dfrac{(q\L \cap \L p)_{i}}{(q \L p)_{i}} \xrightarrow{\,\sim\,} \sHomZ_{\L}(p\L, q\L(i)).
    \end{align}

    \item Assume that $\HomZ_{\L}(p\L, q\L (i)) \not =0$ for an integer $i$. 
    Then $\sHomZ_{\L}(p\L, q\L (i)) \not =0$ if and only if $i < l(q)$.
    In this case, $\P^{\mathbb{Z}}(p\L, q\L (i)) = 0$, and the isomorphism $\theta_{p, q}^\prime: q\L \cap \L p \xrightarrow{\,\sim\,} \Hom_{\L}(p\L, q\L)$ induces an isomorphism
    \begin{align} \label{eq_8}
      \Theta^\prime_{p\L, q\L(i)}:  K p_{p\L, q\L(i)} = (q\L \cap \L p)_{i} \xrightarrow{\,\sim\,} \HomZ_{\L}(p\L, q\L (i)) = \sHomZ_{\L}(p\L, q\L (i)),
    \end{align} 
    where $p_{p\L, q\L(i)} = q^\prime p$ is the product of the left divisor $q^\prime$ of $q$ of length $i$ and the perfect path $p$.
    
\end{enumerate}
\end{lemma}
\begin{proof}
We first prove (1). 
A direct calculation yields the following commutative square of $K$-vector spaces:
\[\xymatrix{
  q\L \cap \L p \ar[r]_-{\sim}^-{\theta_{p, q}^\prime} &  \Hom_{\L}(p\L, q\L)\\
 (q\L \cap \L p)_{i} \ar[r] \ar@{^{(}->}[u] & \HomZ_{\L}(p\L, q\L (i)), \ar@{^{(}->}[u]
}\]  
where the vertical maps are the canonical inclusions.
It is clear that the lower horizontal $K$-linear map is injective.
For any $f \in \HomZ_\L(p\L, q\L(i))$, there exists a homogeneous element $x$ of $\L$ such that $ f(p) = xp$. Note that $xp$ belongs to $q\L(i)_{l(p)} = q\L_{i+l(p)}$.
Denote by $\deg a$ the degree of a homogeneous element $a$.
Then we have $\deg x + l(p) = \deg x +\deg p = \deg xp = \deg f(p) = i + l(p)$, so that $\deg x = i$  
and hence $xp \in (q\L \cap \L p)_{i}$.
Therefore, the lower horizontal map $(q\L \cap \L p)_{i} \to \HomZ_{\L}(p\L, q\L(i))$ is bijective.
On the other hand, let $\pi : \L \rightarrow q\L$ be the left multiplication by $q$.
Since $q$ is a homogeneous element of $\L$,  we  have that $\pi \in \HomZ_\L(\L, q\L(l(q)))$.
Then an argument as in the proof of \cite[Lemma 2.3]{Chen-Shen-Zhou_2018} shows that
\[
\P^{\mathbb{Z}}(p\L, q\L(i)) = \Im \HomZ_\L(\L(i-l(q)), \pi) = (q\L p)_i, 
\]
where the second equality can be obtained by using the isomorphism 
\[
(\L p)_{i-l(q)+l(p)} \xrightarrow{\sim} \HomZ_\L(p\L, \L(i-l(q)))
\]
induced by the isomorphism $\theta_{p}^\prime:\L p \xrightarrow{\sim} \Hom_{\L}(p\L, \L)$.

For (2), under the assumption that $0\not =\HomZ_{\L}(p\L, q\L (i)) \cong (q\L \cap \L p)_{i}$,  we first show that $\sHomZ_{\L}(p\L, q\L (i)) =0$ if and only if $i \geq  l(q)$.  
If $\sHomZ_{\L}(p\L, q\L (i)) =0$, then $q\L_{i-l(q)} p = (q\L p)_{i} =  (q\L \cap \L p)_{i}  \not = 0$, so that $i \geq  l(q)$. 
Conversely, suppose that $i \geq l(q)$. 
For a non-zero path $q^\prime p \in (q\L \cap \L p)_{i}$ with $q^\prime \in \B_{i}$,  the facts that $l(q^\prime) = i \geq l(q)$ and that $q^\prime p \in q\L$ imply that $q$ is a left divisor of $q^\prime$.
This shows that $(q\L \cap \L p)_{i} \subseteq (q\L p)_{i}$. 
Consequently, $\HomZ_{\L}(p\L, q\L (i)) = (q\L \cap \L p)_{i} = (q\L p)_{i} = \P^{\mathbb{Z}}(p\L, q\L (i))$.

For the last statement, assume that $\sHomZ_{\L}(p\L, q\L (i)) \not=0$. Then $\P^{\mathbb{Z}}(p\L, q\L (i))$ $= (q \L p)_{i} = q\L_{i-l(q)} p = 0$ since $i < l(q)$. 
On the other hand, $q$ contains one and only one path of length $i$ as a left divisor, say $q^\prime$.
In particular, $q^\prime p$ is the only non-zero path belonging to $(q\L \cap \L p)_{i}$.
Thus $(q\L \cap \L p)_{i}$ is spanned by $q^\prime p$.
\end{proof}
%%%%%%%%%%%%%%%%%%%%%%%%%%%%%%%%%%%%
%           Statement ↑
%%%%%%%%%%%%%%%%%%%%%%%%%%%%%%%%%%%%

%%%%%%%%%%%%%%%%%%%%%%%%%%%%%%%%%%%%
%           Subsection ↓
%%%%%%%%%%%%%%%%%%%%%%%%%%%%%%%%%%%%
\subsection{Main result} \label{section_3_2}
%%%%%%%%%%%%%%%%%%%%%%%%%%%%%%%%%%%%

In this subsection, we show the existence of a tilting object in  $\sGprojZ\L$. 
Further, computing the endomorphism algebra of the tilting object, we prove that $\sGprojZ\L$ is triangle equivalent to the bounded derived category of a hereditary algebra of Dynkin type $\A$.

Let $c = r_1\cdots r_n$ be an underlying cycle in $\L$ with each $r_i$ co-elementary.
For any integers  $i \leq j$, we set  
\begin{align}
    [i, j] =[i, j]_c :=r_i r_{i+1} \cdots  r_{j},
\end{align}
where the index of each co-elementary path is considered modulo $n$. 
For example, if $c= r_1r_2 r_3$, then $[-1, 5] = r_{-1} r_0 r_1 r_2 r_3 r_4 r_5  = r_2 r_3 r_1 r_2 r_3 r_1 r_2$ and $[7,7]=r_7 = r_1$.
It is understood that $l([i, j]) = 0$ when  $i > j$.
Also, we define a subset $\mathbb{P}_{\L}(c)$ of $\mathbb{P}_{\L} $ by
\begin{align} \label{eq_14}
    \mathbb{P}_{\L}(c) := \{ p \in \mathbb{P}_{\L} \mid r_1 \preceq p  \}
\end{align}
and set $m_c = m_c^{\L} := |\mathbb{P}_{\L}(c)|$, where $|S|$ denotes the cardinality of a set $S$.
In other words, the subset $\mathbb{P}_\L(c)$ is the vertex set of the connected component of the Hasse quiver $H(\mathbb{P}_\L, \preceq)$ that contains $r_1$.
We also note that $\mathbb{P}_\L(c)$ contains a unique maximal element with respect to $\preceq$, which is nothing but the elementary path containing $r_1$ as a left divisor, which can be written as $[1, m_c]$.   
It is straightforward that $\mathbb{P}_{\L}(c) =\{[1, i] \mid 1 \leq i \leq m_c\}$.
Recall that $\F$ denotes the set of minimal paths in the admissible ideal $I$ of $\L = KQ/I$ and that $X_c = \{ p \in \mathbb{E}_{\L} \mid c_p = c \mbox{ in } \C(\L) \}$.

%%%%%%%%%%%%%%%%%%%%%%%%%%%%%%%%%%%%
%           Statement ↓
%%%%%%%%%%%%%%%%%%%%%%%%%%%%%%%%%%%%
\begin{lemma} \label{u_claim_51}
Under the above situation, we have $X_c = \{ [i, i+m_c-1] \,|\, 1 \leq i \leq n \}$.
In particular, the path $[i, i+m_c]$ belongs to $\F$ for all $i \in \mathbb{Z}$.
\end{lemma} 
\begin{proof}
Let $p$ be an elementary path such that $c_p = c$ in  $\C(\L)$. 
By Theorem \ref{u_claim_25}, $p$ is of the form $[i, i+m-1]$, where $1 \leq i \leq n$ and $m>0$.
Then it is easy to see that the pair $([i, i+m-1], r_{i+m})$ is perfect, so that $[i, i+m]$ belongs to $\F$. 
Thus the non-zero right divisor $[i+1, i+m]$ is perfect by Theorem \ref{u_claim_25} again and must be elementary.
We then see by induction that $[i, i+m-1] \in \mathbb{E}_\L$ for all $1 \leq i \leq n$. 
But then, the fact that $[1, m_c]$ is elementary yields that $m=  m_c$.
\end{proof}
%%%%%%%%%%%%%%%%%%%%%%%%%%%%%%%%%%%%
%           Statement ↑
%%%%%%%%%%%%%%%%%%%%%%%%%%%%%%%%%%%%

Let $\T$ be a triangulated category with suspension functor $\Sigma$.
For a class $\mathcal{X}$ of objects of $\mathcal{T}$, we denote by $\thick_{\T} \mathcal{X}$ the smallest thick subcategory of $\T$ that contains $\mathcal{X}$. 
When $\mathcal{X}$ consists of a single object $X$, we write $\thick_{\T} X$ instead of $\thick_{\T}\{X\}$.
In what follows, we drop the index $\T$ in $\thick_{\T} \mathcal{X}$ when $\T$ is clear from the context. 
Recall that a object $T$ of $\T$ is called a {\it tilting object} if the following two conditions are satisfied:
\begin{enumerate}
    \item[(i)] $\Hom_{\T}(T, \Sigma^{i}T) = 0$ for all $i \not= 0$.
    \item[(ii)] $\thick T = \T$.
\end{enumerate}

A triangulated category is called {\it algebraic} if it is triangle equivalent to the stable category of a Frobenius category.
If $\T$ is an algebraic triangulated Krull-Schmidt category with a tilting object $T$, then $\T$ is triangle equivalent to the perfect derived category $\K^{\rm b}(\proj\End_{\T}T)$ of the endomorphism algebra $\End_{\T}T$ (\cite[4.3]{Keller_1994}).
See also \cite[Theorem 2.11]{Yamaura_2013}.

We define an object $T$ of $\sGprojZ\L$ as
\begin{align} \label{eq_2}
    T := \bigoplus_{c \in \C(\L)}\,\bigoplus_{0 \leq i < l(c)} T_c(i),
\end{align}
where  
\begin{align} \label{eq_15}
    T_{c} := \bigoplus_{p \in \mathbb{P}_{\L}(c)} p\L = \bigoplus_{1 \leq i \leq m_c} [1, i]\L.
\end{align}
It follows from the definition that $T$ depends on the choice of underlying cycles
and is basic in the sense that any two distinct indecomposable summands of $T$ are non-isomorphic. 
We will now show that $T$ is a tilting object of $\sGprojZ\L$. 
Recall that $\mathcal{E}_\L^\mathbb{Z} = \{ p\L (j) \mid p \in \mathbb{E}_\L, j \in \mathbb{Z} \}$ and that $Y_c = \{ r \in \mathbb{E}_{\L}^{\rm co} \mid c_r = c \mbox{ in } \C(\L) \}$ for any $c \in \C(\L)$.
We set 
\begin{align}
    (\mathcal{E}_\L^{\rm co})^\mathbb{Z} := \{ r \L (j) \mid r \in \mathbb{E}_\L^{\rm co}, j \in \mathbb{Z} \}.
\end{align}

%%%%%%%%%%%%%%%%%%%%%%%%%%%%%%%%%%%%
%           Statement ↓
%%%%%%%%%%%%%%%%%%%%%%%%%%%%%%%%%%%%
\begin{lemma} \label{u_claim_30}
We have $\thick T = \sGprojZ\L$, where $T$ is defined as in (\ref{eq_2}).
\end{lemma}
\begin{proof}
It follows from Proposition \ref{u_claim_16} (3) and Lemma \ref{u_claim_24} (1) that $\Sigma (\mathcal{E}_\L^{\rm co})^\mathbb{Z} = \mathcal{E}_\L^{\mathbb{Z}}$.
Therefore, in view of Lemma \ref{u_claim_24} (3), it suffices to show that $(\mathcal{E}_\L^{\rm co})^\mathbb{Z} \subseteq \thick T$. 
For this, we claim that $\{r \L (j) \,|\, r \in Y_c, j \in \mathbb{Z} \} \subseteq \thick \Tilde{T}_c$ for all $c \in \C(\L)$, where $\Tilde{T}_c:=\bigoplus_{0 \leq i < l(c)} T_c(i)$.
Then the formula $\mathbb{E}_{\L}^{\rm co} = \bigcup_{c \in \C(\L)} Y_c$ yields that 
\begin{align}
    (\mathcal{E}_\L^{\rm co})^\mathbb{Z} = \bigcup_{c\in \C(\L)} \{r \L (j) \,|\, r \in Y_c, j \in \mathbb{Z} \} \subseteq \thick T.
\end{align}

Let $c = r_1 \cdots r_n \in \C(\L)$ with $r_i \in \mathbb{E}_\L^{\rm co}$. 
Then $Y_c = \{r_i \,|\, 1 \leq i \leq n\}$.
Consider the following component of the Hasse quiver $H(\mathbb{P}_\L, \preceq)$:
\begin{align}
    [1,m_c] \rightarrow  \cdots \rightarrow [1,3] \rightarrow [1,2] \rightarrow r_1.
\end{align}
This gives rise to the following chain of the canonical inclusions in $\GprojZ\L$:
\begin{align}
    [1,m_c]\L \rightarrow  \cdots \rightarrow [1,3]\L \rightarrow [1,2]\L \rightarrow r_1\L.
\end{align}
For $1 \leq i \leq m_c$, we have the following short exact sequence in $\GprojZ$:
\begin{align} 
    0 
    \rightarrow  
    r_{i+1}\L(-l([1,i])) 
    \xrightarrow{\scriptscriptstyle \begin{bmatrix}-\iota_{i+1}\\\pi_{i+1}\\\end{bmatrix}} 
    e_{s(r_{i+1})}\L(-l([1,i])) \oplus [1, i+1]\L 
    \xrightarrow{\scriptscriptstyle \begin{bmatrix}\pi_i & \hspace{-1.5mm}\iota_i \end{bmatrix}}
   [1, i]\L
    \rightarrow
    0,
\end{align}
where $\iota_{i+1}$ and $\iota_i$ are the canonical inclusions, and $\pi_{i+1}$ and $\pi_i$ are the left multiplications by $[1, i]$. 
Then the short exact sequence induces the triangle
\begin{align} 
    r_{i+1}\L(-l([1,i])) \xrightarrow{\pi_{i+1}}   [1, i+1]\L \xrightarrow{\iota_{i}}  [1, i]\L \to \Sigma r_{i+1}\L(-l([1,i]))
\end{align}
in $\sGprojZ\L$.
On the other hand, Lemma \ref{u_claim_24} yields that $\Sigma^{-1}[1,m_c]\L = r_{m_c+1}\L(-l([1,m_c]))$.
Therefore, we have
\begin{align}
    \{r_{i}\L(-l([1,i-1]))\mid 1 \leq i \leq m_c+1\} \subseteq \thick T_c.
\end{align}
Since the automorphism $(1): \sGprojZ\L \to \sGprojZ\L$ is triangulated, we also obtain that
\begin{align}
    \{r_{i}\L(-l([1,i-1])+j)\mid 1\leq i \leq m_c+1\} \subseteq \thick T_c(j)
\end{align} 
for $0\leq j < l(c)$.
Define a set $\mathcal{M}_{i}^{c}$ by
\begin{align}
    \mathcal{M}_{i}^{c} := \begin{cases}
        \{r_{i}\L(-l([1,i-1])+j) \mid 0\leq j < l(c) \} & \mbox{if $i \geq 1$};\\        
        \{r_{i}\L(l([i, 0])+j) \mid 0\leq j < l(c) \} & \mbox{if $i<1$}.        
    \end{cases}
\end{align}

We now claim that $\bigcup_{i\in \mathbb{Z}}\mathcal{M}_{i}^{c} \subseteq \thick \Tilde{T}_c$.
For $1 \leq i \leq m_c+1$, we already know that $\mathcal{M}_{i}^{c} \subseteq \thick \Tilde{T}_c$.
Moreover, we see that
\begin{align}
  \thick \Tilde{T}_c \supseteq  \Sigma^{2}\mathcal{M}_{i}^{c} = \{ r_{i-m_c-1}\L(l([i-m_c-1, 0])+j) \mid 0\leq j < l(c)\} = \mathcal{M}_{i-m_c-1}^{c}.
\end{align} 
Applying $\Sigma^{2}$ repeatedly, we can deduce that $ \mathcal{M}_{i-k(m_c+1)}^{c} \subseteq \thick \Tilde{T}_c$ for all $k \geq 0$.
This implies that $ \mathcal{M}_{i}^{c} \subseteq \thick \Tilde{T}_c$ for all $i \leq 0$.
Similarly, one can use $\Sigma^{-2}$ to obtain that $ \mathcal{M}_{i}^{c} \subseteq \thick \Tilde{T}_c$ for all $i > m_c +1$.

Next, we claim that $\bigcup_{i\in \mathbb{Z}}\mathcal{M}_{i}^{c} =\{r_i\L(j)\mid 1\leq i \leq n, j \in \mathbb{Z}\}$, and this finishes the proof.
The inclusion $(\subseteq)$ is clear. 
So take an element $r_i\L(j)$ from the right hand side.
Then we have $j+l([1,i-1]) = j_1 l(c)+j_2$ for some $j_1 \in \mathbb{Z}$ and $0 \leq j_2< l(c)$.
Since
\begin{align}
    -l([1,i-1])+ j_1 l(c) = \begin{cases} 
        -l([1,i-j_1n-1]) & \mbox{if $j_1\leq 0$};\\
        l([i-(j_1-1)n, 0]) & \mbox{if $j_1>0$},
    \end{cases}
\end{align}
it follows that $r_i\L(j) = r_i\L(-l([1,i-1])+ j_1 l(c)+j_2)$ belongs to $\mathcal{M}_{i-j_1 n}^{c} \cup \mathcal{M}_{i-(j_1-1)n}^{c}$.
\end{proof}
%%%%%%%%%%%%%%%%%%%%%%%%%%%%%%%%%%%%
%           Statement ↑
%%%%%%%%%%%%%%%%%%%%%%%%%%%%%%%%%%%%

The following lemma will be used to show that $\sHomZ_\L(T, \Sigma^{\not=0} T) = 0$ and to determine the endomorphism algebra $\sEndZ_\L(T)$.

%%%%%%%%%%%%%%%%%%%%%%%%%%%%%%%%%%%%
%           Statement ↓
%%%%%%%%%%%%%%%%%%%%%%%%%%%%%%%%%%%%
\begin{lemma} \label{u_claim_57}
Let $k \in \Z$, $c=r_1 \cdots r_n \in \C(\L)$ with $r_i \in \mathbb{E}_\L^{\rm co}$ and $[i,j], [i^\prime, j^\prime] \in  \mathbb{P}_\L$. 
Then the following conditions are equivalent.
\begin{enumerate}
    \item $\sHomZ_\L([i,j]\L, [i^\prime, j^\prime] \L (k)) \not = 0$.
    \item There exists an $\alpha\in \mathbb{Z}$ such that $i^\prime \leq i +\alpha n \leq j^\prime \leq j+\alpha n  < i^\prime + m_c$ and $k = l([i^\prime, i+\alpha n -1])$.
\end{enumerate}
In this case, $\sHomZ_\L([i,j]\L, [i^\prime, j^\prime] \L (k)) \cong K$.
\end{lemma}
\begin{proof}
We first prove that (1) implies (2).
Under condition (1), Lemma \ref{u_claim_21} (2) implies that 
\begin{align}
    \sHomZ_\L([i,j]\L, [i^\prime, j^\prime] \L (k)) \cong ([i^\prime, j^\prime]\L \cap \L[i,j])_k = K q^\prime[i,j],
\end{align}
where $q^\prime$ is the proper left divisor of $[i^\prime, j^\prime]$ with $l(q^\prime) = k$.
If $k=0$, or equivalently, $q^\prime = e_{s(r_{i^\prime})}$, then we have $i^\prime \equiv i \ (\mod n)$, so that $i^\prime = i+ \alpha n$ for some $\alpha \in \Z$. 
Since $0\not= [i,j] \in [i^\prime, j^\prime]\L$, we can conclude that $i^\prime =i+\alpha n < \  j^\prime \leq j+\alpha n < i^\prime+m_c$. 
Also, one gets that $k = 0 =l([i^\prime, i^\prime-1]) = l([i^\prime, i+\alpha n -1])$.

Assume now that $q^\prime$ is non-trivial. 
Then $q^\prime$ is a perfect path, so that  $q^\prime = [i^\prime, s]$ for some $s \in \mathbb{Z}$ with $i^\prime \leq s \leq j^\prime$ and $s \equiv i-1 \ (\mod n)$.
Taking an $\alpha \in \mathbb{Z}$ such that $s = i+\alpha n -1$, we see that the non-zero path  $q^\prime[i,j] = [i^\prime, i+\alpha n -1][i+\alpha n, j+\alpha n]$ belongs to $[i^\prime, j^\prime]\L$.
Using the fact that $q^\prime$ is a proper left divisor of $[i^\prime, j^\prime]$, we deduce that $i^\prime < i +\alpha n \leq j^\prime \leq j+\alpha n  < i^\prime + m_c$, whereas $k =l(q^\prime) =l([i^\prime, i+\alpha n -1])$.

Conversely, suppose that (2) holds.
Then $[i^\prime, j+\alpha n]$ is a non-zero path in $\L$ that belongs to $[i^\prime, j^\prime]\L \cap \L[i+\alpha n, j+\alpha n]$, and we have that $[i^\prime, j+\alpha n] = [i^\prime, i+\alpha n -1][i+\alpha n, j+\alpha n]$.
Since $k = l([i^\prime, i+\alpha n -1]) < l([i^\prime, j^\prime])$,  Lemma \ref{u_claim_21} (2) concludes that $\sHomZ_\L([i,j]\L, [i^\prime, j^\prime] \L (k)) \not = 0$.
The last statement also follows from the same lemma.
\end{proof}
%%%%%%%%%%%%%%%%%%%%%%%%%%%%%%%%%%%%
%           Statement ↑
%%%%%%%%%%%%%%%%%%%%%%%%%%%%%%%%%%%%

%%%%%%%%%%%%%%%%%%%%%%%%%%%%%%%%%%%%
%           Statement ↓
%%%%%%%%%%%%%%%%%%%%%%%%%%%%%%%%%%%%
\begin{lemma} \label{u_claim_32} 
The object $T$ defined as in (\ref{eq_2}) is a tilting object of $\sGprojZ \L$.
\end{lemma}
\begin{proof}
Thanks to Lemma \ref{u_claim_30}, it is enough to prove that $\sHomZ_\L\left(T, \Sigma^{i} T\right) = 0$ for all $i\not =0$, where $T = \bigoplus_{c \in \C(\L)}\,\bigoplus_{0 \leq i < l(c)} T_c(i)$.
Since  $\sHomZ_\L\left(T_c, \Sigma^{\not=0}T_{c^\prime}(k)\right)=0$ for any $c \not = c^\prime \in \C(\L)$ and any $ k\in \mathbb{Z}$ by Corollary \ref{u_claim_19}, we have to show that $\sHomZ_\L\left(T_c, \Sigma^{\not=0}T_{c}(k)\right)=0$ for any $c \in \C(\L)$ and any $-l(c) < k < l(c)$.

Fix an underlying cycle $c = r_1 \cdots r_n$ with $r_i \in \mathbb{E}_\L^{\rm co}$.
A direct calculation shows that for any $m \in \mathbb{Z}$ and $1 \leq i^\prime \leq m_c$, we have
\begin{align} \label{eq_20}
    \Sigma^{2m+1}[1, i^\prime]\L = [i^\prime -(m+1)m_c -m, -m(m_c+1)]\L(d_{2m+1}),
\end{align}
where
\begin{align}
    d_{2m+1}=\begin{cases}
       l([i^\prime -(m+1)m_c -m, 0]) & \mbox{if $m \geq 0$};  \\
       -l([1,  i^\prime-(m+1)(m_c +1)]) & \mbox{if $m < 0$},  \\
    \end{cases}
\end{align}
and 
\begin{align} \label{eq_21}
    \Sigma^{2m}[1, i^\prime]\L =[-m (m_c +1) +1, i^\prime-m (m_c +1)]\L(d_{2m}),
\end{align}
where
\begin{align}
    d_{2m}=\begin{cases}
       l([-m (m_c +1) +1, 0]) & \mbox{if $m \geq 0$};  \\
       -l([1,  -m(m_c +1)]) & \mbox{if $m < 0$}.  \\
    \end{cases}
\end{align}

We claim that 
\begin{align}
   & \sHomZ_\L([1,i]\L, \Sigma^{2m+1}[1,i^\prime]\L(k))  \\
   =\ & \sHomZ_\L([1,i]\L,[i^\prime -(m+1)m_c -m, -m(m_c+1)]\L(d_{2m+1}+k)) \\
   =\ & 0
\end{align}
for any $1 \leq i, i^\prime \leq m_c$, $m \in \Z $ and $-l(c) <  k < l(c)$.
Assume to the contrary that the Hom space is non-zero for some $1 \leq i, i^\prime \leq m_c$, $m \in \Z $ and $-l(c) <  k < l(c)$.
By Lemma \ref{u_claim_57}, there exists an integer $\alpha$ such that
\begin{enumerate}
    \item[(i)] $i^\prime -(m+1)m_c -m \leq 1 +\alpha n \leq -m (m_c +1) \leq i+\alpha n  < i^\prime -m m_c -m$; 
    \item[(ii)] $d_{2m+1}+k = l([i^\prime -(m+1)m_c -m, \alpha n])$.
\end{enumerate}
If $m \geq 0$, then $\alpha < 0$ by condition (i).
Since 
\begin{align}
    d_{2m+1}-l(c) 
    & = l([i^\prime -(m+1)m_c -m, 0]) -l([1-n, 0]) \\
    & = l([i^\prime -(m+1)m_c -m, -n])  
\end{align}
and 
\begin{align}
    d_{2m+1}+l(c) 
    & = l([i^\prime -(m+1)m_c -m, 0]) +l([1, n]) \\
    & = l([i^\prime -(m+1)m_c -m, n]),
\end{align}
the inequalities $-l(c) <  k < l(c)$ can be deformed as follows.
\begin{align}
    l([i^\prime -(m+1)m_c -m, -n]) < d_{2m+1}+k < l([i^\prime -(m+1)m_c -m, n]).
\end{align}
Then we have $\alpha = 0$ by condition (ii), a contradiction.

Assume now that $m < 0$.
Then $\alpha >0$ by condition (i).
Moreover, the assumption that $\sHomZ_\L([1,i]\L, \Sigma^{2m+1}[1,i^\prime]\L(k)) \not =0$ implies that 
\begin{align}
    0 \leq d_{2m+1}+k <& -l([1,  i^\prime-(m+1)(m_c +1)]) + l(c)\\
    =& l([1, n])-l([1,  i^\prime-(m+1)m_c-m-1])\\
    =& l([i^\prime-(m+1)m_c -m, n]).
\end{align}
Then condition (ii) yields that $\alpha = 0$, a contradiction.

Similarly, one can show that $\sHomZ_\L\left([1,i]\L, \Sigma^{2m}[1,i^\prime]\L(k)\right)  =0$ for any $1 \leq i, i^\prime \leq m_c$, $m \in \Z\backslash\{0\}$ and $-l(c) <  k < l(c)$.
Therefore, we have $\sHomZ_\L\left(T, \Sigma^{i} T\right) = 0$ for all $i\not =0$.
\end{proof}
%%%%%%%%%%%%%%%%%%%%%%%%%%%%%%%%%%%%
%           Statement ↑
%%%%%%%%%%%%%%%%%%%%%%%%%%%%%%%%%%%%

For each  $c \in \C(\L)$, we denote by $\Ac$ the following linear quiver with $m_c$ vertices:
\begin{align}
   \Ac :  1 \to 2 \to \cdots \to m_c.
\end{align}
Further, we denote by $\T^{(k)}$ a direct product of $k$ copies of a triangulated category $\T$.  
The following theorem is the main result of this section.

%%%%%%%%%%%%%%%%%%%%%%%%%%%%%%%%%%%%
%           Statement ↓
%%%%%%%%%%%%%%%%%%%%%%%%%%%%%%%%%%%%
\begin{theorem} \label{u_claim_34}
Let $\L$ be a monomial algebra. 
Then there exists a triangle equivalence
\begin{align}
    \sGprojZ \L \ 
    \cong \prod_{c \in \C(\L)}  \D^{\rm b}\!\left(\mod K\Ac \right)^{(l(c))}.
\end{align}
\end{theorem}
\begin{proof}
By Lemma \ref{u_claim_32}, $T = \bigoplus_{c \in \C(\L)}\,\bigoplus_{0 \leq i < l(c)} T_c(i)$ is a tilting object of $\sGprojZ \L$. 
We now determine its endomorphism algebra $\End_{\sGprojZ \L}(T)$.

Let $c=r_1 \cdots r_n$ be an underlying cycle with $r_i \in \mathbb{E}_\L^{\rm co}$, and let $1 \leq i, i^\prime \leq m_c$ and $k \in \Z$.
It follows from Lemma \ref{u_claim_57} that $\HomZ_\L([1,i]\L, [1,j]\L(k)) \not=0$ if and only if there exists an $\alpha\in \mathbb{Z}$ such that $1 \leq 1 +\alpha n \leq i^\prime \leq i+\alpha n  < 1 + m_c$ and $k = l([1, \alpha n])$.
Thus one obtains that 
    \begin{align} \label{eq_17}
        \sHomZ_\L\left([1, i]\L, [1, i^\prime]\L(k) \right) \cong
        \begin{cases}
            K & \mbox{ if $i^\prime \leq i$ and $k=0$;}\\
            0 & \mbox{ otherwise.}
        \end{cases}
    \end{align}

Now, it follows from (\ref{eq_17}) that the endomorphism algebra $\sEndZ_{\L}\left( T_c \right)$ of the object $T_c = \bigoplus_{1 \leq i \leq m_c} [1,i]\L$ is isomorphic to the algebra of $\left( m_c \times m_c \right)$-upper triangular matrices over $K$
\begin{align}
\begin{pmatrix}
    K & K & \cdots & K \\
     & K & \cdots & K \\
      &   & \ddots &  \vdots \\
     &  &  & K \\
\end{pmatrix},
\end{align}
which is isomorphic to $K \Ac$.
We denote by $\Gamma^{(k)}$ a direct product of $k$ copies of an algebra $\Gamma$.
Then we have the following isomorphisms of algebras:  
\begin{align}
    \End_{\sGprojZ \L}(T) 
    &= \sEndZ_{\L}\left( \bigoplus_{c \in \C(\L)}\,\bigoplus_{0\leq i <l(c)} T_c(i) \right) & & \\[2mm]
    &\cong \prod_{c \in \C(\L)} \sEndZ_{\L}\left(\bigoplus_{0\leq i <l(c)} T_c(i)\right) & \mbox{(by Corollary \ref{u_claim_19})} & \\[3mm]
    &\cong \prod_{c \in \C(\L)}  \sEndZ_{\L}\left(T_c\right)^{(l(c))}  & \mbox{(by (\ref{eq_17}))} & \\[3mm]
    &\cong \prod_{c \in \C(\L)}  \left(K \Ac\right)^{(l(c))}  & & 
\end{align} 
Now, we obtain the following triangle equivalences:
\begin{align}
    \sGprojZ \L 
    \cong \K^{\rm b}\!\left(\proj \End_{\sGprojZ \L}(T) \right)
    \cong \D^{\rm b}\!\left(\mod \End_{\sGprojZ \L}(T) \right)
    \cong \prod_{c \in \C(\L)} \!\D^{\rm b}\!\left(\mod K\Ac \right)^{(l(c))},
\end{align}
where the first equivalence follows from a result of Keller \cite[4.3]{Keller_1994} (cf.~\cite[Theorem 2.11]{Yamaura_2013}), the second equivalence from the fact that the endomorphism algebra $\End_{\sGprojZ \L}(T)$ is hereditary and thus has finite global dimension, and the third equivalence from the fact that $\End_{\sGprojZ \L}(T) \cong \prod_{c \in \C(\L)}  \left(K \Ac\right)^{(l(c))}$ is a direct product of algebras.
This completes the proof.
\end{proof}
%%%%%%%%%%%%%%%%%%%%%%%%%%%%%%%%%%%%
%           Statement ↑
%%%%%%%%%%%%%%%%%%%%%%%%%%%%%%%%%%%%

%%%%%%%%%%%%%%%%%%%%%%%%%%%%%%%%%%%%
%           Statement ↓
%%%%%%%%%%%%%%%%%%%%%%%%%%%%%%%%%%%%
\begin{remark}{\rm \label{remark_1}
When $\L$ is Iwanaga-Gorenstein, the theorem describes the graded singularity category $\D_{\rm sg}(\modZ\L)$ of $\L$.
In particular, the theorem improves and extends \cite[Theorem 5.2.2]{Lu-Zhu_2021}. 
}\end{remark}
%%%%%%%%%%%%%%%%%%%%%%%%%%%%%%%%%%%%
%           Statement ↑
%%%%%%%%%%%%%%%%%%%%%%%%%%%%%%%%%%%%

%%%%%%%%%%%%%%%%%%%%%%%%%%%%%%%%%%%%
%           Statement ↓
%%%%%%%%%%%%%%%%%%%%%%%%%%%%%%%%%%%%
\begin{example}{\rm \label{ex_4}
Let $\L$ be as in Example \ref{ex_3}.
We have seen there that $\C(\L)=\{c=a_{123}, c^\prime=a_{45}\}$.  
Note that $m_{c}=|\mathbb{P}_{\L}(c)| = 4$ and $m_{c^\prime}=|\mathbb{P}_{\L}(c^\prime)| = 3$.
Then Theorem \ref{u_claim_34} yields the following triangle equivalence:
\begin{align}
    \sGprojZ \L\, \cong\,  \D^{\rm b}\!\left(\mod K\Ac \right)^{(3)} \times \D^{\rm b}\!\left(\mod K\A_{c^\prime} \right)^{(2)},
\end{align}
where $\Ac:1 \to 2 \to 3 \to 4$ and $\A_{c^\prime}:1 \to 2 \to 3$.
}\end{example}
%%%%%%%%%%%%%%%%%%%%%%%%%%%%%%%%%%%%
%           Statement ↑
%%%%%%%%%%%%%%%%%%%%%%%%%%%%%%%%%%%%

So far, we have seen that
\begin{align} \label{eq_6}
    \sGprojZ \L = \prod_{c \in \C(\L)}\,\prod_{0 \leq i <l(c)} \thick T_c(i)
\end{align}
and that $\thick T_c(i) \cong \D^{\rm b}(\mod K\Ac)$ for every $c \in \C(\L)$ and $0 \leq i < l(c)$.
Recall that $\thick T_c(i)$ is the smallest thick subcategory of $\sGprojZ \L$ that contains $T_c(i)$.

For any $i \in \mathbb{Z}$, the automorphism $(i):\sGprojZ \L\to \sGprojZ \L$ restricts to an isomorphism  $(i):\thick T_c \to (\thick T_c)(i)$, where $(\thick T_c)(i)$ is the full subcategory of $\sGprojZ \L$ defined by
\begin{align}
    (\thick T_c)(i) =\{ X (i) \mid X \in \thick T_c\}.
\end{align}
It is easy to see that $(\thick T_c)(i) = \thick T_c(i)$.
We end this subsection with the following lemma.

%%%%%%%%%%%%%%%%%%%%%%%%%%%%%%%%%%%%
%           Statement ↓
%%%%%%%%%%%%%%%%%%%%%%%%%%%%%%%%%%%%
\begin{lemma} \label{u_claim_39}
Let $c\in\C(\L)$ and $i, j \in \mathbb{Z}$.
Then $\thick T_c (i) = \thick T_c(j)$ if and only if $i \equiv j \pmod{l(c)}$.
In particular, the automorphism  $(l(c)):\sGprojZ \L\to \sGprojZ \L$ restricts to an automorphism  $(l(c)):\thick T_c \to \thick T_c$.
\end{lemma} 
\begin{proof}
For the equivalence of the conditions, it is enough to show that  $\thick T_c (i) = \thick T_c$ if and only if $i \equiv 0 \pmod{l(c)}$. 
Take co-elementary paths $r_1, \ldots, r_n$ such that $c= r_1 \cdots r_n$.
For the\,\lq\lq if\,\rq\rq\,part, take an integer $i$ such that $i \equiv 0 \pmod{l(c)}$.
A direct computation shows that
\begin{align}
    \mathcal{M}_{kn+1}^{c} =  \{r_{1}\L(-k l(c)+j) \mid 0\leq j < l(c) \} \  \mbox{for any $k \in \mathbb{Z}$},
\end{align}
where $\mathcal{M}_{kn+1}^{c}$ is defined as in the proof of Lemma \ref{u_claim_30}.
Recall that $r_1\L(-k  l(c)+j) \in \thick T_c(j)$.
Then it follows that $[1,s]\L(-k l(c)) \in \thick T_c$ for all $1 \leq s \leq m_c$ and $k \in \mathbb{Z}$, which implies that both $\thick T_c(i)$ and  $\thick T_c(-i)$ are contained in $\thick T_c$.
Thus we have $\thick T_c(i) = \thick T_c$, because 
\begin{align}
    \thick T_c = (\thick T_c)(-i)(i) =(\thick T_c(-i))(i) \subseteq (\thick T_c)(i) = \thick T_c(i).
\end{align}
Conversely, assume that there exists an integer $i$ such that $\thick T_c(i) = \thick T_c$ and such that $i \not\equiv 0 \ (\mod l(c))$. 
Choose integers $k$ and $k^\prime$ with $0 < k^\prime < l(c)$ such that $i=k l(c) + k^\prime$.
Then we have  
\begin{align}
    \thick T_c =(\thick T_c)(k l(c) + k^\prime) = (\thick T_c(k l(c)))(k^\prime) = \thick T_c(k^\prime),
\end{align}
where we use the\,\lq\lq if\,\rq\rq\,part to obtain the third equality.
But then, since 
\begin{align}
    \mbox{$\HomZ_\L(\thick T_c, \thick T_c(j)) = 0$ for any $0< j <l(c)$,}
\end{align}
we have  $\thick T_c =0$, a contradiction. 
\end{proof}
%%%%%%%%%%%%%%%%%%%%%%%%%%%%%%%%%%%%
%           Statement ↑
%%%%%%%%%%%%%%%%%%%%%%%%%%%%%%%%%%%%

%%%%%%%%%%%%%%%%%%%%%%%%%%%%%%%%%%%%
%           Subsection ↓
%%%%%%%%%%%%%%%%%%%%%%%%%%%%%%%%%%%%
\subsection{Auslander-Reiten triangles} \label{section_3_3}
%%%%%%%%%%%%%%%%%%%%%%%%%%%%%%%%%%%%

We will recall several facts on Auslander-Reiten triangles from \cite{Happel_Book,Reiten-VandenBergh_2002}. 
Let $\mathcal{T}$ be a Hom-finite $K$-linear Krull-Schmidt triangulated category with suspension functor $\Sigma$.
A triangle $A \xrightarrow{f} B \xrightarrow{g} C \xrightarrow{h} \Sigma A$ in $\T$ is called {\it Auslander-Reiten} if the following three conditions are satisfied:
\begin{enumerate}[label=(\roman*)]
    \item[(AR1)] $A$ and $C$ are indecomposable.
    \item[(AR2)] $h \not= 0$.
    \item[(AR3)] If $\alpha: W \rightarrow C$ is not a retraction, then $\alpha$ factors through $g$.
\end{enumerate}
We say that {\it $\T$ has right ({\rm resp}.~left) Auslander-Reiten triangles}  if for any indecomposable object $C$ (resp.~$A$), there is a triangle satisfying the above conditions. 
We say that {\it $\T$ has Auslander-Reiten triangles} if $\T$ has both left and right Auslander-Reiten triangles.
A {\it right Serre functor} is an additive endofunctor $F: \T \rightarrow \T$ together with bi-functorial isomorphisms
\begin{align}
    \eta_{A, B} : \Hom_{\T}(A, B) \rightarrow D\Hom_{\T}(B, FA)
\end{align}
for any $A, B  \in \T$, where $D=\Hom_{K}(-, K)$ is the $K$-duality.
The right Serre functor is unique up to natural isomorphism (if it exists).
{\it Left Serre functors} are defined dually.
A {\it Serre functor} is a functor that is both left and right Serre, or equivalently, a right Serre functor that is an equivalence.
Serre functors are known to be triangulated. 
It is also known that $\T$ has Auslander-Reiten triangles if and only if there exists a Serre functor $F$ for $\T$.
In this case, the composition $\tau = \Sigma^{-1}F$ satisfies that $\tau C = A$ in any Auslander-Reiten triangle $A \to B \to C \to \Sigma A$. 
The autoequivalence $\tau$ is called the {\it Auslander-Reiten translation} for $\T$.

The bounded derived category $\D^{\rm b}(\mod \Gamma )$ of an algebra $\Gamma$ of finite global dimension has Auslander-Reiten triangles (\cite[Examples 3.2]{Bondal-Kapranov_1990}), so that $\sGprojZ\L$ has Auslander-Reiten triangles by Theorem \ref{u_claim_34}.
We will now construct an Auslander-Reiten triangle $A \to B \to C \to \Sigma A$ for any indecomposable object $C \in \sGprojZ\L$.  
To do this, we can focus on indecomposable objects of the form $p\L$ with $p \in \mathbb{P}_\L$. 
The reason for this is that the degree-shifted triangle $A(i) \xrightarrow{f(i)} B(i) \xrightarrow{g(i)} C(i) \xrightarrow{h(i)} \Sigma A(i)$ of an Auslander-Reiten triangle $A \xrightarrow{f} B \xrightarrow{g} C \xrightarrow{h} \Sigma A$  is also Auslander-Reiten for any $i \in \mathbb{Z}$.

Given a perfect path $p$, there exists an underlying cycle $c=r_1 \cdots r_n$ with $r_i \in \mathbb{E}_\L^{\rm co}$ and integers $1 \leq i \leq n$ and $m>0$ such that $p = r_i r_{i+1} \cdots r_{i+m-1} = [i, i+m-1]$.
Letting 
\begin{align} \label{eq_19}
    A_p := [i+1, i+m]\L (-l(r_i)) 
\end{align}
and $    B_p := [i+1, i+m-1]\L(-l(r_i)) \oplus [i,i+m]\L$,
we have the following short exact sequence in $\GprojZ\L$:
\begin{align} 
    0 
    \rightarrow  
    A_p
    \xrightarrow{\scriptscriptstyle \begin{bmatrix}\iota^\prime\\-\pi^\prime\\\end{bmatrix}} 
    B_p
    \xrightarrow{\scriptscriptstyle \begin{bmatrix}\pi & \hspace{-1.5mm}\iota\end{bmatrix}}
    p\L 
    \rightarrow
    0,
\end{align}
where $\iota^\prime = \iota^\prime_p$ and  $\iota = \iota_p$ are the canonical inclusions, and $\pi^\prime = \pi^\prime_p$ and  $\pi = \pi_p$ are the left multiplications by $r_i$.
We now show that the following induced triangle of $\sGprojZ\L$ is Auslander-Reiten:
\begin{align} \label{eq_5}
    A_p
    \xrightarrow{\scriptscriptstyle \begin{bmatrix}\iota^\prime\\-\pi^\prime\\\end{bmatrix}} 
    B_p
    \xrightarrow{\scriptscriptstyle \begin{bmatrix}\pi & \hspace{-1.5mm}\iota\end{bmatrix}}
    p\L 
    \xrightarrow{\scriptscriptstyle \xi}
    \Sigma A_p,
\end{align}
where $\xi = \xi_p$ is the image of the non-zero path $p_{p\L, \Sigma A_p}$ under the $K$-linear isomorphism 
    \begin{align} 
     \Theta^\prime_{p\L, \Sigma A_p}: K p_{p\L, \Sigma X_p} \xrightarrow{\,\sim\,} \sHomZ_\L(p\L, \ \Sigma A_p).
    \end{align}
Note that $\Sigma A_p = [i+m-m_c, i]\L(l([i+m-m_c, i-1]))$ and 
\begin{align}
    p_{p\L, \Sigma A_p}= [i+m-m_c, i-1]p =[i+m-m_c, i+m-1].
\end{align}

%%%%%%%%%%%%%%%%%%%%%%%%%%%%%%%%%%%%
%           Statement ↓
%%%%%%%%%%%%%%%%%%%%%%%%%%%%%%%%%%%%
\begin{proposition} \label{u_claim_36}
The triangle (\ref{eq_5}) is Auslander-Reiten.
\end{proposition} 
\begin{proof} 
We only have to show that (AR3) holds for the triangle  (\ref{eq_5}).
Let $\alpha : W \to p\L$ be a non-zero morphism in $\sGprojZ\L$ that is not a retraction.  
Without loss of generality, we may assume that $W = q\L(-t)$ for some $q \in \mathbb{P}_\L$ and $t \geq 0$.
If $t=0$, then $\alpha$ is (represented by) the canonical inclusion $q\L \hookrightarrow p\L$. 
Thus $p$ is a left divisor of $q$.
But then $p$ must be proper since $\alpha$ is not a retraction.
This implies that the perfect path $q$ contains $p r_{i+m}  = [i, i+m-1] r_{i+m} =[i, i+m]$ as a left divisor. 
Namely, $\Im \alpha = q\L \subseteq [i, i+m]\L$.
In other words, we have the following commutative diagram in $\sGprojZ\L$:
\begin{align}
    \xymatrix@C=15mm@R=7mm@M=2mm{
     & W \ar[d]^-{\alpha} \ar[ld]_-{\scriptscriptstyle \begin{bmatrix}0\\ \alpha \end{bmatrix}} \\
    B_p \ar[r]_-{\scriptscriptstyle \begin{bmatrix}\pi & \hspace{-1.5mm}\iota \end{bmatrix}} & p\L \\
    }
\end{align}

Now, suppose that $t>0$. 
Since $\sHomZ_\L(q\L, p\L(t)) = \sHomZ_\L(q\L(-t), p\L)\not = 0$ by assumption, it follows that
\begin{align}
    p = [i, i+t^\prime-1] [i+t^\prime, i+m-1] \quad \mbox{and} \quad  q= [i+t^\prime, i+m-1] q^\prime
\end{align}
for some $[i, i+t^\prime-1] \in \B_{t}, [i+t^\prime, i+m-1] \in \B_{>0}$ and $q^\prime \in \B$ with $p q^\prime = [i, i+t^\prime-1]q \not= 0$ in $\L$, where $t^\prime >0$.
By Lemma \ref{u_claim_21}, there exists a $\lambda \in K$ such that $\alpha(q) = \lambda [i, i+t^\prime-1]q$.
If $q^\prime \in\B_{>0}$, then $q^\prime$ becomes a perfect path and hence contains the co-elementary path $r_{i+m}$ as a left divisor, so that  $[i,i+t^\prime-1]q$ contains $[i,i+m]$ as a left divisor. This means that $\Im \alpha = [i,i+t^\prime-1]q\L \subseteq [i,i+m]\L$.
Therefore, we obtain the desired commutative diagram as in the above.

Let us consider the remaining case where $q^\prime \in \B_{0}$, that is, $p = [i, i+t^\prime-1]q$.
Then we have $[i+1, i+t^\prime-1]q = [i+1, i+m-1]$ and $l([i+1, i+t^\prime-1]) = t-l(r_i)$.
Define $\alpha^\prime \in \sHomZ_\L(W, [i+1, i+m-1]\L(-l(r_{i})))$ by 
$\alpha^\prime(q) = \lambda [i+1,i+t-1] q$.
Then we get the following commutative diagram in $\sGprojZ\L$:
\begin{align}
    \xymatrix@C=15mm@R=7mm@M=2mm{
     & W \ar[d]^-{\alpha} \ar[ld]_-{\scriptscriptstyle \begin{bmatrix}\alpha^\prime\\ 0 \end{bmatrix}} \\
    B_p \ar[r]_-{\scriptscriptstyle \begin{bmatrix}\pi & \hspace{-1.5mm}\iota \end{bmatrix}} & p\L \\
    }
\end{align}
This completes the proof.
\end{proof}
%%%%%%%%%%%%%%%%%%%%%%%%%%%%%%%%%%%%
%           Statement ↑
%%%%%%%%%%%%%%%%%%%%%%%%%%%%%%%%%%%%

We will now construct the Auslander-Reiten translation $\tau : \sGprojZ\L \to \sGprojZ\L$ for  $\sGprojZ\L$.
It is worth noting that if the Auslander-Reiten translation exists, it is uniquely determined up to functorial isomorphism.
Therefore, the Auslander-Reiten translations on the equivalent triangulated categories $\sGprojZ\L$ and $\prod_{c \in \C(\L)}  \D^{\rm b}\!\left(\mod K\Ac \right)^{(l(c))}$ agree up to functorial isomorphism.
For any $c \in \C(\L)$ and $0\leq i < l(c)$, the fact that $\thick T_c(i) \cong \D^{\rm b}(\mod K\Ac)$ implies that $\thick T_c(i)$ admits an Auslander-Reiten translation $\tau_{c, i}: \thick T_c(i) \to \thick T_c(i)$.
We put $\tau_c := \tau_{c, 0}$. 
Then it is easy to verify the existence of the following commutative square of triangulated categories:
\begin{align}
    \xymatrix@C=10mm@R=10mm@M=2mm{
    \thick T_c(i) \ar[r]^-{\tau_{c, i}} \ar[d]_-{(-i)} & \thick T_c(i)  \\
    \thick T_c \ar[r]^-{\tau_{c}} & \thick T_c \ar[u]_-{(i)}.
    }
\end{align}
It also follows from (\ref{eq_6}) that $\tau: \sGprojZ\L \to \sGprojZ\L$ can be decomposed as follows:
\begin{align}
    \tau = \prod_{c \in \C(\L)}\,\prod_{0\leq i <l(c)}\tau_{c, i}.
\end{align}
Thus, in order to describe $\tau$, it is enough to determine $\tau_{c}:\thick T_c \to \thick T_c$ for all $c \in \C(\L)$. 
To do this, it suffices to specify the actions of $\tau_{c}$ on indecomposable objects and the morphisms between them (since $\thick T_c$ is Krull-Schmidt). 
This can be accomplished by using Lemma \ref{u_claim_21}, Proposition \ref{u_claim_36}, and \cite[Section I.2]{Reiten-VandenBergh_2002}. 
We leave the details to the reader.
It can be seen that the autoequivalence $\tau_c: \thick T_c \to \thick T_c $ is given as follows. 
For any indecomposable objects $p\L(j), q\L(k) \in \thick T_c$,\vspace{1mm}
\begin{enumerate}
  \setlength{\itemsep}{1mm} 

    \item[$\bullet$] $\tau_c(p\L(j)) = A_p (j)$, where $A_p$ is defined as in (\ref{eq_19}). 
\item[$\bullet$] The action $\tau_c : \sHomZ_\L(p\L(j), q\L(k)) \to \sHomZ_\L(\tau_c(p\L(j)), \tau_c(q\L(k)))$ on morphisms is defined to be the map making the following square of $K$-vector spaces commute:
\begin{align}
    \xymatrix{
    Kp_{p\L(j), q\L(k)}\ar[r] \ar[d]_-{\Theta^\prime_{p\L(j), q\L(k)}} & Kp_{\tau_c (p\L(j)), \tau_c (q\L(k))} \ar[d]^-{\Theta^\prime_{\tau_c (p\L(j)), \tau_c (q\L(k))}}\\
    \sHomZ_\L(p\L(j), q\L(k)) \ar[r]^-{\tau_c} & \sHomZ_\L(\tau_c (p\L(j)), \tau_c (q\L(k))),
    }
\end{align} 
where the upper horizontal map sends $p_{p\L(j), q\L(k)}$ to $p_{\tau_c (p\L(j)), \tau_c (q\L(k))}$.
\end{enumerate}

The following proposition is a consequence of the definition of $\tau_c$.
Recall that $|c|$ denotes the number of co-elementary subpaths of $c$.

%%%%%%%%%%%%%%%%%%%%%%%%%%%%%%%%%%%%
%           Statement ↓
%%%%%%%%%%%%%%%%%%%%%%%%%%%%%%%%%%%%
\begin{proposition} \label{u_claim_40}
Let $c \in \C(\L)$. 
Then the Auslander-Reiten translation $\tau_c: \thick T_c \to \thick T_c$ satisfies that $\tau_c ^{|c|} \cong (-l(c))$ as autoequivalences.
\end{proposition}
\begin{proof}
For any $A \in \ind \thick T_c$, there exists a perfect path $p$ and an integer $j$ such that $A \cong p\L(j)$.
We see from Corollary \ref{u_claim_19} that 
\begin{align}
    \sHomZ_\L(\thick T_c, \thick T_{c^\prime}(i)) = 0 \ \mbox{ for any $i \in \mathbb{Z}$}
\end{align}
whenever $c \not= c^\prime$ in $\C(\L)$.
Thus the fact that $0\not=A \in \thick T_c \cap \thick T_{c_p}(j)$ implies that $c_p = c$ in $\C(\L)$.
Let $c= r_1 \cdots r_n$ with $r_i \in \mathbb{E}_\L^{\rm co}$.
Then $p = r_i  \cdots r_{i+m-1} = [i, i+m-1]$ for some $1 \leq i \leq n$ and $m >0$.
Consequently, we have 
\begin{align}
    \tau_c^{|c|}A 
    &\cong \tau_c^n [i, i+m-1]\L(j)\\ 
    &= \tau_c^{n-1} [i+1, i+m]\L(-l(r_i)+j) \\
    &= \tau_c^{n-2} [i+2, i+m+1]\L(-l([i, i+1])+j) \\
    &\ \ \vdots\\
    &= [i+n, i+m+n-1]\L(-l([i, i+n-1])+j) \\
    &= [i, i+m-1]\L(j)(-l(c)) \\
    &\cong A(-l(c)),
\end{align}
which is functorial in $A$. 
This completes the proof. 
\end{proof}
%%%%%%%%%%%%%%%%%%%%%%%%%%%%%%%%%%%%
%           Statement ↑
%%%%%%%%%%%%%%%%%%%%%%%%%%%%%%%%%%%%

For a Hom-finite $K$-linear Krull-Schmidt category $\C$, we denote by $\G(\C)$ the Auslander-Reiten quiver of $\C$.
We refer to \cite{ShipingLiu_2010} for its definition. 
The following presents an example of the Auslander-Reiten quiver of $\sGprojZ\L$.
It follows from Theorem \ref{u_claim_34} that the Auslander-Reiten quiver of $\sGprojZ\L$ coincides with that of $\prod_{c \in \C(\L)}  \D^{\rm b}\!\left(\mod K\Ac \right)^{(l(c))}$.

%%%%%%%%%%%%%%%%%%%%%%%%%%%%%%%%%%%%
%           Statement ↓
%%%%%%%%%%%%%%%%%%%%%%%%%%%%%%%%%%%%
\begin{example}{\rm \label{ex_6}
Let $\L$ be as in Example \ref{ex_3}.
We know that $\C(\L)=\{c=a_{123}, c^\prime=a_{45}\}$ with $m_c= 4$ and $m_{c^\prime} = 3$ and $\mathbb{E}_{\L}^{\rm co}=\{r_1 = a_{12}, r_2 =a_{3}, r^\prime_1 =a_{45}\}$.
Then $T_c = \bigoplus_{1 \leq i \leq 4} [1, i]_c \L$ and $T_{c^\prime} = \bigoplus_{1 \leq i \leq 3} [1, i]_{c^\prime} \L$, where $c=r_1r_2$ and $c^\prime= r^\prime_1$. 
For $0\leq i < l(c)= 3$, the Auslander-Reiten quiver of $\thick T_c(i)$  is depicted in Figure \ref{fig_1}.
Further, for $0\leq i < l(c^\prime)= 2$, the Auslander-Reiten quiver of $\thick T_{c^\prime}(i)$ is depicted in Figure \ref{fig_2}.
Finally, the Auslander-Reiten quiver of $\sGprojZ\L$ is given by a disjoint union of the $\G(\thick T_c(i))$ and the $\G(\thick T_{c^\prime}(i))$.
}\end{example}
%%%
\begin{figure}[htbp] 
%%%%%%%%%%%%%%%%%%%%%%%%%%%%%
%\small
%%%%%%%%%%%%%%%%%%%%%%%%%%%%%
% diagram ↓
%%%%%%%%%%%%%%%%%%%%%%%%%%%%%
\begin{tikzpicture}%[scale=1]
%%%
\node[outer sep=2pt] (-3) at (-2.6,-2.7) {$\cdots$ \ \ };
\node[outer sep=2pt] (-1) at (-2.6,-0.0) {$\cdots$  \ \ };
\node[outer sep=2pt] (0) at (-1.7,1.35) {$r_4\Lambda(i-5)$};
\node[outer sep=2pt] (2) at (-1.2,-1.35) {$[3,5]_c\Lambda(i-3)$};
\node[outer sep=2pt] (3) at (0.1,0) {$[3,4]_c\Lambda(i-3)$};
\node[outer sep=2pt] (4) at (1.4,1.35) {$r_3\Lambda(i-3)$};
\node[outer sep=2pt] (5) at (0.5,-2.7) {$[2,5]_c\Lambda(i-2)$};
\node[outer sep=2pt] (6) at (1.8,-1.35) {$[2,4]_c\Lambda(i-2)$};
\node[outer sep=2pt] (7) at (3.1,0) {$[2,3]_c\Lambda(i-2)$};
\node[outer sep=2pt] (8) at (4.4,1.35) {$r_2\Lambda(i-2)$};
\node[outer sep=2pt] (9) at (3.4,-2.7) {$[1,4]_c\Lambda(i)$};
\node[outer sep=2pt] (10) at (4.7,-1.35) {$[1,3]_c\Lambda(i)$};
\node[outer sep=2pt] (11) at (6.0,0) {$[1,2]_c\Lambda(i)$};
\node[outer sep=2pt] (12) at (7.3,1.35) {$r_1\Lambda(i)$};
\node[outer sep=2pt] (13) at (6.3,-2.7) {$[0,3]_c\Lambda(i+1)$};
\node[outer sep=2pt] (14) at (7.6,-1.35) {$[0,2]_c\Lambda(i+1)$};
\node[outer sep=2pt] (15) at (8.9,0) {$[0,1]_c\Lambda(i+1)$};
\node[outer sep=2pt] (21) at (9.2,-2.7) {$[-1,2]_c\Lambda(i+3)$};
\node[outer sep=2pt] (22) at (10.1,-1.35) { \ \ $\cdots$};
\node[outer sep=2pt] (24) at (10.1,1.35) { \ \ $\cdots$};
%%%
\draw[->] (0) -- (3);
\draw[->] (2) -- (3);
\draw[->] (2) -- (5);
\draw[->] (3) -- (4);
\draw[->] (3) -- (6);
\draw[->] (4) -- (7);
\draw[->] (5) -- (6);
\draw[->] (6) -- (7);
\draw[->] (6) -- (9);
\draw[->] (7) -- (8);
\draw[->] (7) -- (10);
\draw[->] (8) -- (11);
\draw[->] (9) -- (10);
\draw[->] (10) -- (11);
\draw[->] (10) -- (13);
\draw[->] (11) -- (12);
\draw[->] (11) -- (14);
\draw[->] (12) -- (15);
\draw[->] (13) -- (14);
\draw[->] (14) -- (21);
\draw[->] (14) -- (15);
% $\tau$-orbit
\draw[dashed] (-3) -- (5) -- (9) -- (13) --  (21);
\draw[dashed]  (2) -- (6) -- (10) -- (14) --  (22);
\draw[dashed] (-1) -- (3) -- (7) -- (11) -- (15);
\draw[dashed] (0) -- (4) -- (8) -- (12) --  (24);
\end{tikzpicture}
%%%%%%%%%%%%%%%%%%%%%%%%%%%%%
% diagram ↑
%%%%%%%%%%%%%%%%%%%%%%%%%%%%%
%\normalsize
%%%%%%%%%%%%%%%%%%%%%%%%%%%%%
    \caption{The Auslander-Reiten quiver of $\thick T_c(i)$, where the valuation of each arrow is $(1,1)$}
    \label{fig_1} 
\end{figure}
%%%
\begin{figure}[htbp] 
%%%%%%%%%%%%%%%%%%%%%%%%%%%%%
%\small
%%%%%%%%%%%%%%%%%%%%%%%%%%%%%
% diagram ↓
%%%%%%%%%%%%%%%%%%%%%%%%%%%%%
\begin{tikzpicture}
%%%
\node[outer sep=2pt] (-1) at (-0.2,-0.0) {$\cdots$ \ \ };
\node[outer sep=2pt] (4) at (0.8,1.35) {$r^\prime_3\Lambda(i-4)$};
\node[outer sep=2pt] (6) at (1.2,-1.35) {$[2,4]_{c^\prime}\Lambda(i-2)$};
\node[outer sep=2pt] (7) at (2.4,0) {$[2,3]_{c^\prime}\Lambda(i-2)$};
\node[outer sep=2pt] (8) at (3.8,1.35) {$r^\prime_2\Lambda(i-2)$};
\node[outer sep=2pt] (10) at (3.9,-1.35) {$[1,3]_{c^\prime}\Lambda(i)$};
\node[outer sep=2pt] (11) at (5.2,0) {$[1,2]_{c^\prime}\Lambda(i)$};
\node[outer sep=2pt] (12) at (6.5,1.35) {$r^\prime_1\Lambda(i)$};
\node[outer sep=2pt] (14) at (6.7,-1.35) {$[0,2]_{c^\prime}\Lambda(i+2)$};
\node[outer sep=2pt] (15) at (8.1,0) {$[0,1]_{c^\prime}\Lambda(i+2)$};
\node[outer sep=2pt] (16) at (9.3,1.35) {$r^\prime_0\Lambda(i+2)$};
\node[outer sep=2pt] (22) at (9.9,-1.35) {$[-1,1]_{c^\prime}\Lambda(i+4)$};
\node[outer sep=2pt] (23) at (10.7,0) { \ \ $\cdots$};
%%%
\draw[->] (4) -- (7);
\draw[->] (6) -- (7);
\draw[->] (7) -- (8);
\draw[->] (7) -- (10);
\draw[->] (8) -- (11);
\draw[->] (10) -- (11);
\draw[->] (11) -- (12);
\draw[->] (11) -- (14);
\draw[->] (12) -- (15);
\draw[->] (14) -- (15);
\draw[->] (15) -- (16);
\draw[->] (15) -- (22);
% $\tau$-orbit
\draw[dashed] (4) -- (8) -- (12) -- (16);
\draw[dashed] (-1)  -- (7) -- (11) -- (15) -- (23);
\draw[dashed]  (6) -- (10) -- (14) --  (22);
\end{tikzpicture}
%%%%%%%%%%%%%%%%%%%%%%%%%%%%%
% diagram ↑
%%%%%%%%%%%%%%%%%%%%%%%%%%%%%
%\normalsize
%%%%%%%%%%%%%%%%%%%%%%%%%%%%%
    \caption{The Auslander-Reiten quiver of $\thick T_{c^\prime}(i)$, where the valuation of each arrow is $(1,1)$}
    \label{fig_2} 
\end{figure}
%%%%%%%%%%%%%%%%%%%%%%%%%%%%%%%%%%%%
%           Statement ↑
%%%%%%%%%%%%%%%%%%%%%%%%%%%%%%%%%%%%

We end this section with some observations concerning gradings on $\L$.
By a positive grading on $\L$, we mean a grading that makes $\L$ into a positively graded algebra. 
There are many positive gradings on $\L$ other than assigning degree one to each arrow.
Any positive grading on $\L = KQ/I$ is induced by a grading on $KQ$. 
No matter how $KQ$ is graded, we denote by $\degQ x$ the degree of a homogeneous element $x$ in $KQ$. 
For the positive grading on $\L$ defined by setting each arrow to be degree one, we have $\degQ p = l(p)$ for $p \in \B$.
For an arbitrary positive grading on $\L$ such that $\degQ c > 0$ for all $c \in \C(\L)$, one can get results analogous to results in this section. 
We state the analogous results in the following theorem without proof.

%%%%%%%%%%%%%%%%%%%%%%%%%%%%%%%%%%%%
%           Statement ↓
%%%%%%%%%%%%%%%%%%%%%%%%%%%%%%%%%%%%
\begin{theorem} \label{u_claim_59}
Let $\L$ be a monomial algebra endowed with a positively graded algebra structure such that $\degQ c > 0$ for all $c \in \C(\L)$.
Then the following statements hold.
\begin{enumerate}
    \item The object $T^\prime$ defined as follows is a tilting object of $\sGprojZ\L$:
\begin{align} 
    T^\prime := \bigoplus_{c \in \C(\L)}\,\bigoplus_{0 \leq i < \degQ c} T_c(i),
\end{align}
where $T_c$ is defined as in (\ref{eq_15}).
    \item There exists a triangle equivalence 
    \begin{align}
    \sGprojZ \L \ 
    \cong \prod_{c \in \C(\L)}  \D^{\rm b}\!\left(\mod K\Ac \right)^{(\degQ c)}.
    \end{align}
    \item For $c \in \C(\L)$, the automorphism  $(\degQ c):\sGprojZ \L\to \sGprojZ \L$ restricts to an automorphism  $(\degQ c):\thick T_c \to \thick T_c$, which satisfies that $\tau_c ^{|c|} \cong (-\degQ c)$ as autoequivalences of $\thick T_c$, where $\tau_c$ is the Auslander-Reiten translation of $\thick T_c$.
\end{enumerate}
\end{theorem}
%%%%%%%%%%%%%%%%%%%%%%%%%%%%%%%%%%%%
%           Statement ↑
%%%%%%%%%%%%%%%%%%%%%%%%%%%%%%%%%%%%

%%%%%%%%%%%%%%%%%%%%%%%%%%%%%%%%%%%%
%           Section ↑
%%%%%%%%%%%%%%%%%%%%%%%%%%%%%%%%%%%%

%%%%%%%%%%%%%%%%%%%%%%%%%%%%%%%%%%%%
%           Section ↓
%%%%%%%%%%%%%%%%%%%%%%%%%%%%%%%%%%%%
\section{Realizing stable categories of Gorenstein-projective modules as orbit categories} \label{section_4}
%%%%%%%%%%%%%%%%%%%%%%%%%%%%%%%%%%%%

The aim of this section is to realize $\sGproj \L$ as the stable module category of a self-injective Nakayama algebra.
Moreover, we examine $\sGproj \L$  in the case where $\L$ belongs to one of the following classes of algebras: monomial algebras without overlaps, Nakayama algebras, and $1$-Iwanaga-Gorenstein monomial algebras.

We know from Lemmas \ref{claim_38}  and \ref{u_claim_48} that the induced functor $\Tilde{F}_G: \sGprojZ\L \to \sGproj\L$ is a $G$-covering, where $G$ is the cyclic group generated by the automorphism $(1): \sGprojZ\L \to \sGprojZ\L$. 
Consequently, by Section \ref{preliminaries_2}, the orbit category $\sGprojZ\L/(1)$ becomes a triangulated category and is triangle equivalent to $\sGproj\L$. 
Moreover, the canonical functor $P: \sGprojZ\L \to \sGprojZ\L/(1)$ is a triangulated functor.

It is worth pointing out that for any perfect path $p$,  the canonical image 
\begin{align} 
    [i+1, i+m]\L
    \,\to\,  
    [i+1, i+m-1]\L\oplus [i,i+m]\L
    \,\to\,
    p\L 
    \,\to\, 
    [i+m-m_c, i]\L
\end{align}
in $\sGproj\L$ of the Auslander-Reiten triangle (\ref{eq_5}) in $\sGprojZ\L$ is still Auslander-Reiten.

%%%%%%%%%%%%%%%%%%%%%%%%%%%%%%%%%%%%
%           Subsection ↓
%%%%%%%%%%%%%%%%%%%%%%%%%%%%%%%%%%%%
\subsection{Main result} \label{section_4_1}
%%%%%%%%%%%%%%%%%%%%%%%%%%%%%%%%%%%%

This subsection investigates $\sGproj \L$ via the orbit category induced by $\sGprojZ\L$ and the automorphism $(1)$.

For a class $\mathcal{X}$ of objects of $\sGprojZ\L$, we denote by $P(\mathcal{X})$ the full subcategory of $\sGprojZ\L/(1)$ given by $P(\mathcal{X}) = \{ P(X) \mid X \in  \mathcal{X} \}$.
Then, for any $c \in \C(\L)$ and $i \in \mathbb{Z}$, 
we have
\begin{align}
     P(\thick T_c(i)) = P((\thick T_c)(i)) = (P(\thick T_c))(i) = P(\thick T_c),
\end{align}
where the third equality is obtained from the fact that the induced automorphism $(1) : \sGprojZ\L/(1) \to \sGprojZ\L/(1)$ is isomorphic to the identity functor $\id_{\sGprojZ\L/(1)}$.
Therefore,  (\ref{eq_6}) yields the following description of $\sGprojZ \L/(1)$ as additive categories:
\begin{align} \label{eq_9}
    \sGprojZ \L/(1) = \prod_{c \in \C(\L)} P(\thick T_c).
\end{align}
As the following lemma shows, this is a decomposition into triangulated subcategories.

%%%%%%%%%%%%%%%%%%%%%%%%%%%%%%%%%%%%
%           Statement ↓
%%%%%%%%%%%%%%%%%%%%%%%%%%%%%%%%%%%%
\begin{lemma} \label{u_claim_41}
Under the above situation, we have $P(\thick T_c) = \thick P(T_c)$.
\end{lemma}
\begin{proof}
The fact that the canonical functor $P: \sGprojZ\L \to \sGprojZ\L/(1)$ is triangulated implies that $P(\thick T_c)$ $\subseteq \thick P(T_c)$.
To complete the proof, it suffices to show that $P(\thick T_c)$ is a thick subcategory of $\sGprojZ\L/(1)$ that contains $P(T_c)$.
It is clear that $P(T_c) \in P(\thick T_c)$.
Let $\Sigma$ and $\overline{\Sigma}$ be the suspension functors of $\sGprojZ\L$ and $\sGprojZ\L/(1)$, respectively.
Then we have
\begin{align}
    \overline{\Sigma} P(\thick T_c) = P\Sigma(\thick T_c) = P(\thick T_c).
\end{align}

Next, we claim that $X \oplus Y \in P(\thick T_c)$ implies $X \in P(\thick T_c)$.
Let $X$ and $Y \in \sGprojZ\L/(1)$ satisfy that $X \oplus Y \in P(\thick T_c)$. 
We may assume that $X$ is non-zero in $P(\thick T_c)$.
Then, since $\Hom_{\sGprojZ\L/(1)}(X, X\oplus Y) \not = 0$, the direct summand $X$ must belong to $P(\thick T_c)$ because of (\ref{eq_9}).

Finally, let $X \to Y \to Z \xrightarrow{h} \overline{\Sigma}X$ be a triangle of $\sGprojZ\L/(1)$ with $X, Y \in P(\thick T_c)$.
We may assume that $h \not=0$. 
But then the fact that $\overline{\Sigma}X \in P(\thick T_c)$ implies that $Z \in P(\thick T_c)$ because of (\ref{eq_9}) again.
Therefore, we obtain the desired equality.
\end{proof}
%%%%%%%%%%%%%%%%%%%%%%%%%%%%%%%%%%%%
%           Statement ↑
%%%%%%%%%%%%%%%%%%%%%%%%%%%%%%%%%%%%

For an underlying cycle $c$,  the canonical functor $P: \sGprojZ\L \to \sGprojZ\L/(1)$ restricts to a functor $P_c : \thick T_c \to 
P(\thick T_c)$, which is dense.
Also, by Lemma \ref{u_claim_39}, we have the induced automorphism $(l(c)) : \thick T_c \to \thick T_c$.
Let $G_c$ be the cyclic group generated by the automorphism $(l(c))$.
Then it is easily seen that $P_c : \thick T_c \to P(\thick T_c)$ is a $G_c$-invariant functor whose invariance adjuster is induced by that of the $G$-covering $P: \sGprojZ\L \to \sGprojZ\L/(1)$; see Section \ref{preliminaries_2} for the definitions of $G$-invariant functors and invariance adjusters.

%%%%%%%%%%%%%%%%%%%%%%%%%%%%%%%%%%%%
%           Statement ↓
%%%%%%%%%%%%%%%%%%%%%%%%%%%%%%%%%%%%
\begin{lemma} \label{u_claim_42}
Let $c \in \C(\L)$. Then there exists a triangle equivalence 
\begin{align}
    P\!\left(\thick T_c\right)  \cong \D^{\rm b}(\mod K\Ac)/\tau^{|c|},   
\end{align}    
where $\D^{\rm b}(\mod K\Ac)/\tau^{|c|}$ is the triangulated orbit category induced by $\D^{\rm b}(\mod K\Ac)$ and its Auslander-Reiten translation $\tau$ in the sense of \cite{Keller_2005}.
\end{lemma}

\begin{proof}
For any  $X, Y \in \thick T_c$, we have
\begin{align}
    \Hom_{P(\thick T_c)}(P_cX,  P_cY) 
    &=  \Hom_{\sGprojZ\L/(1)}(PX, PY) & & \\[1mm]
    &\cong \bigoplus_{i \in \mathbb{Z}} \sHomZ_{\L}( X, Y(i)) & & \\
    &= \bigoplus_{i \in \mathbb{Z}} \sHomZ_{\L}(X, Y(il(c))) & (\mbox{by } (\ref{eq_6}) \mbox{ and Lemma \ref{u_claim_39}}) &\\
    &= \bigoplus_{i \in \mathbb{Z}} \Hom_{\thick T_c}(X, Y(il(c))).& &
\end{align}
Thus the fact that $P_c$ is dense yields that the induced functor $P_c :  \thick T_c \to P(\thick T_c)$ is a $G_c$-covering. 
Recall from Section \ref{section_3} that $\thick T_c$ is triangle equivalent to $\D^{\rm b}(\mod K \Ac)$  and that the Auslander-Reiten translation $\tau_c:\thick T_c \to \thick T_c$ satisfies $\tau_c^{-|c|} \cong (l(c))$ as autoequivalences.
Therefore, we obtain the following triangle equivalences:
\begin{align}
    P(\thick T_c) \cong \thick T_c/(l(c)) \cong \D^{\rm b}(\mod K\Ac)/\tau^{-|c|} \cong \D^{\rm b}(\mod K\Ac)/\tau^{|c|},
\end{align}
where $\tau$ is the Auslander-Reiten translation for $\D^{\rm b}(\mod K\Ac)$.
Here, we use the fact that Auslander-Reiten translation is unique up to functorial isomorphism.
\end{proof}
%%%%%%%%%%%%%%%%%%%%%%%%%%%%%%%%%%%%
%           Statement ↑
%%%%%%%%%%%%%%%%%%%%%%%%%%%%%%%%%%%%

For $c\in \C(\L)$, we denote by $\L_c$ the connected self-injective Nakayama algebra 
\begin{align} \label{eq_11}
    KQ_c \left/ R_c^{m_c+1}\right.,
\end{align}
where $Q_c$ is the following cyclic quiver with $|c|$ vertices:
\begin{align}
  Q_c : \xymatrix@M=2.5mm{
    1 \ar[r] & 2 \ar[r] & \cdots \ar[r] & |c|, \ar@/^15pt/[lll]
    }\\[-3mm]
\end{align}
and $R_c$ is the arrow ideal of $KQ_c$. 
Recall that $\Ac$ is the linear quiver  
\begin{align}
 1 \rightarrow 2 \rightarrow \cdots \rightarrow m_c.
\end{align}
The following theorem is the main result of this section.

%%%%%%%%%%%%%%%%%%%%%%%%%%%%%%%%%%%%
%           Statement ↓
%%%%%%%%%%%%%%%%%%%%%%%%%%%%%%%%%%%%
\begin{theorem} \label{u_claim_43}
Let $\L$ be a monomial algebra. Then there exist  triangle equivalences
\begin{align}
 \sGproj \L \  \cong \prod_{c \in \C(\L)} \D^{\rm b}(\mod K\Ac)/\tau^{|c|} \ \cong   \prod_{c \in \C(\L)} \smod  \L_c,   
\end{align}
where $\D^{\rm b}(\mod K\Ac)/\tau^{|c|}$ is the triangulated orbit category induced by $\D^{\rm b}(\mod K\Ac)$ and its Auslander-Reiten translation $\tau$ in the sense of \cite{Keller_2005}, and $\L_c$ is defined as in (\ref{eq_11}).
\end{theorem}  
\begin{proof}
It follows from (\ref{eq_9}) and Lemma \ref{u_claim_42} that  there exist triangle equivalences
\begin{align}
    \sGproj \L 
    \cong   \sGprojZ \L/(1) 
    = \prod_{c \in \C(\L)} P(\thick T_c)
    \cong \prod_{c \in \C(\L)} \D^{\rm b}(\mod K\Ac)/\tau^{|c|}.
\end{align}
Now, let us consider the self-injective Nakayama algebra $\L_c$. 
Since $Q_c$ is a cyclic quiver, all underlying cycles are equivalent, which means that $\C(\L_c)$ is a singleton, say $\C(\L_c)=\{c^\prime\}$. 
Then the fact that all arrows of $Q_c$ are perfect in $\L_c$ implies that $|c^\prime| = |c|$.
Moreover, it follows from the definition of $\L_c$ that $m_{c^\prime}^{\L_c} = |\mathbb{P}_{\L_c}(c^\prime)| = m_c^{\L}$.
Thus the above equivalences applied to $\L_c$ yield that  
\begin{align}
    \smod \L_c = \sGproj \L_c  \cong  \D^{\rm b}(\mod K\A_{c^\prime})/\tau^{|c^\prime|} \cong  \D^{\rm b}(\mod K\Ac)/\tau^{|c|}.
\end{align}
This completes the proof.
\end{proof}
%%%%%%%%%%%%%%%%%%%%%%%%%%%%%%%%%%%%
%           Statement ↑
%%%%%%%%%%%%%%%%%%%%%%%%%%%%%%%%%%%%

%%%%%%%%%%%%%%%%%%%%%%%%%%%%%%%%%%%%
%           Statement ↓
%%%%%%%%%%%%%%%%%%%%%%%%%%%%%%%%%%%%
\begin{remark}{\rm \label{remark_3}
As in the graded case, when $\L$ is Iwanaga-Gorenstein, Theorem \ref{u_claim_43} characterizes the singularity category $\D_{\rm sg}(\mod\L)$ of $\L$ and extends \cite[Theorem 6.3.6]{Lu-Zhu_2021}, where the authors describe $\D_{\rm sg}(\mod\L)$ explicitly for $1$-Iwanaga-Gorenstein monomial algebras $\L$.  
In Section \ref{section_4_4} below, we recover this result from Theorem \ref{u_claim_43}.
}\end{remark}
%%%%%%%%%%%%%%%%%%%%%%%%%%%%%%%%%%%%
%           Statement ↑
%%%%%%%%%%%%%%%%%%%%%%%%%%%%%%%%%%%%

%%%%%%%%%%%%%%%%%%%%%%%%%%%%%%%%%%%%
%           Statement ↓
%%%%%%%%%%%%%%%%%%%%%%%%%%%%%%%%%%%%
\begin{example}{\rm \label{ex_5}
Let $\L$ be as in Example \ref{ex_3}.
Recall that $\C(\L)=\{c=a_{123}, c^\prime=a_{45}\}$ with $m_c= 4$ and $m_{c^\prime} = 3$ and $\mathbb{E}_{\L}^{\rm co}=\{r_1 = a_{12}, r_2 =a_{3}, r^\prime_1 =a_{45}\}$.
In particular,  $|c| =  2$ and $|c^\prime| = 1$. 
Then Theorem \ref{u_claim_43} yields that
\begin{align}
    \sGproj \L\, 
    \cong\,  \D^{\rm b}(\mod K \Ac)/\tau^{2} \times \D^{\rm b}(\mod K \A_{c^\prime})/\tau,
\end{align}
where $\Ac : 1 \to 2 \to 3 \to 4$ and $\A_{c^\prime} : 1 \to 2 \to 3$.
We note that 
\begin{align}
    \D^{\rm b}(\mod K\Ac)/\tau^{2} \cong \smod \L_{c} \quad  \mbox{and} \quad  \D^{\rm b}(\mod K\A_{c^\prime})/\tau \cong \smod \L_{c^\prime},
\end{align}
where $\L_{c} = KQ_{c} / R_c^{5}$ and  $\L_{c^\prime} = KQ_{c^\prime} / R_{c^\prime}^{4}$ with $Q_c: \xymatrix{1  \ar@<0.6ex>[r] & 2 \ar@<0.6ex>[l]}$ and $Q_{c^\prime}:  \xymatrix{1 \ar@(ur,dr)}$\quad\ \ .
On the other hand, the Auslander-Reiten quiver of $\sGproj\L$ is given as follows:
%%%%%%%%%%%%%%%%%%%%%%%%%%%%%
% diagram ↓
%%%%%%%%%%%%%%%%%%%%%%%%%%%%%
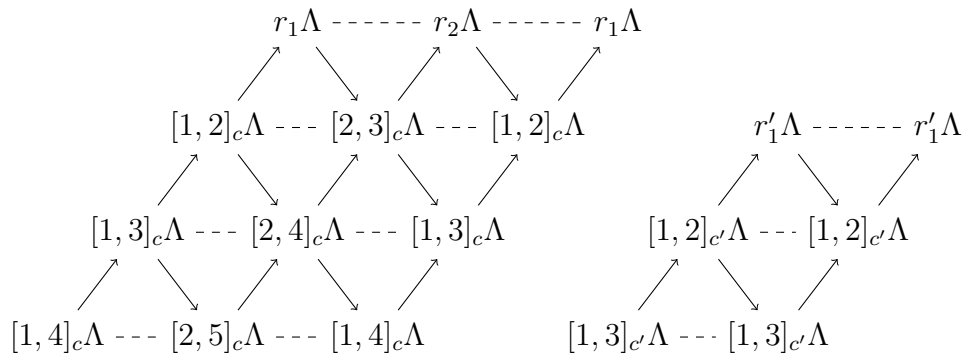
\[
\begin{tikzpicture}[scale=.99]
\node (1) at (-2.14,-2.8) {$[1,4]_c\L$};
\node (2) at (-1.07,-1.4) {$[1,3]_c\L$};
\node (3) at (0,0) {$[1,2]_c\L$};
\node (4) at (1.07,1.4) {$r_1\L$};
\node (5) at (0,-2.8) {$[2,5]_c\L$};
\node (6) at (1.07,-1.4) {$[2,4]_c\L$};
\node (7) at (2.14,0) {$[2,3]_c\L$};
\node (8) at (3.21,1.4) {$r_2\L$};
\node (9) at (2.14,-2.8) {$[1,4]_c\L$};
\node (10) at (3.21,-1.4) {$[1,3]_c\L$};
\node (11) at (4.28,0) {$[1,2]_c\L$};
\node (12) at (5.35,1.4) {$r_1\L$};
\node (13) at (5.35,-2.8) {$[1,3]_{c^\prime}\L$};
\node (14) at (6.42,-1.4) {$[1,2]_{c^\prime}\L$};
\node (15) at (7.49,0) {$r^\prime_1\L$};
\node (16) at (7.49,-2.8) {$[1,3]_{c^\prime}\L$};
\node (17) at (8.56,-1.4) {$[1,2]_{c^\prime}\L$};
\node (18) at (9.63,0) {$r^\prime_1\L$};
\draw[->] (1) -- (2);
\draw[->] (2) -- (3);
\draw[->] (2) -- (5);
\draw[->] (3) -- (4);
\draw[->] (3) -- (6);
\draw[->] (4) -- (7);
\draw[->] (5) -- (6);
\draw[->] (6) -- (7);
\draw[->] (6) -- (9);
\draw[->] (7) -- (8);
\draw[->] (7) -- (10);
\draw[->] (8) -- (11);
\draw[->] (9) -- (10);
\draw[->] (10) -- (11);
\draw[->] (11) -- (12);
\draw[->] (13) -- (14);
\draw[->] (14) -- (15);
\draw[->] (14) -- (16);
\draw[->] (15) -- (17);
\draw[->] (16) -- (17);
\draw[->] (17) -- (18);
% $\tau$-orbit
\draw[dashed] (1) -- (5) -- (9);
\draw[dashed] (2) -- (6) -- (10);
\draw[dashed] (3) -- (7) -- (11);
\draw[dashed] (4) -- (8) -- (12);
\draw[dashed] (13) -- (16);
\draw[dashed] (14) -- (17);
\draw[dashed] (15) -- (18);
\end{tikzpicture}
\]
%%%%%%%%%%%%%%%%%%%%%%%%%%%%%
% diagram ↑
%%%%%%%%%%%%%%%%%%%%%%%%%%%%%
where the valuation of each arrow is $(1,1)$, and in each component, the extreme left and extreme right full subquivers with the same vertices are identified.
}\end{example}
%%%%%%%%%%%%%%%%%%%%%%%%%%%%%%%%%%%%
%           Statement ↑
%%%%%%%%%%%%%%%%%%%%%%%%%%%%%%%%%%%%

In the next three subsections, we apply Theorem \ref{u_claim_43} to certain three classes of monomial algebras.

%%%%%%%%%%%%%%%%%%%%%%%%%%%%%%%%%%%%
%           Subsection ↓
%%%%%%%%%%%%%%%%%%%%%%%%%%%%%%%%%%%%
\subsection{The case of monomial algebras without overlaps} \label{section_4_2}
%%%%%%%%%%%%%%%%%%%%%%%%%%%%%%%%%%%%

Recall that there exists no overlap in $\L$ if there exists no overlap between any perfect paths. 
As mentioned at the end of Section \ref{preliminaries}, Chen, Shen and Zhou \cite[Proposition 5.9]{Chen-Shen-Zhou_2018} provided conditions under which there is no overlap in $\L$ in terms of the stable category $\sGproj \L$.  
We now use the Hasse quivers $H(\mathbb{P}_{\L}, \preceq)$ and $H(\mathbb{P}_{\L}, \leq)$ to give further equivalent conditions.

%%%%%%%%%%%%%%%%%%%%%%%%%%%%%%%%%%%%
%           Statement ↓
%%%%%%%%%%%%%%%%%%%%%%%%%%%%%%%%%%%%
\begin{proposition} \label{u_claim_45}
The following conditions are equivalent.
\begin{enumerate}
    \item There exists no overlap in $\L$.
    \item Every perfect path is both elementary and co-elementary.
    \item Every perfect path is both a sink and a source in $H(\mathbb{P}, \preceq)$.
    \item Every perfect path is both a sink and a source in  $H(\mathbb{P}, \leq)$.
    
\end{enumerate}
In this case, 
\begin{align}
    \sGprojZ \L \cong  \D^{\rm b}\!\left(\mod K \right)^{(l_\L)} \ 
      \mbox{and} \ \  
    \sGproj \L \cong \prod_{c \in \C(\L)}\D^{\rm b}(\mod K)/\tau^{|c|} \cong \prod_{c \in \C(\L)} \smod  \L_c 
\end{align}
as triangulated categories, where $l_\L:=\sum_{c \in \C(\L)} l(c)$, and each $\L_c$ is radical square zero. 
\end{proposition}  
\begin{proof} 
The equivalence of (2) through (4) is a consequence of Proposition \ref{u_claim_16}.
We now show the equivalence of (1) and (2).
The implication (2) $\Rightarrow$ (1) immediately follows from Lemma \ref{u_claim_28}. 
Conversely, assume that there is no overlap in $\L$. 
Then any perfect path $p$ contains no perfect left divisor and no perfect right divisor, which means that $p \in \mathbb{E}_\L \cap \mathbb{E}^{\rm co}_\L$.
For the last statement, together with the fact that $m_c= 1$ for any $c \in \C(\L)$, Theorems \ref{u_claim_34} and \ref{u_claim_43} yield the desired triangle equivalences.
\end{proof}
%%%%%%%%%%%%%%%%%%%%%%%%%%%%%%%%%%%%
%           Statement ↑
%%%%%%%%%%%%%%%%%%%%%%%%%%%%%%%%%%%%

%%%%%%%%%%%%%%%%%%%%%%%%%%%%%%%%%%%%
%           Subsection ↓
%%%%%%%%%%%%%%%%%%%%%%%%%%%%%%%%%%%%
\subsection{The case of Nakayama algebras} \label{section_4_3}
%%%%%%%%%%%%%%%%%%%%%%%%%%%%%%%%%%%%

Next, suppose that $\L=KQ/I$ is a connected Nakayama algebra.
It is known that $Q$ is either a linear quiver or a cyclic quiver.
In the former case, the algebra $\L$ is always CM-free.
On the other hand, in case $\L$ is not CM-free, the fact that $Q$ is a cyclic quiver implies that $\C(\L)$ is a singleton, say $\C(\L) = \{c\}$. 
Moreover, we have $l(c)=|Q_0|$, the number of the vertices of $Q$.
Hence we obtain the following consequence of  Theorems \ref{u_claim_34} and \ref{u_claim_43}; compare  \cite[Proposition 1]{Ringel_2013}.

%%%%%%%%%%%%%%%%%%%%%%%%%%%%%%%%%%%%
%           Statement ↓
%%%%%%%%%%%%%%%%%%%%%%%%%%%%%%%%%%%%
\begin{proposition} \label{u_claim_44}
Let $\L=KQ/I$ be a connected Nakayama algebra that is not CM-free.
Then $\C(\L)$ is a singleton, say $\C(\L) = \{c\}$, and we  have
\begin{align}
\sGprojZ \L \ \cong \D^{\rm b}\left(\mod K\Ac \right)^{(|Q_0|)} \  \mbox{and} \quad  \sGproj \L \cong \D^{\rm b}(\mod K\Ac)/\tau^{|c|} \cong \smod  \L_c    
\end{align}
as triangulated categories.
\end{proposition} 
%%%%%%%%%%%%%%%%%%%%%%%%%%%%%%%%%%%%
%           Statement ↑
%%%%%%%%%%%%%%%%%%%%%%%%%%%%%%%%%%%%

This proposition extends \cite[Corollary 4.3.6]{Lu-Zhu_2021} to arbitrary Nakayama algebras.

%%%%%%%%%%%%%%%%%%%%%%%%%%%%%%%%%%%%
%           Subsection ↓
%%%%%%%%%%%%%%%%%%%%%%%%%%%%%%%%%%%%
\subsection{The case of 1-Iwanaga-Gorenstein monomial algebras} \label{section_4_4}
%%%%%%%%%%%%%%%%%%%%%%%%%%%%%%%%%%%%

First of all, we characterize underlying cycles whose subarrows are all perfect.

%%%%%%%%%%%%%%%%%%%%%%%%%%%%%%%%%%%%
%           Statement ↓
%%%%%%%%%%%%%%%%%%%%%%%%%%%%%%%%%%%%
\begin{proposition} \label{u_claim_54}
The following conditions are equivalent for any $c=a_1 \cdots a_n \in \C(\L)$, where $a_i$ is an arrow for all $i$.
\begin{enumerate}
    \item Each $a_i$ is perfect.
    \item $|c| = l(c)$.
    \item The following two conditions are satisfied:
    \begin{enumerate}
        \item Any non-zero non-trivial subpath $p$ of an $m$-th power $c^m$ of $c$ with the property that $pq \in \F$ for some $q \in \B$ is perfect.
        \item  $a_i \cdots a_{i+m_c} \in \F$ for all $1 \leq i \leq n$, where the index $j$ of $a_j$ is considered modulo $n$.
    \end{enumerate}
\end{enumerate}
\end{proposition}
\begin{proof}  
The equivalence of (1) and (2) and the implication from (3) to (1) are both clear.
To prove that (1) implies (3), assume that $a_i \in \mathbb{P}_\L$ for all $i$.
Since each $a_i$ is co-elementary, every non-zero non-trivial subpath $p$ of $c^m$ is perfect by Theorem \ref{u_claim_25}.
This shows (a).
For (b), let us recall that $Y_c = \{ r \in \mathbb{E}_{\L}^{\rm co} \mid c_r = c \mbox{ in } \C(\L) \}$.
Then, since $Y_c =\{ a_i \mid 1 \leq i \leq n\}$, it follows from Lemma \ref{u_claim_51} that $a_i \cdots a_{i+m_c} \in \F$ for all $1 \leq i \leq n$.
\end{proof}
%%%%%%%%%%%%%%%%%%%%%%%%%%%%%%%%%%%%
%           Statement ↑
%%%%%%%%%%%%%%%%%%%%%%%%%%%%%%%%%%%%

We see from \cite[Proposition 5.1]{Chen-Shen-Zhou_2018} that for a quadratic monomial algebra, any underlying cycles satisfy the equivalent conditions in Proposition \ref{u_claim_54}.
The following lemma shows that this is also the case for $1$-Iwanaga-Gorenstein monomial algebras.

%%%%%%%%%%%%%%%%%%%%%%%%%%%%%%%%%%%%
%           Statement ↓
%%%%%%%%%%%%%%%%%%%%%%%%%%%%%%%%%%%%
\begin{proposition} \label{u_claim_55} 
Let $\L$ be a $1$-Iwanaga-Gorenstein monomial algebra. 
Then any underlying cycle satisfies the equivalent conditions in Proposition \ref{u_claim_54}.
\end{proposition}
\begin{proof} 
Let $c= a_1 \cdots a_n \in \C(\L)$ with  $a_i \in Q_1$.
We claim that each $a_i$ is perfect. 
Fix an index $i$ with $1 \leq i \leq n$. 
Since $\L$ is $1$-Iwanaga-Gorenstein, the cyclic module $a_i \L$ is either projective or non-projective Gorenstein-projective. 
Assume that $a_i\L$ is projective. 
Then $R(a_i)$ is empty (because for any non-zero path $p$, we have the following short exact sequence in $\mod\L$:
\begin{align}
    0 \to \bigoplus_{q \in R(p)} q \L \xrightarrow{\iota} e_{t(p)}\L \xrightarrow{\pi} p\L \to 0,
\end{align}
where $\iota$ is the canonical inclusion, and $\pi$ is the left multiplication by $p$, which is a projective cover of $p\L$). 
Thus we can conclude that $a_i$ cannot be a left divisor of any path in $\F$. 
On the other hand, since $a_i$ is a subarrow of $c$ with $R(a_i) = \varnothing$, there exists a co-elementary subpath $r$ of $c$ such that $r = r^\prime a_i r^{\prime\prime}$ for some $r^{\prime}\in \B_{>0}$ and  $r^{\prime\prime} \in \B$.
But then, \cite[Theorem 5.1.4]{Lu-Zhu_2021} implies that the left divisor $r^\prime$ of the co-elementary path $r$ is a perfect path, a contradiction.
Therefore, $a_i\L$ is non-projective Gorenstein-projective.
Now, it follows from \cite[Proposition 4.6]{Chen-Shen-Zhou_2018} that $R(a_i) =\{p\}$ for some $p\in \mathbb{P}_\L$, which implies that $a_i p \in \F$, so that $a_i$ is perfect by \cite[Theorem 5.1.4]{Lu-Zhu_2021} again. 
\end{proof}
%%%%%%%%%%%%%%%%%%%%%%%%%%%%%%%%%%%%
%           Statement ↑
%%%%%%%%%%%%%%%%%%%%%%%%%%%%%%%%%%%%

Let $\L=KQ/I$ be a $1$-Iwanaga-Gorenstein monomial algebra. 
Recall from \cite{Lu-Zhu_2021} that a {\it repetition-free cycle} is a cycle $c=a_1 \cdots a_n$ with $a_i \in Q_1$ such that there exists an integer $r \geq 2$ such that $a_i \cdots a_{i+r-1}$ belongs to $\F$ for all $1 \leq i \leq n$, where we consider the index $j$ of $a_j$ modulo $n$.
We refer to the integer $r$, uniquely determined,  as the {\it length of relations for $c$}.
Combining the above two propositions and \cite[Theorem 5.1.4]{Lu-Zhu_2021},  we obtain that underlying cycles $c$ in  $\L$ are precisely repetition-free cycles and satisfy that  $|c| = l(c)$, and the length of relations for $c$ equals $m_c +1$.
Therefore, Theorem \ref{u_claim_43} specialized to $\L$ is nothing but \cite[Theorem 6.3.6]{Lu-Zhu_2021}.

%%%%%%%%%%%%%%%%%%%%%%%%%%%%%%%%%%%%
%           Section ↑
%%%%%%%%%%%%%%%%%%%%%%%%%%%%%%%%%%%%

%%%%%%%%%%%%%%%%%%%%%%%%%%%%%%%%%%%%
%          ↓ Acknowledgments
%%%%%%%%%%%%%%%%%%%%%%%%%%%%%%%%%%%%
\section*{Acknowledgments} 
The authors deeply thank the referee for their careful reading and for providing many valuable comments and suggestions on the mathematics, which have improved the readability of the paper.
%%%%%%%%%%%%%%%%%%%%%%%%%%%%%%%%%%%%
%         ↑ Acknowledgments
%%%%%%%%%%%%%%%%%%%%%%%%%%%%%%%%%%%%

%\bibliographystyle{alpha}
\bibliographystyle{plain}
\bibliography{ref}

\end{document}